\documentclass[preprint,12pt,sort&compress]{elsarticle}

\usepackage{geometry}
\usepackage{lineno}
\modulolinenumbers[5]

\usepackage[utf8x]{inputenc}
\usepackage[english]{babel}
\usepackage[T1]{fontenc}

\usepackage{amsmath}
\usepackage{amssymb}
\usepackage{graphicx}
\usepackage{subfig}

\usepackage[pdfborder={0 0 0},colorlinks,allcolors=blue]{hyperref}
\usepackage{xcolor}

\usepackage{pgfplots}
\usetikzlibrary{arrows}
\usetikzlibrary{calc}
\pgfplotsset{compat=1.16}

\usepackage{amsthm}
\theoremstyle{definition}

\newtheorem{example}{Example}[section]
\newtheorem{remark}{Remark}[section]

\usepackage{tcolorbox}

\usepackage{siunitx} 

\makeatletter
\newcommand*\bigcdot{\mathpalette\bigcdot@{.8}}
\newcommand*\bigcdot@[2]{\mathbin{\vcenter{\hbox{\scalebox{#2}{$\m@th#1\bullet$}}}}}
\makeatother

\usepackage{bm} 

\usepackage{cleveref}
\hypersetup{colorlinks=true,allcolors=blue,pdfauthor=author}
\crefname{equation}{Equation}{Equations}
\crefname{figure}{Figure}{Figures}
\crefname{section}{Section}{Sections}
\crefname{remark}{Remark}{Remarks}
\crefname{appendix}{}{}

\usepackage{xparse}

\NewDocumentCommand \an{ m }{%
    \langle {#1} \rangle%
}

\newcommand*{\medcap}{\mathbin{\scalebox{1.5}{\ensuremath{\cap}}}}%

\newcommand\notype[1]{\unskip}

\title{Coupling of IGA and Peridynamics for Air-Blast Fluid-Structure Interaction Using an Immersed Approach}

\begin{document}

\begin{frontmatter}



\author[brown]{Masoud Behzadinasab\corref{mycorrespondingauthor}}
\ead{masoud\_behzadinasab@brown.edu}
\author[suny,iacs]{Georgios Moutsanidis}
\author[sandia]{Nathaniel Trask}
\author[ut]{John T. Foster}
\author[brown]{Yuri Bazilevs\corref{mycorrespondingauthor}}
\ead{yuri\_bazilevs@brown.edu}

\address[brown]{School of Engineering, Brown University, 184 Hope St., Providence, RI 02912, USA}
\address[suny]{Department of Civil Engineering, Stony Brook University, Stony Brook, NY 11794, USA}
\address[iacs]{Institute for Advanced Computational Science, Stony Brook University, Stony Brook, NY 11794, USA}
\address[sandia]{Center for Computing Research, Sandia National Laboratories, Albuquerque, NM 87185, USA\textsuperscript{1}\footnote{\textsuperscript{1}Sandia National Laboratories is a multi-mission laboratory managed and operated by National Technology and Engineering Solutions of Sandia, LLC., a wholly owned subsidiary of Honeywell International, Inc., for the U.S. Department of Energy’s National Nuclear Security Administration under contract DE-NA0003525. This paper describes objective technical results and analysis. Any subjective views or opinions that might be expressed in the paper do not necessarily represent the views of the U.S. Department of Energy or the United States Government.}}
\address[ut]{Hildebrand Department of Petroleum \& Geosystems Engineering, The University of Texas at Austin, Austin, TX 78712, USA}

\cortext[mycorrespondingauthor]{Corresponding authors.}

\begin{abstract}
We present a novel formulation based on an immersed coupling of Isogeometric Analysis (IGA) and Peridynamics (PD) for the simulation of fluid–structure interaction (FSI) phenomena for air blast. We aim to develop a practical computational framework that is capable of capturing the mechanics of air blast coupled to solids and structures that undergo large, inelastic deformations with extreme damage and fragmentation. An immersed technique is used, which involves an {\em a priori} monolithic FSI formulation with the implicit detection of the fluid-structure interface and without limitations on the solid domain motion. The coupled weak forms of the fluid and structural mechanics equations are solved on the background mesh. Correspondence-based PD is used to model the meshfree solid in the foreground domain. We employ the Non-Uniform Rational B-Splines (NURBS) IGA functions in the background and the Reproducing Kernel Particle Method (RKPM) functions for the PD solid in the foreground. We feel that the combination of these numerical tools is particularly attractive for the problem class of interest due to the higher-order accuracy and smoothness of IGA and RKPM, the benefits of using immersed methodology in handling the fluid-structure coupling, and the capabilities of PD in simulating fracture and fragmentation scenarios. Numerical examples are provided to illustrate the performance of the proposed air-blast FSI framework. 
\end{abstract}


\end{frontmatter}


\section{Introduction}
\label{sec:intro}

In~\cite{bazilevs2017new1,bazilevs2017new2}, the authors developed an immersed fluid--structure interaction (FSI) formulation for air blast. The formulation made use of the weak forms of the fluid and structural mechanics equations and two meshes, background and foreground.

The coupled solution was approximated using the basis functions of a fixed background mesh for both the fluid and solid parts of the coupled FSI problem resulting in an {\em a priori} strongly coupled formulation and elimination of any mesh distortion issues. It was demonstrated that the use of high-order accurate and smooth background dicsretizations, such as those in spline-based Isogeometric Analysis (IGA)~\cite{hughes2005isogeometric,cottrell2009isogeometric}, resulted in much higher solution quality, especially for the solid, due to the higher-order accuracy of the strain-rate approximation and the circumvention of cell-crossing instabilities occurring in the traditional Material-Point Methods (MPMs)~\cite{steffen2008analysis}. This finding has led to the development of the Isogeometric MPM in~\cite{moutsanidis2020iga}.

The foreground mesh in the originally developed air-blast FSI framework was mainly used to carry out the solid-domain quadrature, track the solid current configuration, and store the solid constitutive law history variables. In the follow-up work of~\cite{moutsanidis2018hyperbolic}, the foreground discretization comprised of the Reproducing Kernel Particle Method (RKPM) functions~\cite{liu1995reproducing,chen1996reproducing} was used to directly solve, in the weak form, the time-dependent partial differential equation (PDE) governing the phase-field variable of a brittle fracture model.

While the resulting coupled FSI formulation presents a very promising approach, especially for applications involving blast, fracture, and fragmentation, the framework could benefit from further improvements outlined in what follows.

At the engineering scale, blast on structures involves sophisticated structural dicsretizations that often make use of several structural element types ranging from 3D continuum solids to discrete point masses. These are implemented in sophisticated structural solvers that took decades to develop and that have a large user base. As a result, a more {\em modular} approach is desired where the foreground discretization is directly used to compute the nodal forces on the foreground mesh, which are subsequently transferred to the background mesh in order the maintain the the advantages of the {\em a priori} strong coupling, as per the original formulation. In the present paper we address this challenge, and, as one instantiation of the proposed approach, demonstrate how to use a state-based Peridynamic (PD) formulation of a solid~\cite{silling2007peridynamic} into the FSI framework.

The PD theory~\cite{silling2000reformulation,silling2007peridynamic} has been in development for two decades. The PD integro-differential governing equations of solid mechanics naturally accommodate nonlocal physics as well as low-regularity solutions commonly encountered in fracture mechanics problems. PD has been extensively applied to simulate crack growth in a wide range of materials and structures under different loading conditions~\cite{gerstle2007peridynamic,foster2009dynamic,foster2010viscoplasticity,agwai2011predicting,weckner2013viscoelastic,bobaru2015cracks,yaghoobi2017fracture,jafarzadeh2018peridynamic,behzadinasab2018peridynamics,kramer2019third,chen2019peridynamic,behzadinasab2019third,gao2020peridynamics,butt2021peridynamic}. A limited number of PD formulations for FSI phenomena may be found in the literature, e.g., hydraulic fracturing~\cite{ouchi2015fully,oterkus2017fully,ni2020hybrid}, ice-water interaction~\cite{liu2020modeling}, ice-structure interaction~\cite{vazic2019peridynamic,lu2020peridynamic}, and fluid-elastic structure interaction~\cite{gao2020fluid,sun2020smoothed}. In our proposed FSI methodology, a stable correspondence-based PD formulation~\cite{behzadinasab2020semi} is utilized, which enables direct use of classical continuum material and damage models. A great deal of numerical complexity is eliminated as there is no need to explicitly track crack interface (advantage of using PD) or fluid-structure interfaces (advantage of using immersed techniques). The developed computational methodology is applied to the simulation of fracture in brittle and ductile solids in an air-blast environment.

While constraining the solid to the background kinematics gives the desired strong coupling, modeling of fragmentation becomes challenging since the smooth background discretion is not designed to excel in approximating discontinuities in the solution. In practice, when continuum-damage or phase-field approaches are employed to simulate fracture, the damage zones, whose size scales with that of the background elements, appear to be artificially thick leading to excessive solid damage. The size of the damage zones may be reduced by refining the background mesh, however, this leads to a significant increase in the computational costs. While we do not address this issue in the present work in a rigorous way, and leave it for the future developments leveraging such techniques as developed in~\cite{wang2021consistent}, we demonstrate a way to use the foreground discretization to model fracture and fragmentation due to blast in a more effective manner.

The paper is outlined as follows. In \cref{sec:formulation}, the continuum- and discrete-level FSI formulations, which include the fluid and solid subproblems in the weak form and their coupling, are presented. In \cref{sec:implementation}, numerical examples are presented to verify the accuracy of the coupled IGA-PD FSI formulation and to demonstrate the ability of the new formulation to simulate fracture and fragmentation in brittle and ductile solid structures subjected to blast loading. \cref{sec:conclusions} makes concluding remarks and outlines future research directions to address the remaining numerical challenges.

\section{Mathematical formulation}
\label{sec:formulation}

In this section, we develop an immersed FSI formulation that accommodates compressible flow and an inelastic solid or structural model coming from a standalone weak formulation and its discretization. We first provide a general continuum variational framework and then specialize to the coupling of IGA on the background and PD on the foreground grids at the discrete level. The resulting methodology resembles the classical Immersed Boundary Method~\cite{peskin1972flow} and a more recent Immersed Finite Element Method~\cite{zhang2004immersed}, but without the use of {\em ad hoc} smoothed delta functions and with the force transfer motivated by variational arguments. In what follows, all the vectors are column vectors unless otherwise noted. Bold symbols denote tensors of rank one (i.e., vectors) or higher. Underscored terms indicate the PD field variables at the bond level.

\subsection{Weak form of the coupled FSI problem}
\label{sec:weakCoupled}

We take, as the starting point, the coupled FSI formulation proposed in~\cite{bazilevs2017new1}, where the fluid mechanics part of the problem is governed by the Navier--Stokes equations for compressible flow and the structural mechanics part is modeled as an isothermal large-deformation inelastic solid. Let $\Omega$ denote the combined fluid and solid domain, and let $\Omega^f$ and $\Omega^s$ denote the individual fluid and solid subdomains in the spatial configuration, such that $\Omega^f\bigcup\Omega^s = \Omega$ and $\Omega^f\bigcap\Omega^s = \emptyset$. Both the fluid and solid problems are formulated in a variational setting, where we define the following semilinear forms and linear functionals that comprise the weak forms of the fluid and solid subproblems
\begin{align}
\label{MfCont}
M^f_{\omega} (\mathbf{W},\mathbf{Y}) = \int_{\omega} \mathbf{W} \cdot \mathbf{A}^f_0 \mathbf{Y}_{,t}~d\omega, 
\end{align}
\begin{align}
\label{BfCont}
B^f_{\omega} (\mathbf{W},\mathbf{Y}) = \int_{\omega} \mathbf{W} \cdot \mathbf{A}^f_i  \mathbf{Y}_{,i}~d\omega - \int_{\omega} \mathbf{W}_{,i} \cdot (\mathbf{F}^{p}_{i} - \mathbf{F}^{d}_{i})~d\omega,
\end{align}
\begin{align}
\label{FfCont}
F^f_{\omega} (\mathbf{W}) = \int_{\omega} \mathbf{W} \cdot \mathbf{S}^f~d\omega + \int_{\Gamma^f_H} \mathbf{W} \cdot \mathbf{H}^f~d\Gamma,
\end{align}
\begin{align}
\label{MsCont}
M^s_{\omega}(\mathbf{W},\mathbf{Y}) = \int_{\omega} \mathbf{W} \cdot \mathbf{A}^s_0 \mathbf{Y}_{,t}~d\omega,
\end{align}
\begin{align}
\label{BsCont}
B^s_{\omega}(\mathbf{W},\mathbf{Y}) = \int_{\omega} \mathbf{W} \cdot \mathbf{A}^s_i  \mathbf{Y}_{,i}~d\omega + \int_{\omega} \mathbf{W}_{,i} \cdot \mathbf{F}^{\sigma}_{i}~d\omega,
\end{align}
\begin{align}
\label{FsCont}
F^s_{\omega}(\mathbf{W}) = \int_{\omega} \mathbf{W} \cdot \mathbf{S}^s~d\omega + \int_{\Gamma^s_H} \mathbf{W} \cdot \mathbf{H}^s~d\Gamma.
\end{align}
Here, $\mathbf{Y}$ denotes a set of pressure-primitive variables~\cite{hauke1998comparative,hauke2001simple},
\begin{equation}
\mathbf{Y} = 
\begin{bmatrix}
p \\
\mathbf{v} \\
T
\end{bmatrix},
\label{eqn:Y}
\end{equation}
where $p$ is the pressure, $\mathbf{v}$ is the material-particle velocity vector, and $T$ is the temperature. $\mathbf{Y}$ and $\mathbf{W}$, the vector-valued trial and test functions, respectively, are the members of $\mathcal{S}$ and $\mathcal{V}$, the corresponding trial and test function spaces, respectively, \textit{defined on all of} $\Omega$. The superscripts $f$ and $s$ correspond to the fluid and solid quantities, respectively. $\Gamma^f_H$ and $\Gamma^s_H$ are the subsets of the fluid- and solid-domain boundaries where natural boundary conditions are imposed, and $\mathbf{H}^f$ and $\mathbf{H}^s$ contain the prescribed values of the natural boundary conditions.  $\mathbf{A}_0$ and $\mathbf{A}_i$ are the so-called Euler jacobian matrices, $\mathbf{F}^{p}_{i}$, $\mathbf{F}^{d}_{i}$, and $\mathbf{F}^{\sigma}_{i}$ are the pressure, viscous/thermal, and stress fluxes, respectively, and $\mathbf{S}$ is the volume source. (See~\cite{bazilevs2017new1} and references therein for the details.) The subscript $\omega$ on the semilinear forms and linear functionals denotes the domain of integration, comma denotes partial differentiation, and $i=1,\dots,d$, where $d=2,3$ is the space dimension.

With the above definitions, following~\cite{bazilevs2017new1} the coupled FSI formulation may be written as: Find $\mathbf{Y} \in \mathcal{S}$, such that $\forall \mathbf{W} \in \mathcal{V}$,
\begin{align}
\label{EqFSIAlt}
 M^f_{\Omega}(\mathbf{W},\mathbf{Y}) &+ B^f_{\Omega}(\mathbf{W},\mathbf{Y}) - F^f_{\Omega}(\mathbf{W}) \nonumber \\
& \textcolor{red}{+} \nonumber \\
 M^s_{\Omega^s}(\mathbf{W},\mathbf{Y}) &+ B^s_{\Omega^s}(\mathbf{W},\mathbf{Y}) - F^s_{\Omega^s}(\mathbf{W}) \nonumber \\
& \textcolor{red}{-} \nonumber \\
\Big(M^f_{\Omega^s}(\mathbf{W},\mathbf{Y}) &+ B^f_{\Omega^s}(\mathbf{W},\mathbf{Y}) - F^f_{\Omega^s}(\mathbf{W})\Big) \nonumber \\
&\textcolor{red}{=} \nonumber \\
& 0,
\end{align}
where the integration over the fluid mechanics domain is replaced by the integration over the combined domain \textit{minus} that over the solid domain. This form of the coupled problem statement is convenient for the application of an immersed approach to the discretization of the coupled FSI equations (see, e.g.,~\cite{casquero2015a}). It was shown in the original reference~\cite{bazilevs2017new1} that the above coupled formulation satisfies a.) The fluid and solid governing equations in their respective subdomains; and b.) The compatibility of the kinematics and tractions at the fluid-solid interface.

The formulation given by Equation~(\ref{EqFSIAlt}) may be employed directly in the case the weak form of the solid problem can be written in terms of the test and trial functions coming from the spaces $\mathcal{V}$ and $\mathcal{S}$ defined over $\Omega$. However, in order to have more flexibility for the solid problem  formulation and discretization, the above formulation needs to be modified, which we do as follows. First, we assume that the solid formulation is given in the weak form that is defined over a different set of test and trial functions\footnote{e.g., shell-structure formulations using displacement and rotational unknowns. While natural to shells, rotational unknowns are not natural to fluids or solids.}: Find $\tilde{\mathbf{Y}} \in \tilde{\mathcal{S}}$, such that $\forall \tilde{\mathbf{W}} \in \tilde{\mathcal{V}}$,
\begin{equation}
\label{EqStructAlt}
\mathcal{M}^s_{\Omega^s}(\tilde{\mathbf{W}},\tilde{\mathbf{Y}}) + \mathcal{B}^s_{\Omega^s}(\tilde{\mathbf{W}},\tilde{\mathbf{Y}}) - \mathcal{F}^s_{\Omega^s}(\tilde{\mathbf{W}}) = 0.
\end{equation}
We also assume that we can develop a relationship between the functions in $\mathcal{S}$ and $\tilde{\mathcal{S}}$, and the functions in $\mathcal{V}$ and $\tilde{\mathcal{V}}$, by defining the mappings
\begin{equation}
\label{EqProjSol}
\tilde{\mathbf{Y}} = \Pi \mathbf{Y}
\end{equation}
and
\begin{equation}
\label{EqProjTest}
\tilde{\mathbf{W}} = \Pi \mathbf{W}.
\end{equation}
In the discrete setting these mappings may correspond to interpolation or projection, and will be detailed in what follows. 

These assumptions allow us to re-write the solid problem as: Find $\mathbf{Y} \in \mathcal{S}$, such that $\forall \mathbf{W} \in \mathcal{V}$,
\begin{equation}
\label{EqStructAltYW}
\mathcal{M}^s_{\Omega^s}(\Pi \mathbf{W},\Pi \mathbf{Y}) + \mathcal{B}^s_{\Omega^s}(\Pi \mathbf{W},\Pi \mathbf{Y}) - \mathcal{F}^s_{\Omega^s}(\Pi \mathbf{W}) = 0.
\end{equation}
As a consequence, the coupled FSI formulation given by Equation~(\ref{EqFSIAlt}) may be modified as: Find $\mathbf{Y} \in \mathcal{S}$, such that $\forall \mathbf{W} \in \mathcal{V}$,
\begin{align}
\label{EqFSIAltYW}
M^f_{\Omega}(\mathbf{W},\mathbf{Y}) &+ B^f_{\Omega}(\mathbf{W},\mathbf{Y}) - F^f_{\Omega}(\mathbf{W}) \nonumber \\
& \textcolor{red}{+} \nonumber \\
\mathcal{M}^s_{\Omega^s}(\Pi \mathbf{W},\Pi \mathbf{Y}) & + \mathcal{B}^s_{\Omega^s}(\Pi \mathbf{W},\Pi \mathbf{Y}) - \mathcal{F}^s_{\Omega^s}(\Pi \mathbf{W}) \nonumber \\
& \textcolor{red}{-} \nonumber \\
\Big( M^f_{\Omega^s}(\mathbf{W},\mathbf{Y}) &+ B^f_{\Omega^s}(\mathbf{W},\mathbf{Y}) - F^f_{\Omega^s}(\mathbf{W}) \Big) \nonumber \\
& \textcolor{red}{=} \nonumber \\ 
&0,
\end{align}
which, as will be shown in what follows, enables us to flexibly incorporate a much larger class of solid and structural formulations in our immersed FSI framework.

\subsection{Correspondence-based PD framework for solid modeling}
\label{sec:correspondencePD}

In this section we focus on the PD formulation of the solid. The rate form of the energy balance law on $\Omega^s$ may be expressed as~\cite{silling2007peridynamic,behzadinasab2021shell}:
\begin{equation} 
\int_{\Omega^s} \rho^s \, \dot{\mathbf{v}} \cdot \mathbf{v} \, {\rm d}\Omega^s + \int_{\Omega^s} \dot{\mathcal{U}} \, {\rm d}\Omega^s = \int_{\Omega^s} \mathbf{s} \cdot \mathbf{v} \, {\rm d}\Omega^s .
\label{eqn:energy_balace}
\end{equation}
Here, the first and second terms on the left hand side denote the rate of the total kinetic and internal energy, respectively, while the right hand side gives the total supplied power. Here, $\rho^s$ is the solid mass density in the current configuration, $\mathbf{v}$ is the velocity field, and $\mathbf{s}$ is the volume source term. The overdot symbol indicates the time derivative. In the PD formulation the strain-power density $\dot{\mathcal{U}}$ at a material point $\mathbf{x}$ is regularized as~\cite{silling2010peridynamic,behzadinasab2020peridynamic}
\begin{equation}   
\dot{\mathcal{U}}(\mathbf{x}) = \int_{\mathcal{H}(\mathbf{x})} \alpha \an{\mathbf{x}-\mathbf{x}'} \, \dot{\mathcal{U}} \an{\mathbf{x}-\mathbf{x}'} \, {\rm d}\mathbf{x}' .
\label{eqn:U_PD_x}
\end{equation}
Here, $\mathcal{H}(\mathbf{x})$ is the family set (also known as the horizon) of $\mathbf{x}$ defined as
\begin{equation} 
\mathcal{H}(\mathbf{x}) = \left\{ \mathbf{x}' \ | \ \mathbf{x}' \in \mathcal{H}(\mathbf{x}) \medcap \Omega^s, \, 0 < | \mathbf{x}' - \mathbf{x} | \leq \delta \right\} ,
\label{eqn:family}
\end{equation}
where $\delta$ is the horizon size, ${\an{\mathbf{x}-\mathbf{x}'}}$ denotes a PD bond between $\mathbf{x}$ and $\mathbf{x}'$ in the horizon, and the strain-power density ${\dot{\mathcal{U}}} \an{\mathbf{x}-\mathbf{x}'}$ is given as a bond-level quantity. The normalized weighting function $\alpha \an{\mathbf{x}-\mathbf{x}'}$ is constrained to satisfy the equation
\begin{equation}
\int_{\mathcal{H}(\mathbf{x})} \alpha \an{\mathbf{x}-\mathbf{x}'} \, {\rm d}\mathbf{x}' = 1.
\label{eqn:identity}
\end{equation}
For simplicity of notation, to distinguish the material-point- and bond-level quantities, the bond-associated fields are marked using underscores in the remainder of this paper. For example, \cref{eqn:U_PD_x} may be re-written as
\begin{equation}   
\dot{\mathcal{U}} = \int_{\mathcal{H}} \underline{\alpha} \, \underline{\dot{\mathcal{U}}} \, {\rm d}\mathcal{H} ,
\label{eqn:U_PD}
\end{equation}
where $\dot{\mathcal{U}}=\dot{\mathcal{U}}(\mathbf{x})$, $\underline{\alpha}=\alpha \an{\mathbf{x}-\mathbf{x}'}$, $\underline{\dot{\mathcal{U}}}=\dot{\mathcal{U}} \an{\mathbf{x}-\mathbf{x}'}$, and $\mathcal{H}=\mathcal{H}(\mathbf{x})$. The bond-associated strain-power density is given by
\begin{equation} 
\underline{\dot{\mathcal{U}}} = \underline{\boldsymbol{\sigma}} : \underline{\mathbf{L}} ,
\label{eqn:U_classical}
\end{equation}
where $\underline{\boldsymbol{\sigma}}$ and $\underline{\mathbf{L}}$ are the bond-associated power-conjugate Cauchy stress and velocity gradient tensors, respectively. The total stress power, using \cref{eqn:energy_balace,eqn:U_PD,eqn:U_classical}, may be re-expressed as
\begin{equation}   
\int_{\Omega^s} \dot{\mathcal{U}} \, {\rm d}\Omega^s = \int_{\Omega^s} \int_{\mathcal{H}} \underline{\alpha} \ \underline{\boldsymbol{\sigma}} : \underline{\mathbf{L}} \, {\rm d}\mathcal{H} \, {\rm d}\Omega^s 
= \int_{\Omega^s} \int_{\mathcal{H}} \underline{\mathbf{T}} \cdot \underline{\mathbf{v}} \, {\rm d}\mathcal{H} \, {\rm d}\Omega^s ,
\label{eqn:U_stress_based}
\end{equation}
where $\underline{\mathbf{T}} = \mathbf{T} \an{\mathbf{x}-\mathbf{x}'}$ is the so-called PD force state (with units of force per unit volume squared). The PD velocity state $\underline{\mathbf{v}}$ is defined as
\begin{equation}   
\underline{\mathbf{v}} = \mathbf{v}' - \mathbf{v} = \mathbf{v}(\mathbf{x}') - \mathbf{v}(\mathbf{x}).
\label{eqn:v_state}
\end{equation}
In the correspondence PD framework, the spatial gradients are computed using an integral form. Using the computed velocity gradient, a classical constitutive law is employed to evaluate the Cauchy stress tensor. As a result, the way $\underline{\mathbf{L}}$ is calculated determines the relation between the force state and Cauchy stress at the bond level, which we detail in the Appendix. 

Re-writing \cref{eqn:energy_balace} using \cref{eqn:U_stress_based,eqn:v_state} yields
\begin{equation}
\int_{\Omega^s} \left( \rho^s \, \dot{\mathbf{v}} - \int_{\mathcal{H}} \left( \underline{\mathbf{T}} - \underline{\mathbf{T}'} \right) \, {\rm d}\mathcal{H} - \mathbf{s} \right) \cdot \mathbf{v} \, {\rm d}\Omega^s = 0,
\label{eqn:PD_energy_balace}
\end{equation}
where $\underline{\mathbf{T}'} = \mathbf{T} \an{\mathbf{x}'-\mathbf{x}}$. Noting that \cref{eqn:PD_energy_balace} must hold {\em for all} choices of $\mathbf{v}$ enables us to recast the PD linear-momentum balance equation in the weak form given by Equation~(\ref{EqStructAltYW}) with the following definitions of the semi-linear and linear forms:
\begin{equation}
\label{EqWFPDM}
\mathcal{M}^s_{\Omega^s}(\tilde{\mathbf{W}},\tilde{\mathbf{Y}}) = \int_{\Omega^s}  \tilde{\mathbf{w}} \cdot \rho^s \, \dot{\tilde{\mathbf{v}}} \, {\rm d}\Omega^s,
\end{equation}
\begin{equation}
\label{EqWFPDB}
\mathcal{B}^s_{\Omega^s}(\tilde{\mathbf{W}},\tilde{\mathbf{Y}}) = \int_{\Omega^s}  \tilde{\mathbf{w}} \cdot \int_{\mathcal{H}} \left( \underline{\mathbf{T}} - \underline{\mathbf{T}'} \right) \, {\rm d}\mathcal{H} \, {\rm d}\Omega^s,
\end{equation}
and
\begin{equation}
\label{EqWFPDF}
\mathcal{F}^s_{\Omega^s}(\tilde{\mathbf{W}}) = \int_{\Omega^s}  \tilde{\mathbf{w}} \cdot \mathbf{s} \, {\rm d}\Omega^s.
\end{equation}
It is important to note that no higher regularity than $L^2$ is required for $\tilde{\mathbf{v}}$ and $\tilde{\mathbf{w}}$ in order for the weak form of the PD linear-momentum balance equation to be well defined. Also note that the solid is assumed to be isothermal and the mass balance is satisfied in the Lagrangian description at each material point, resulting in the PD formulation that only makes use of the momentum-balance PDE written in terms of the velocity unknowns. The extension to a thermally coupled solid does not present a conceptual difficulty; however, it is not pursued here.

\subsection{Discrete form of the coupled FSI problem}

To develop a discrete version of the FSI formulation, we first define $\mathcal{S}^h \subset \mathcal{S}$ and $\mathcal{V}^h \subset \mathcal{V}$, the background-domain finite-dimensional trial- and test-function spaces, respectively. The discrete trial and test functions may be expressed as
\begin{equation}
\mathbf{Y}^h(\mathbf{x}) = \sum_{B=1}^{\mathcal{N}_{cp}} \mathbf{Y}_B \, N_B(\mathbf{x})
\label{eqn:backgroundtri}
\end{equation}
and
\begin{equation}
\mathbf{W}^h(\mathbf{x}) = \sum_{A=1}^{\mathcal{N}_{cp}} \mathbf{W}_A \, N_A(\mathbf{x}),
\label{eqn:backgroundtest}
\end{equation}
where $\mathbf{Y}_B$ and $\mathbf{W}_A$ are the vector-valued control-point degrees of freedom (DOFs) and weights, respectively, $N(\mathbf{x})$'s in this work are assumed to be the B-Spline basis functions defined everywhere in $\Omega$, and $\mathcal{N}_{cp}$ is the dimension of the B-Spline space. Note that all the components of the trial and test-function vectors are approximated using the same basis functions.

We also define the finite-dimensional trial and test function spaces $\tilde{\mathcal{S}}$ and $\tilde{\mathcal{V}}$, respectively, for the PD solid. The discrete trial and test functions may be expressed as
\begin{equation}
\tilde{\mathbf{Y}}^h(\mathbf{x}) = \sum_{Q=1}^{\mathcal{N}_{mp}} \tilde{\mathbf{Y}}_Q \, \chi_Q(\mathbf{x})
\label{eqn:foregroundtri}
\end{equation}
and
\begin{equation}
\tilde{\mathbf{W}}^h(\mathbf{x}) = \sum_{P=1}^{\mathcal{N}_{mp}} \tilde{\mathbf{W}}_P \, \chi_P(\mathbf{x}),
\label{eqn:foregroundtest}
\end{equation}
where the solid domain is represented using a finite number $\mathcal{N}_{mp}$ of material points or PD nodes, $\tilde{\mathbf{Y}}_Q$ and $\tilde{\mathbf{W}}_P$ are the nodal DOFs and weights, respectively, and $\chi_P(\mathbf{x})$ is a characteristic function of a PD node $P$ that attains a value of unity at $\mathbf{x}_P \in \Omega^s$, the spatial location of the PD node, is zero at all other PD nodes, and satisfies
\begin{equation}
\int_{\Omega^s} \chi_P(\mathbf{x}) \, {\rm d}\Omega^s = V_P,
\end{equation}
where $V_P$ is the local volume of the PD node. Note that the exact shape of the characteristic functions is not necessary for our purposes because, ultimately, we will make use of nodal quadrature to evaluate the integrals in the PD formulation.

We now define the mapping $\Pi$ between the discrete background and foreground functions using nodal interpolation. Namely, we constrain the PD nodal DOFs and weights to the background test and trial functions, respectively, as 
\begin{equation}
\tilde{\mathbf{Y}}_Q = \mathbf{Y}^h(\mathbf{x}_Q) = \sum_{B=1}^{\mathcal{N}_{cp}} \mathbf{Y}_B \, N_B(\mathbf{x}_Q)
\label{eq:trialprojcoeff}
\end{equation}
and
\begin{equation}
\tilde{\mathbf{W}}_P = \mathbf{W}^h(\mathbf{x}_P) = \sum_{A=1}^{\mathcal{N}_{cp}} \mathbf{W}_A \, N_A(\mathbf{x}_P),
\label{eq:testprojcoeff}
\end{equation}
which yields
\begin{equation}
\tilde{\mathbf{Y}}^h (\mathbf{x}) = \Pi \mathbf{Y}^h(\mathbf{x}) = \sum_{Q=1}^{\mathcal{N}_{mp}} \left(~\sum_{B=1}^{\mathcal{N}_{cp}} \mathbf{Y}_B \, N_B(\mathbf{x}_Q)~\right) \, \chi_Q(\mathbf{x})
\label{eq:trialprojfun}
\end{equation}
and
\begin{equation}
\tilde{\mathbf{W}}^h(\mathbf{x}) = \Pi \mathbf{W}^h(\mathbf{x}) = \sum_{P=1}^{\mathcal{N}_{mp}} \left(~\sum_{A=1}^{\mathcal{N}_{cp}} \mathbf{W}_A \, N_A(\mathbf{x}_P)~\right) \, \chi_P(\mathbf{x}),
\label{eq:testprojfun}
\end{equation}
and defines the mapping $\Pi$.

With the above definitions of the test and trial functions, the spatially-discretized immersed FSI formulation may be stated as: Find $\mathbf{Y}^h \in \mathcal{S}^h$, such that $\forall \mathbf{W}^h \in \mathcal{V}^h$, 
\begin{align}
\label{eqn:FSIDisc}
M^f_{\Omega}(\mathbf{W}^h,\mathbf{Y}^h) & + B^f_{\Omega}(\mathbf{W}^h,\mathbf{Y}^h) - 
F^f_{\Omega}(\mathbf{W}^h) + B^{st}_{\Omega}(\mathbf{W}^h,\mathbf{Y}^h) + B^{dc}_{\Omega}(\mathbf{W}^h,\mathbf{Y}^h) \nonumber \\
& \textcolor{red}{+} \nonumber \\
\mathcal{M}^s_{\Omega^s}(\tilde{\mathbf{W}}^h,\tilde{\mathbf{Y}}^h) & + \mathcal{B}^s_{\Omega^s}(\tilde{\mathbf{W}}^h,\tilde{\mathbf{Y}}^h) - \mathcal{F}^s_{\Omega^s}(\tilde{\mathbf{W}}^h) \nonumber \\
& \textcolor{red}{-} \nonumber \\
(~M^f_{\Omega^s}(\mathbf{W}^h,\mathbf{Y}^h) & + B^f_{\Omega^s}(\mathbf{W}^h,\mathbf{Y}^h) - 
F^f_{\Omega^s}(\mathbf{W}^h) + B^{st}_{\Omega^s}(\mathbf{W}^h,\mathbf{Y}^h) + B^{dc}_{\Omega^s}(\mathbf{W}^h,\mathbf{Y}^h)~) \nonumber \\
& \textcolor{red}{=} \nonumber \\ 
&0.
\end{align}
Note that, in the space-discrete case, we introduce the SUPG stabilization ($B^{st}_{\Omega}(\mathbf{W}^h,\mathbf{Y}^h)$)~\cite{brooks1982streamline,le1993supg,tezduyar2006stabilization,hughes2010stabilized}  and  discontinuity-capturing ($B^{dc}_{\Omega}(\mathbf{W}^h,\mathbf{Y}^h)$)~\cite{hughes1986new,tezduyar2006computation,rispoli2009computation,rispoli2015particle} operators to the compressible flow formulation to address the convective instability of the Galerkin technique in the regime of convection dominance and to provide additional dissipation in the shock regions. Both the stabilization and shock-capturing operators are proportional to the strong-form PDE residual of the compressible-flow Navier--Stokes equations, which retains the method consistency. More details on the definition of these operators that are employed in the present formulation of the compressible flow may be found in~\cite{hauke1998comparative,xu2017compressible}.

Introducing the definitions of the test and trial functions from Equations~(\ref{eqn:backgroundtri})-(\ref{eqn:backgroundtest}),~(\ref{eqn:foregroundtri})-(\ref{eqn:foregroundtest}), and~(\ref{eq:trialprojcoeff})-(\ref{eq:testprojfun}) into the space-discrete weak form given by Equation~(\ref{eqn:FSIDisc}) gives
\begin{align}
\label{eqn:FSIDiscWeights}
\sum_{A=1}^{N_{cp}} W_{AI}~[~M^f_{\Omega}(N_A \mathbf{e}_I,\mathbf{Y}^h) & + B^f_{\Omega}(N_A \mathbf{e}_I,\mathbf{Y}^h) - 
F^f_{\Omega}(N_A \mathbf{e}_I) + 
B^{st}_{\Omega}(N_A \mathbf{e}_I,\mathbf{Y}^h) + 
B^{dc}_{\Omega}(N_A \mathbf{e}_I,\mathbf{Y}^h) \nonumber \\
& \textcolor{red}{+} \nonumber \\
N_A(\mathbf{x}_P)~(~\mathcal{M}^s_{\Omega^s}(\chi_P \mathbf{e}_I,\tilde{\mathbf{Y}}^h) & + \mathcal{B}^s_{\Omega^s}(\chi_P \mathbf{e}_I,\tilde{\mathbf{Y}}^h) - \mathcal{F}^s_{\Omega^s}(\chi_P \mathbf{e}_I)~) \nonumber \\
& \textcolor{red}{-} \nonumber \\
(~M^f_{\Omega^s}(N_A \mathbf{e}_I,\mathbf{Y}^h) & + B^f_{\Omega^s}(N_A \mathbf{e}_I,\mathbf{Y}^h) - 
F^f_{\Omega^s}(N_A \mathbf{e}_I) + 
B^{st}_{\Omega^s}(N_A \mathbf{e}_I,\mathbf{Y}^h) + B^{dc}_{\Omega^s}(N_A \mathbf{e}_I,\mathbf{Y}^h)~)~] \nonumber \\
& \textcolor{red}{=} \nonumber \\ 
&0,
\end{align}
where $I = 1,2,\dots,d+2$ is the nodal or control-point DOF index and $\mathbf{e}_I$ is a Cartesian basis vector of dimension $(d+2)$ with one in the $I^{th}$ entry and zero in the remaining entries. 

Because Equation~(\ref{eqn:FSIDiscWeights}) holds \textit{for all} background-domain discrete weights $W_{AI}$, we arrive at the following vector form of the semi-discrete FSI problem: Find the background-domain control variables $\mathbf{Y}_A$ such that
\begin{equation}
R^f_{AI} + R^s_{AI} = 0 , \quad A=1,...,\mathcal{N}_{cp} , \quad I=1,2,...,d+2 ,
\label{eqn:FSIDiscWeightsReduced}
\end{equation}
where
\begin{align}
R^f_{AI} = M^f_{\Omega}(N_A \mathbf{e}_I,\mathbf{Y}^h) & + B^f_{\Omega}(N_A \mathbf{e}_I,\mathbf{Y}^h) - 
F^f_{\Omega}(N_A \mathbf{e}_I) + 
B^{st}_{\Omega}(N_A \mathbf{e}_I,\mathbf{Y}^h) + 
B^{dc}_{\Omega}(N_A \mathbf{e}_I,\mathbf{Y}^h) \nonumber \\
& \textcolor{red}{-} \nonumber \\
(~M^f_{\Omega^s}(N_A \mathbf{e}_I,\mathbf{Y}^h) & + B^f_{\Omega^s}(N_A \mathbf{e}_I,\mathbf{Y}^h) - 
F^f_{\Omega^s}(N_A \mathbf{e}_I) + 
B^{st}_{\Omega^s}(N_A \mathbf{e}_I,\mathbf{Y}^h) + B^{dc}_{\Omega^s}(N_A \mathbf{e}_I,\mathbf{Y}^h)~)
\end{align}
and
\begin{equation}
R^s_{AI} = \sum_{P=1}^{N_mp} N_A(\mathbf{x}_P) \tilde{R}^s_{PI},
\label{eqn:ForceTransfer}
\end{equation}
where
\begin{equation}
\tilde{R}^s_{PI} = \mathcal{M}^s_{\Omega^s}(\chi_P \mathbf{e}_I,\tilde{\mathbf{Y}}^h) + \mathcal{B}^s_{\Omega^s}(\chi_P \mathbf{e}_I,\tilde{\mathbf{Y}}^h) - \mathcal{F}^s_{\Omega^s}(\chi_P \mathbf{e}_I).
\label{eqn:SolidPDRes}
\end{equation}
Equation~(\ref{eqn:FSIDiscWeightsReduced}) states that at each control point the generalized force vectors coming from the fluid and solid parts of the problem balance; Equation~(\ref{eqn:SolidPDRes}) states that the PD solid nodal forces may be computed \textit{directly from the PD formulation on the foreground mesh} provided the kinematic data can be interpolated on the PD nodes from the background solution field; and Equation~(\ref{eqn:ForceTransfer}) defines a consistent transformation of the force vector from the foreground to the background grid. 

Introducing the definitions of the semi-linear and linear forms of the PD formulation given by Equations~(\ref{EqWFPDM})-(\ref{EqWFPDF}) into Equation~(\ref{eqn:SolidPDRes}) yields an explicit expression for the PD force vector on the foreground grid:
\begin{equation}
\tilde{\mathbf{R}}^s = [\tilde{R}_{PI}] = 
\begin{bmatrix}
0 \\
\left( \rho^s_P \, \dot{\mathbf{v}}_P - \int_{\mathcal{H}_P} \left( \underline{\mathbf{T}} - \underline{\mathbf{T}'} \right) \, {\rm d}\mathcal{H} - \mathbf{s}_P \right) V_P \\
0
\end{bmatrix},
\label{eqn:PDNodalForce}
\end{equation}
where $\mathcal{H}_P$ is the horizon of the PD node located at $\mathbf{x}_P$ and the integral is evaluated using nodal quadrature. Note that in the present model the PD solid will only contribute to the linear-momentum equation balance of the discrete FSI formulation given by Equation~(\ref{eqn:FSIDiscWeightsReduced}).

\remark{The formulation developed in this section presents a minimally intrusive approach to immersed FSI, and opens the door to using not just PD, but any other suitable discrete formulation of solid and structural mechanics as part of the FSI coupling framework. All that is required is the ability to develop a mapping of the background solution to the DOFs employed by the solid and structural mechanics formulation.}

\subsubsection{Time integration and the lumped mass matrix}

The semi-discrete FSI equations are integrated in time using the lumped-mass explicit Generalized-$\alpha$ predictor–multi-corrector procedure~\cite{chung1993time} adopted for the immersed FSI formulation and detailed in~\cite{bazilevs2017new1}. This approach requires the computation of a lumped mass matrix for the coupled system. A consistent mass matrix for the coupled FSI problem is obtained by adding the contributions from the fluid and solid subdomains on the background mesh. Denoting by $\tilde{M}^s_{PQ}$ the coefficients of the consistent mass matrix from the PD foreground mesh and using the relations developed in this section, one can show that its counterpart on the background mesh is obtained through the transformation
\begin{equation}
M^s_{AB} = \sum_{P=1}^{N_{mp}} \sum_{Q=1}^{N_{mp}}  N_A(\mathbf{x}_P) \tilde{M}^s_{PQ} N_B(\mathbf{x}_Q). 
\end{equation}
The diagonal entries of the lumped mass $L^s_{A}$ may be obtained using a classical row-sum technique, namely $L^s_{A} = \sum_{B=1}^{N_{cp}} M^s_{AB}$. The partition-of-unity property of the B-Spline functions implies $\sum_{B=1}^{N_{cp}} N_B(\mathbf{x}_Q) = 1$ for each $Q$. This, in turn, gives
\begin{equation}
L^s_{A} = \sum_{P=1}^{N_{mp}} N_A(\mathbf{x}_P) \sum_{Q=1}^{N_{mp}} \tilde{M}^s_{PQ} = \sum_{P=1}^{N_{mp}} N_A(\mathbf{x}_P) \tilde{L}^s_{P}, 
\label{eq:lumpedbg}
\end{equation}
where $\tilde{L}^s_{P}$ are the diagonal entries of the lumped mass on the PD foreground mesh given by
\begin{equation}
\tilde{L}^s_{P} = \rho^s_P V_P.
\end{equation}
Equation~(\ref{eq:lumpedbg}) states the diagonal entries of the lumped mass matrix transform in the same way from the foreground to the background mesh as the nodal forces (see Equation~(\ref{eqn:ForceTransfer})). Because the present PD formulation only solves the momentum balance equation, $L^s_A$'s are added only to the diagonal entries $2,\dots,d+1$ of the $(d+2)\times (d+2)$ control-point block of the background lumped mass matrix.

\section{Numerical examples}
\label{sec:implementation}

Three 2D numerical simulations are provided here to showcase the capabilities of the proposed IGA-PD framework in blast loading and fragmentation events. The first example serves as a verification case for the developed formulation in solving problems involving large deformation and plasticity. The next two cases are demonstrative examples that include fracture propagation and fragmentation in brittle and ductile solids subjected to blast loading. In all the presented computations, $C^1$-continuous quadratic NURBS functions are used for the background solution. Unless otherwise noted, RK functions with quadratic consistency, rectangular support, and the BA stabilization technique~\cite{breitzman2018bond,behzadinasab2020semi} are employed in the PD formulation. The PD support size $\delta$ is chosen with respect to the discretization size $\Delta x$ and the order of the kernel function $n$ is chosen such that $\delta / \Delta x = n+1$~\cite{behzadinasab2021unifiedI}. The air properties are used for the fluid with constant viscocity  $\mu = 1.81 \times 10^{-5} \, {\rm kg}/({\rm m \, s})$, Prandtl number 0.72, and adiabatic index $\gamma=1.4$.

\subsection{Chamber detonation}
\label{sec:detonation}

\begin{figure*}[!hbpt]
  \centering
  \subfloat[][]{\includegraphics[width=0.7\textwidth,trim={0cm 0cm 0cm 0cm},clip]{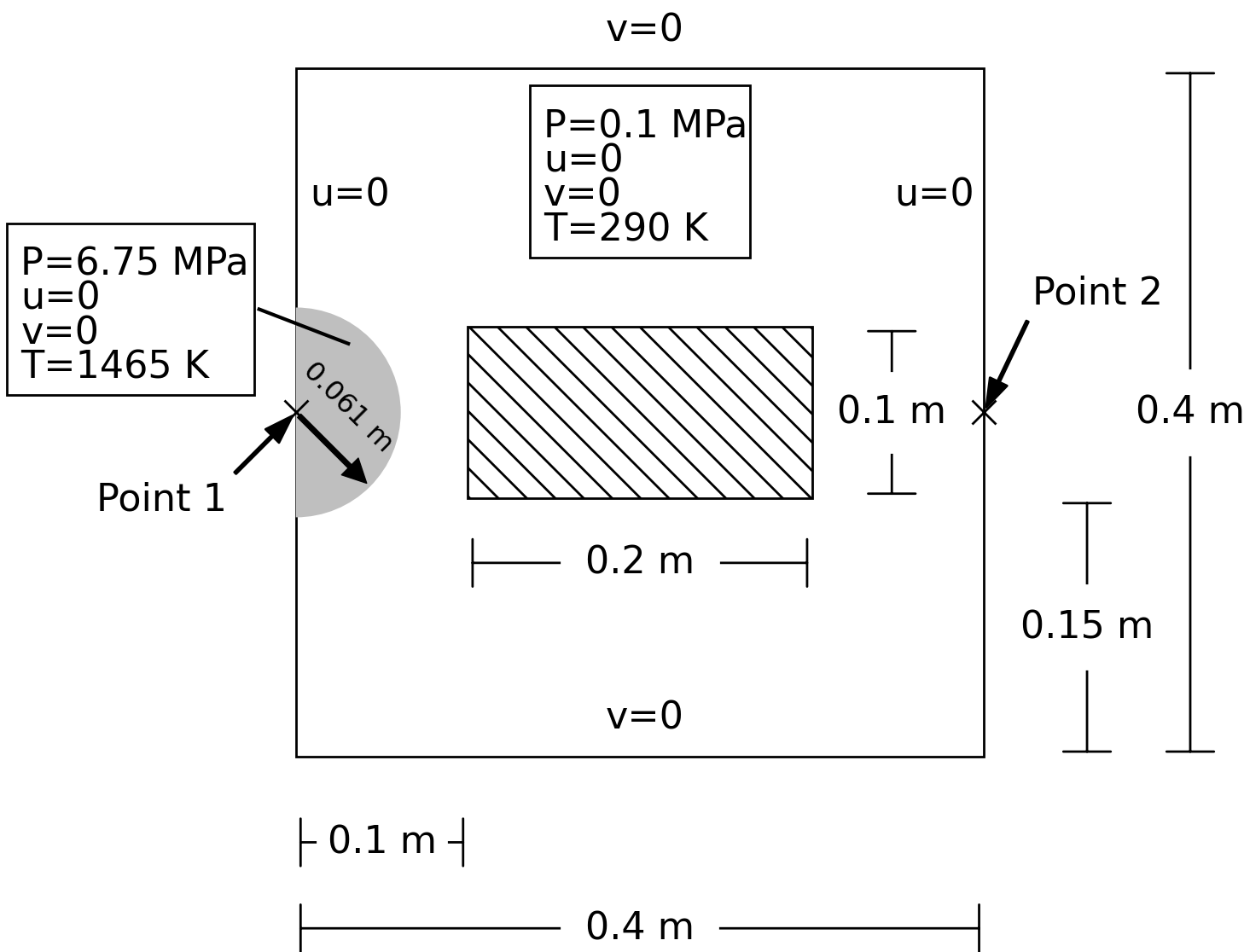}}
  \caption{Chamber detonation. Problem setup and geometry.}
  \label{fig:detonation_setup}
\end{figure*}

In this FSI example, a steel bar undergoes a detonation blast loading condition. The problem description is shown in \cref{fig:detonation_setup}, where a $0.2 \, {\rm m} \, \times \, 0.1 \, {\rm m}$ block is located at the center of a closed chamber with dimensions $0.4 \, {\rm m} \, \times \, 0.4 \, {\rm m}$.
The bar thickness is set to $3.5 \, {\rm mm}$. Isotropic linear hardening rule is used for the solid material with Young's modulus ${\rm E}=200 \, {\rm GPa}$, Poisson's ratio $\nu=0.29$, yield stress $\sigma_Y=0.4 \, {\rm GPa}$, hardening modulus $H = 0.1 \, {\rm GPa}$, and initial density $\rho^s=7870 \, {\rm kg}/{\rm m}^3$. Initially, the air in the chamber is at rest with $p = 0.1 \, {\rm MPa}$ and $T = 290 \, {\rm K}$. The detonation condition is enforced by setting higher-than-ambient values of the pressure and temperature, i.e., we set $p = 6.75 \, {\rm MPa}$ and $T = 1465 \, {\rm K}$ in a semi-circular area with a radius of $6.1 \, {\rm mm}$, centered on the left wall. No-penetration and free-slip boundary conditions are prescribed at the chamber walls. 

The IGA-PD simulations are compared with the reference immersed simulations in~\cite{bazilevs2017new1,bazilevs2017new2}, which we refer to as IGA-IGA. Three different discretization levels are considered: coarse - Fluid: $40 \times 40$ elements; Solid: $60 \times 30$ elements (PD nodes); medium - Fluid: $80 \times 80$ elements; Solid: $120 \times 60$ elements (PD nodes); fine - Fluid: $120 \times 120$ elements; Solid: $180 \times 90$ elements (PD nodes). In the PD case, each foreground element is replaced by a meshfree node at its centroid with an equivalent volume. The time step used for the coarse, medium, and fine cases are $0.5 \, {\rm ms}$, $0.25 \, {\rm ms}$, and $0.166 \, {\rm ms}$, respectively. 

\begin{figure*}[!hbpt]
  \centering
  \subfloat{\includegraphics[width=0.32\textwidth,trim={0cm 0cm 0cm 0cm},clip]{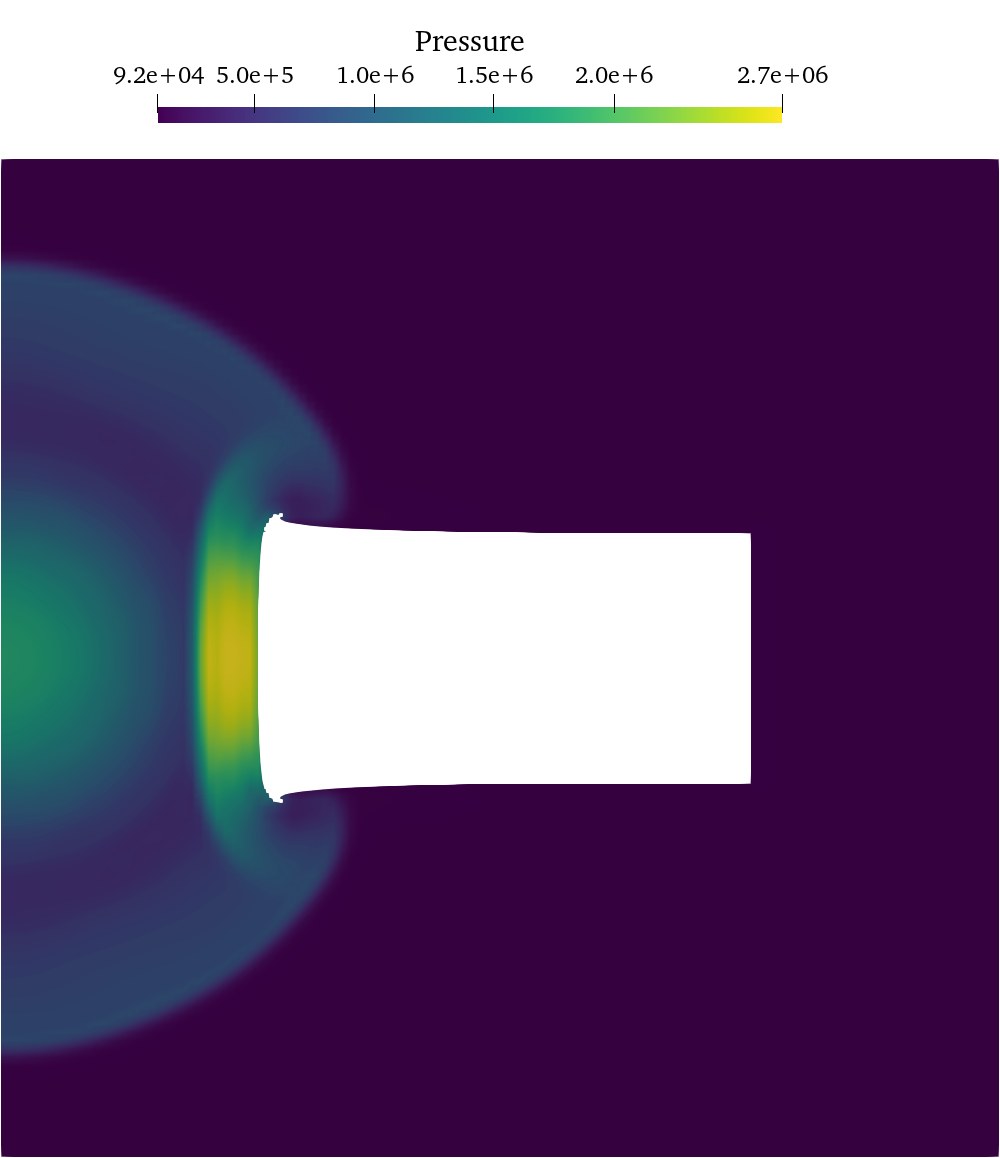}}
  \hspace{4pt}
  \subfloat{\includegraphics[width=0.32\textwidth,trim={0cm 0cm 0cm 0cm},clip]{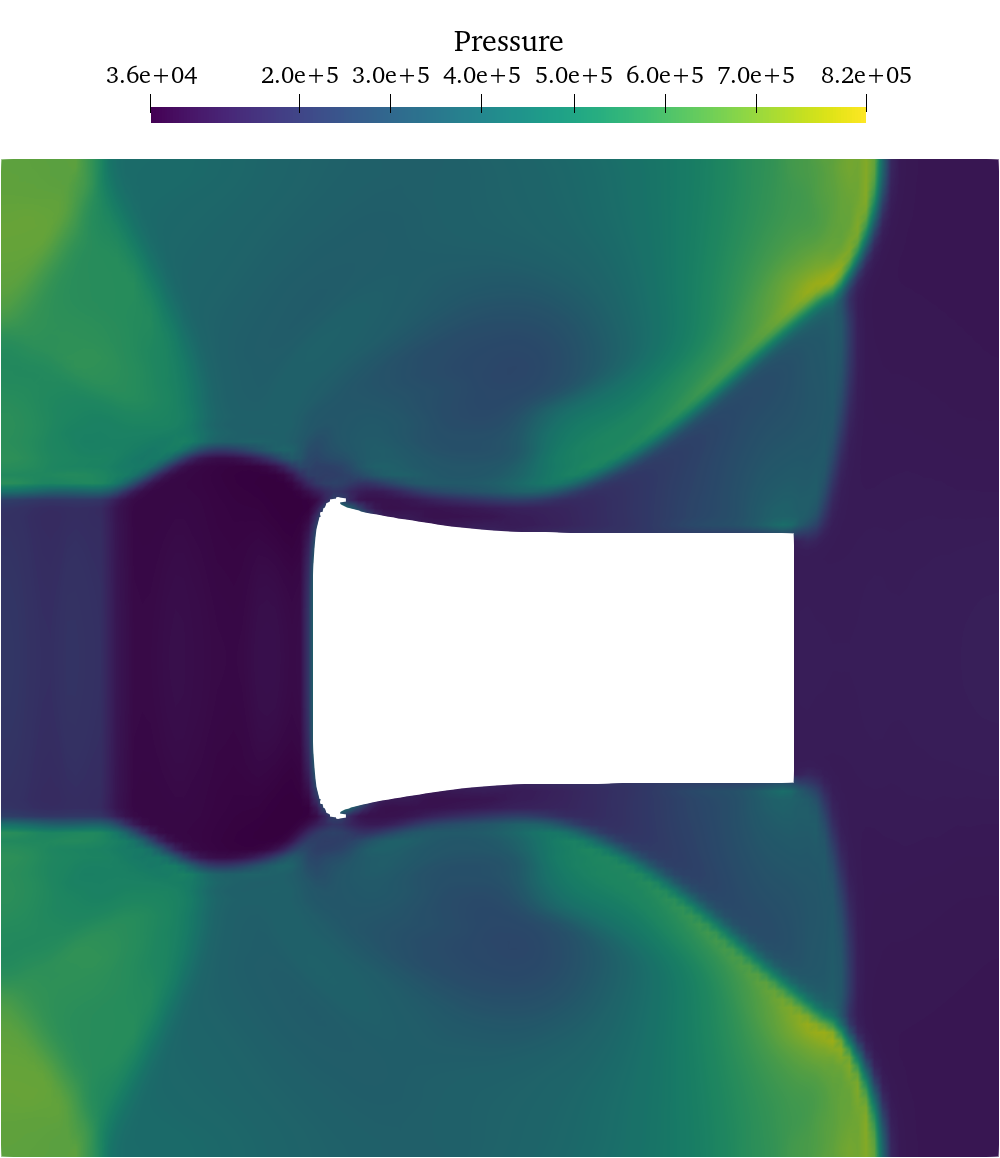}}
  \hspace{4pt}
  \subfloat{\includegraphics[width=0.32\textwidth,trim={0cm 0cm 0cm 0cm},clip]{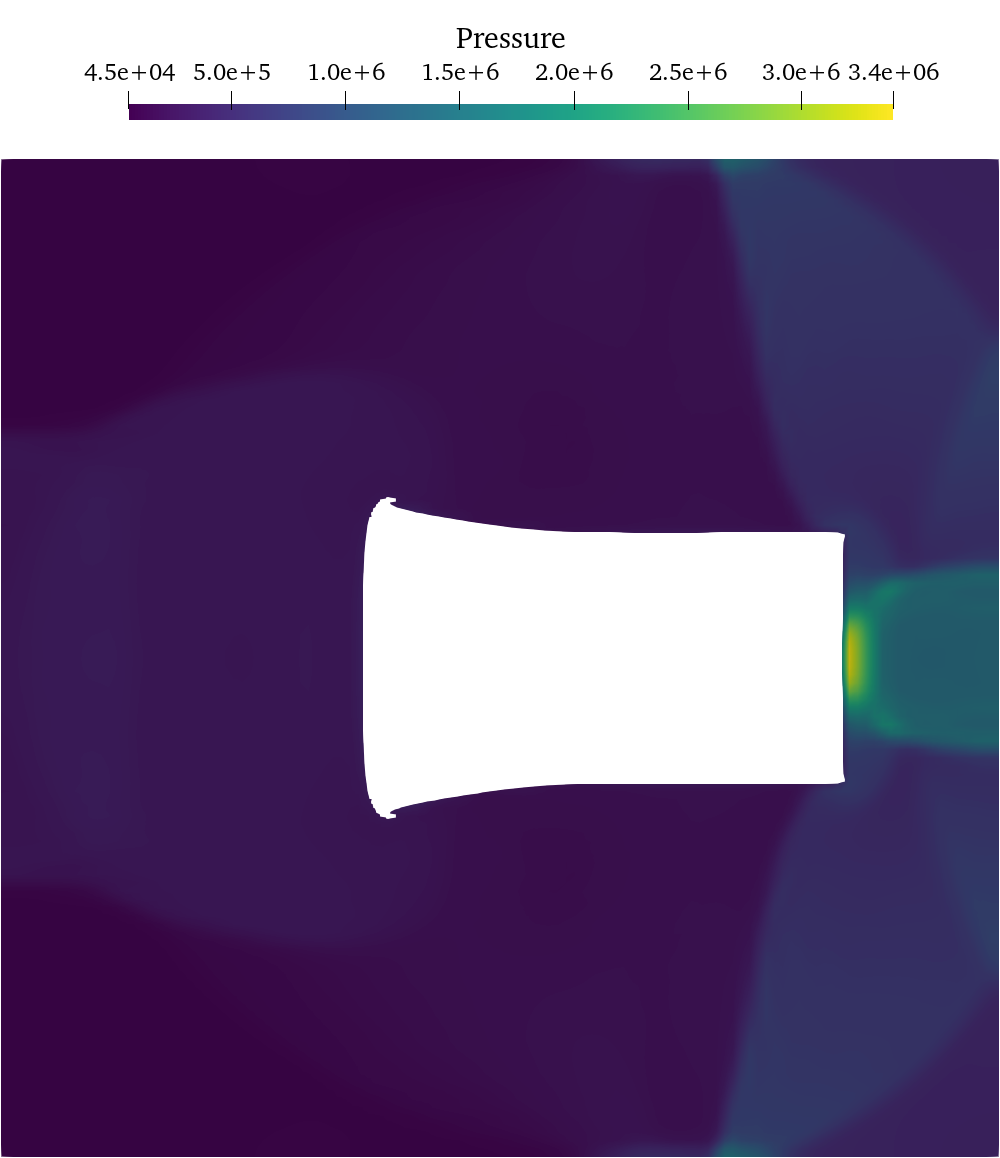}}

  \setcounter{subfigure}{0}
  \subfloat[][$t=\SI{0.1}{ms}$]{\includegraphics[width=0.32\textwidth,trim={0cm 0cm 0cm 5cm},clip]{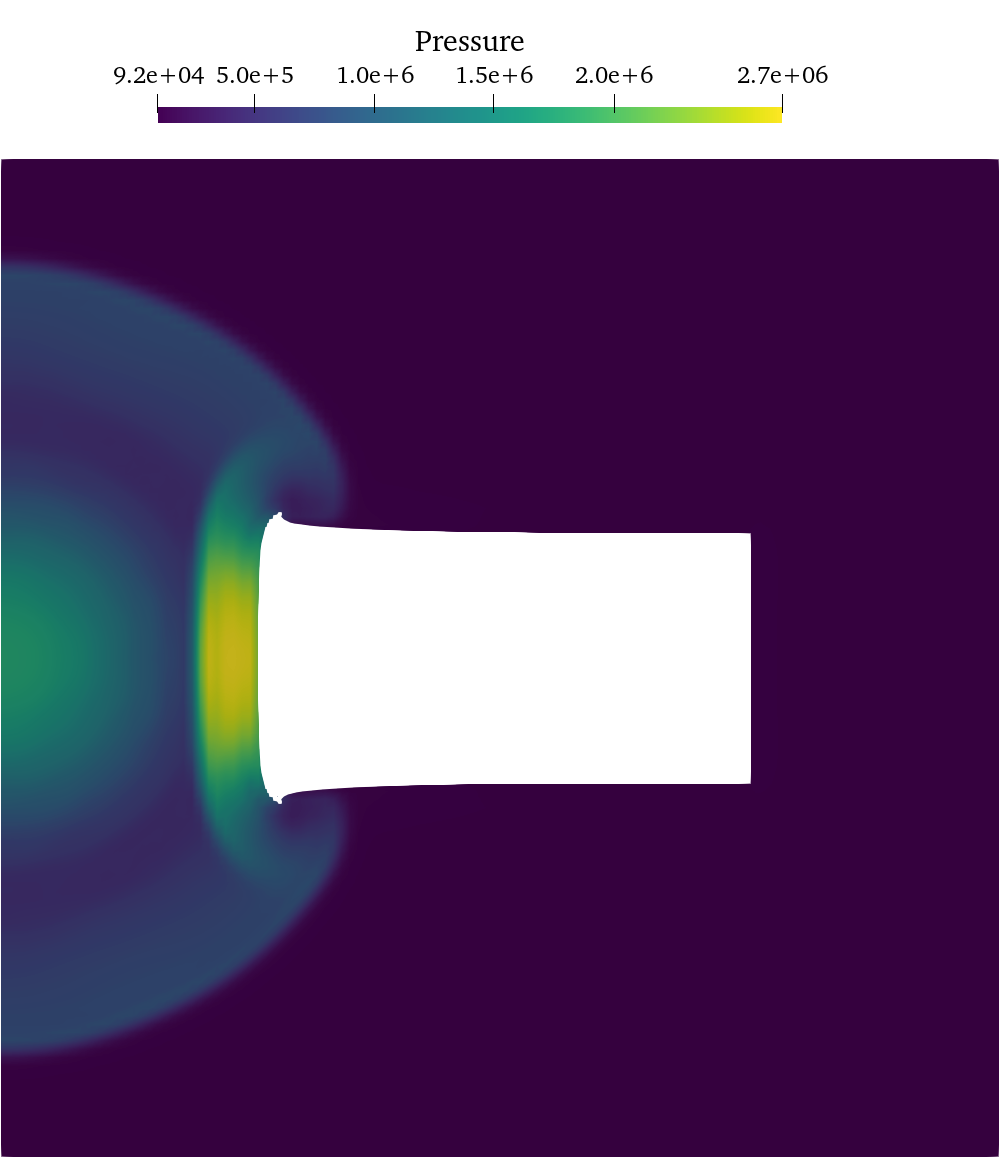}}
  \hspace{4pt}
  \subfloat[][$t=\SI{0.4}{ms}$]{\includegraphics[width=0.32\textwidth,trim={0cm 0cm 0cm 5cm},clip]{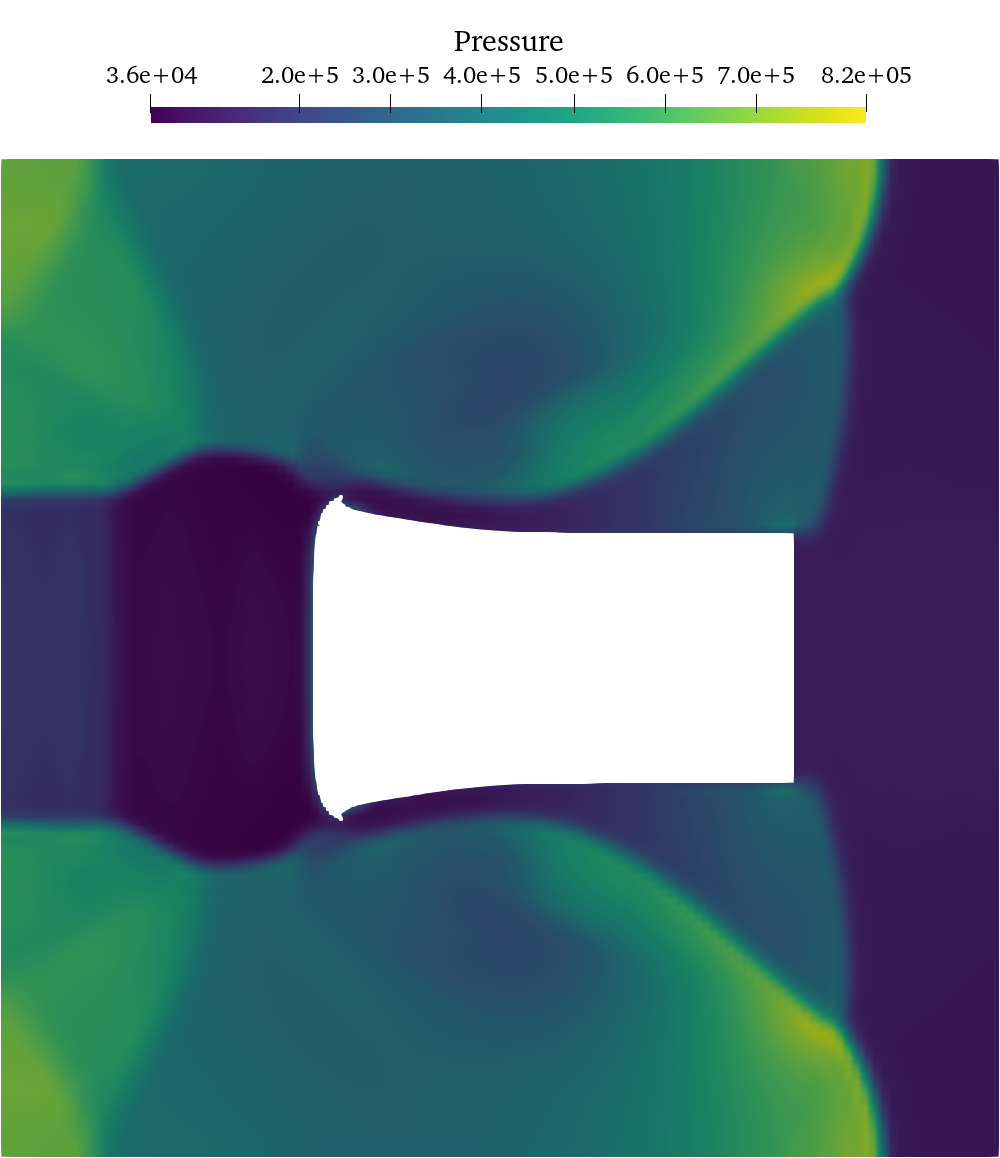}}
  \hspace{4pt}
  \subfloat[][$t=\SI{0.7}{ms}$]{\includegraphics[width=0.32\textwidth,trim={0cm 0cm 0cm 5cm},clip]{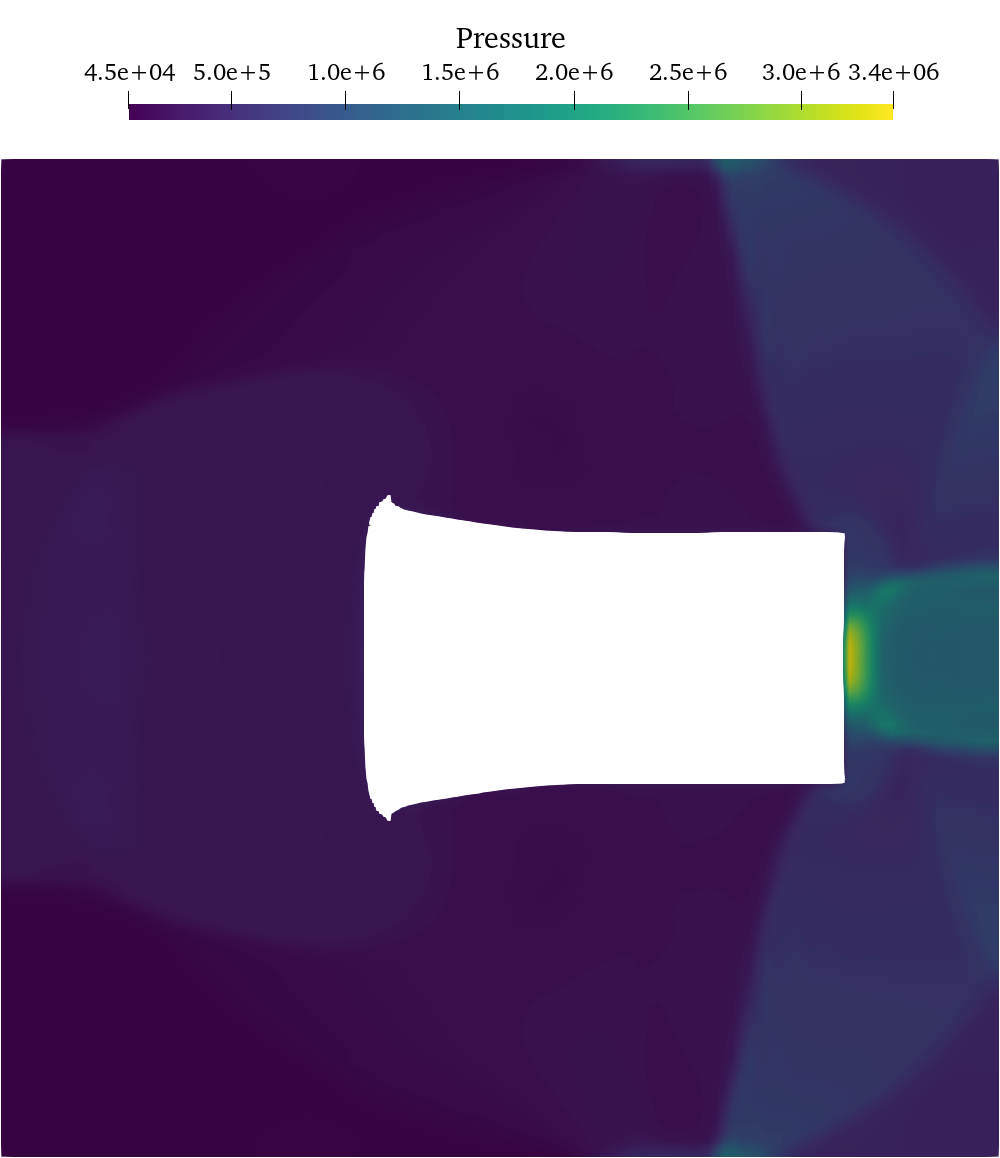}}
  \caption{Chamber detonation. Pressure at different time instants. Top and bottom include IGA solid and PD solid simulations, respectively. (a): 0.1 ms, (b): 0.4 ms and (c): 0.7 ms.}
  \label{fig:detonation_pressure}
\end{figure*}

\begin{figure*}[!hbpt]
  \centering
  \subfloat[][IGA]{\includegraphics[height=0.19\textwidth,height=0.31\textwidth,trim={0cm 0cm 6cm 0cm},clip]{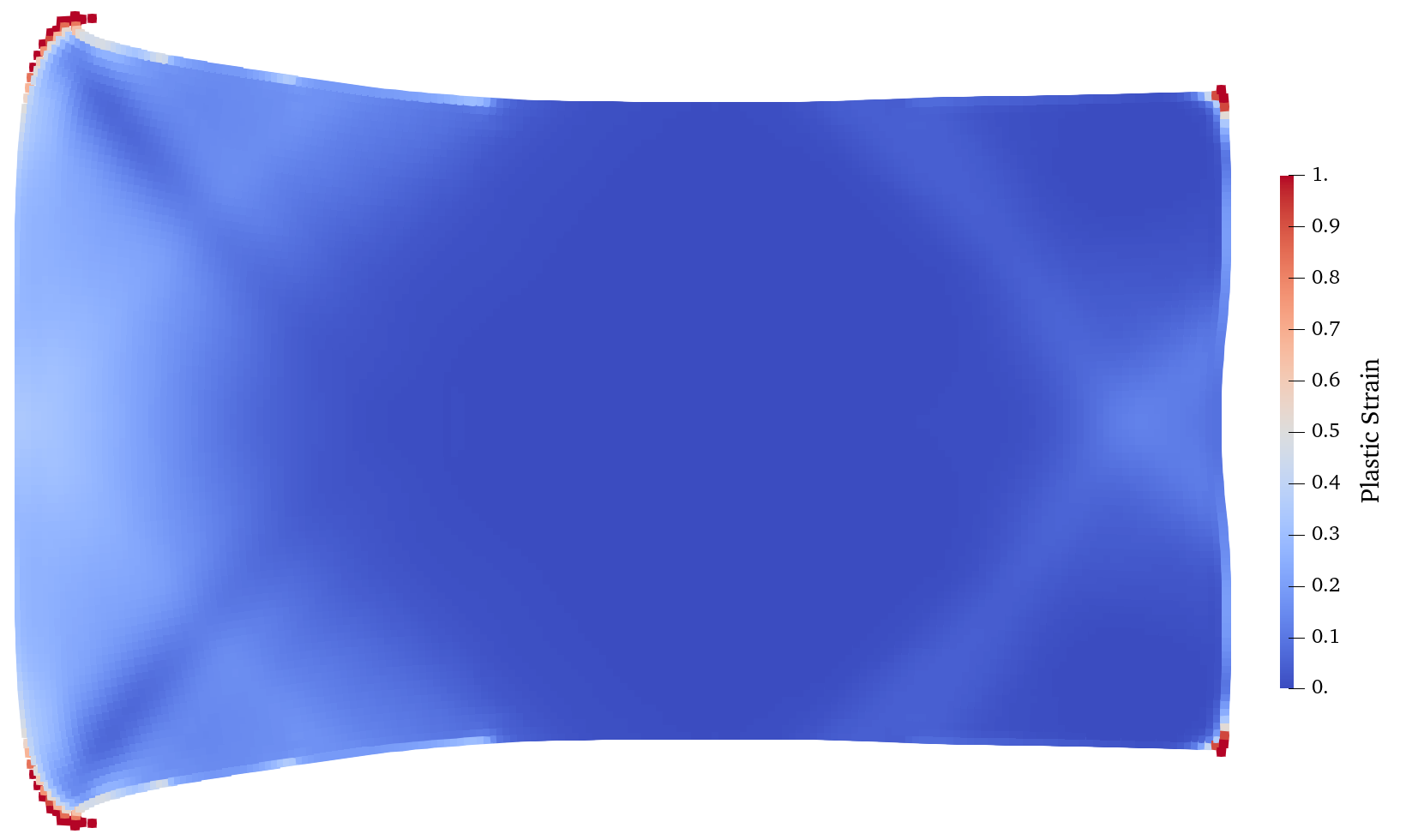}}
  \hfill
  \subfloat[][PD]{\includegraphics[height=0.19\textwidth,height=0.31\textwidth,trim={0cm 0cm 0cm 0cm},clip]{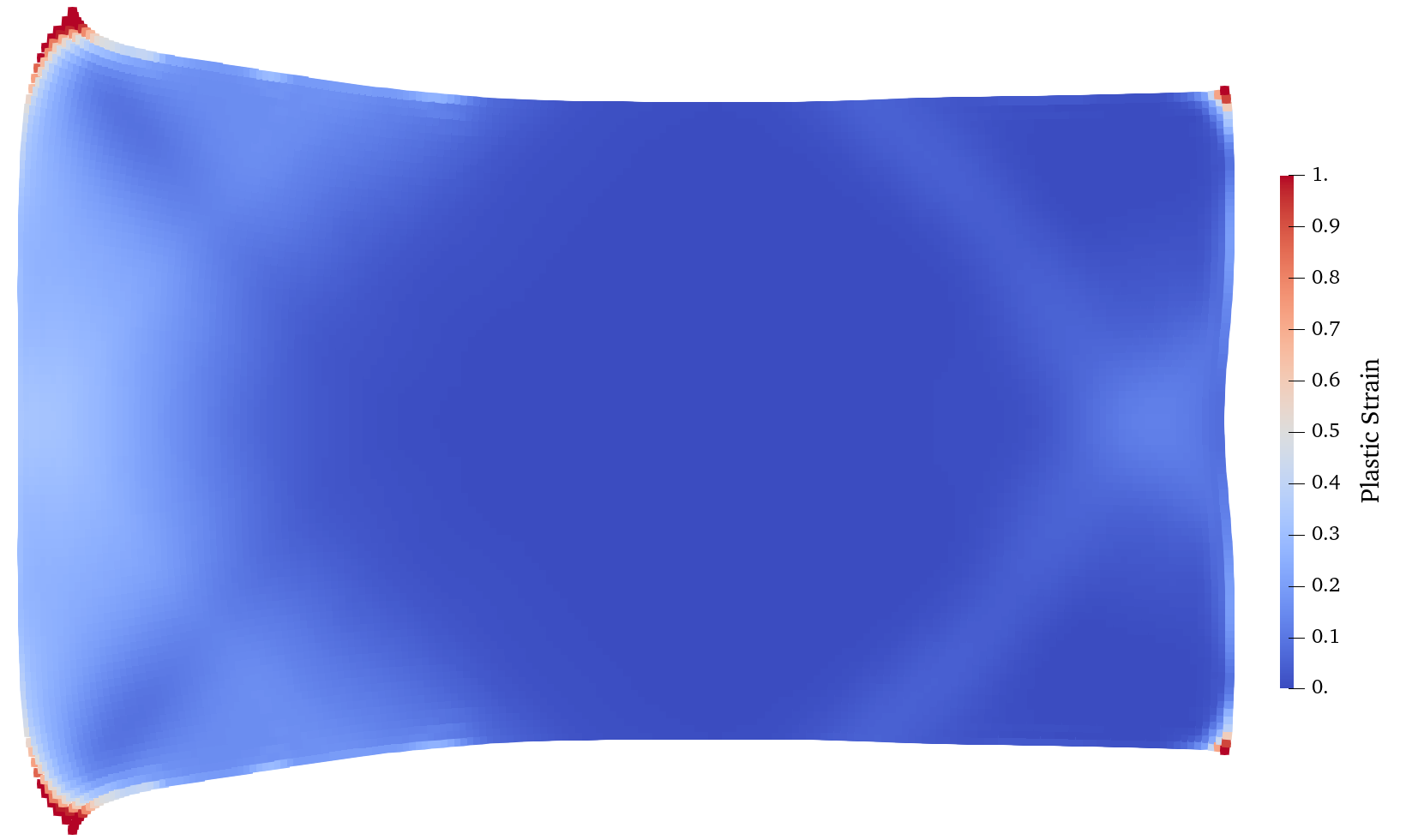}}
  \caption{Chamber detonation. Final configuration of the solid structure compared between the IGA solid and PD solid. Contours indicate the plastic strain levels.}
  \label{fig:detonation_eqps}
\end{figure*}

The bond-associative, quadratic ($n=2$) PD solid and the quadratic IGA solid are used in the computations. Air pressure at different time instants during the simulation is compared between the IGA-PD and reference simulations in \cref{fig:detonation_pressure} using the finest mesh.
The plastic strain contours in the final configuration of the specimen are shown in \cref{fig:detonation_eqps}. A good agreement between the results is obtained. Note the large deformation at the bar corners and the {\em mushrooming} at the left edge of the bar. Also note that the shock waves bounce off the right wall, impact the specimen, and leave a permanent indentation on the right edge of the bar.

\begin{figure*}[!hbpt]
  \centering
  \subfloat[][]{\includegraphics[width=0.32\textwidth,height=0.31\textwidth,trim={0cm 0cm 0cm 0cm},clip]{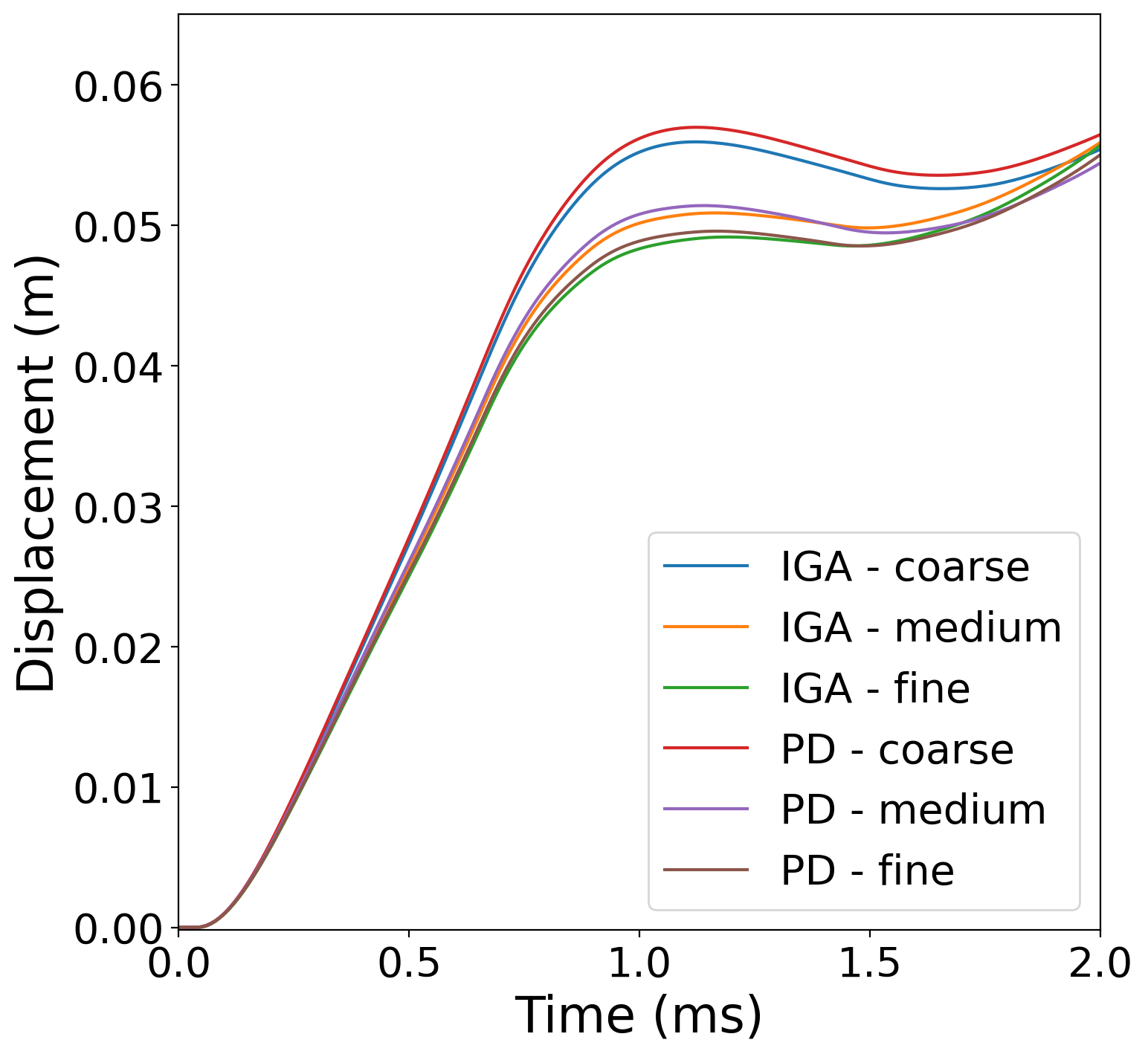}}
  \hspace{4pt}
  \subfloat[][]{\includegraphics[width=0.32\textwidth,height=0.31\textwidth,trim={0cm 0cm 0cm 0cm},clip]{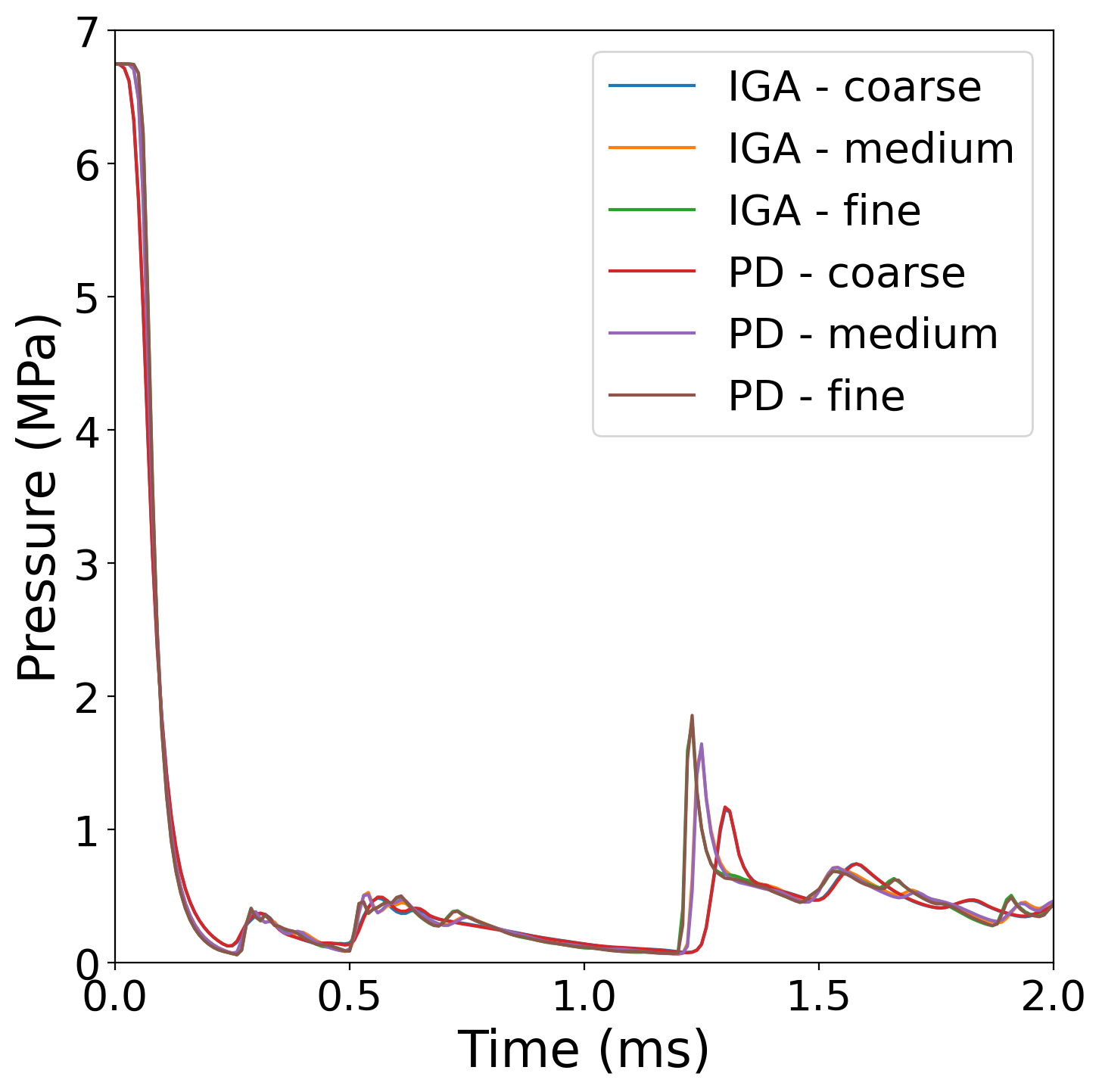}}
  \hspace{4pt}
  \subfloat[][]{\includegraphics[width=0.32\textwidth,height=0.31\textwidth,trim={0cm 0cm 0cm 0cm},clip]{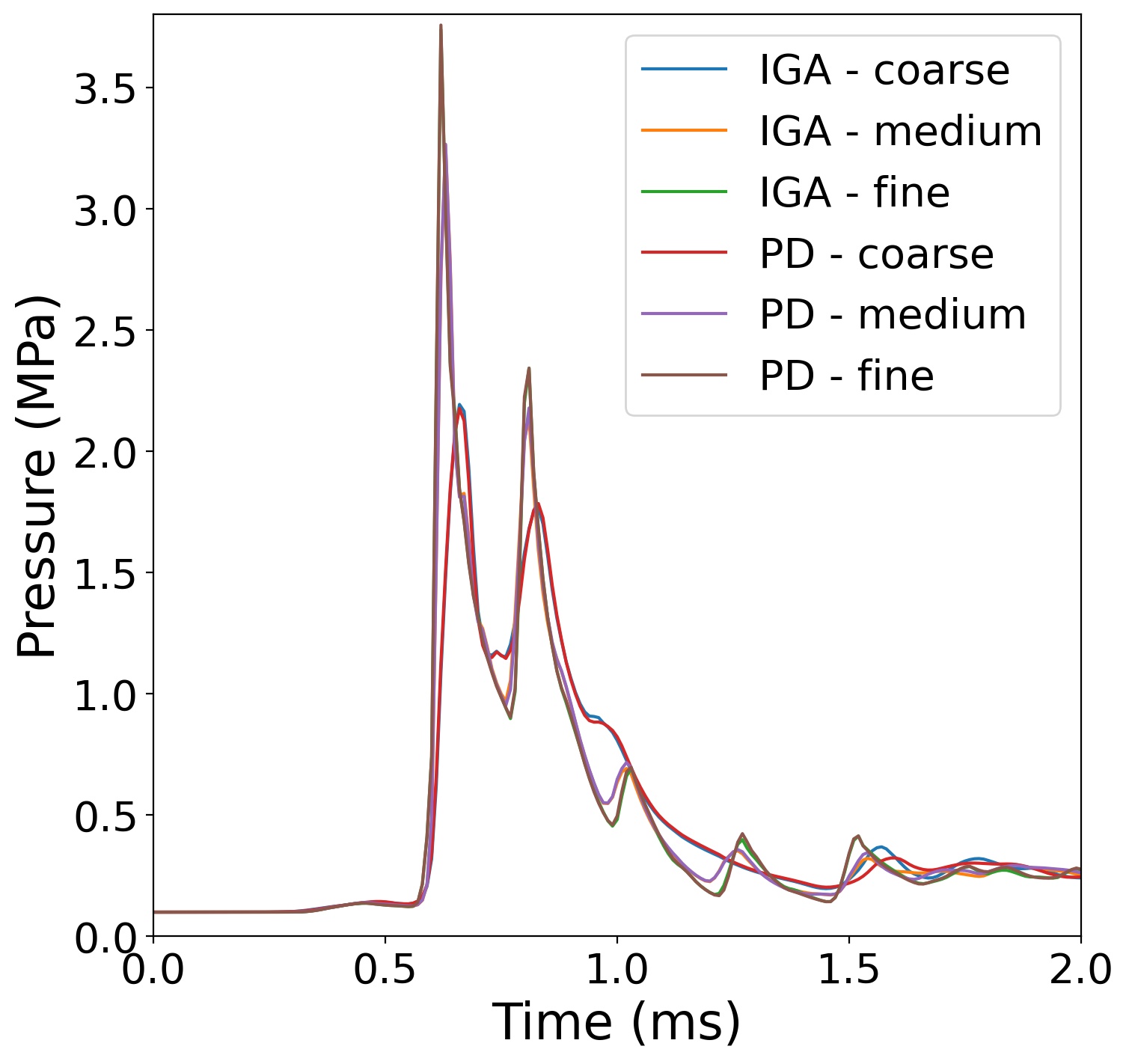}}
  \caption{Chamber detonation. Comparison of simulation results of different discretization levels of IGA solid and PD solid. The background mesh is refined at the same rate as the foreground discretization. (a): Horizontal displacement of the bar center of mass. (b): Pressure at the center of detonation. (c): Pressure at the center of the right wall.}
  \label{fig:detonation_result_mesh}
\end{figure*}

The discretization effects are studied in this problem. The time history of the center-of-mass displacement of the bar, pressure at the detonation center (Point 1 in \cref{fig:detonation_setup}), and pressure at the right-wall center (Point 2 in \cref{fig:detonation_setup}) are shown in \cref{fig:detonation_result_mesh} for different discretization levels of IGA-PD and IGA-IGA. The results indicate good convergence of the IGA-PD solution with mesh refinement, and a good agreement between the IGA-PD and IGA-IGA results.

\begin{figure*}[!hbpt]
  \centering
  \subfloat[][]{\includegraphics[width=0.32\textwidth,height=0.31\textwidth,trim={0cm 0cm 0cm 0cm},clip]{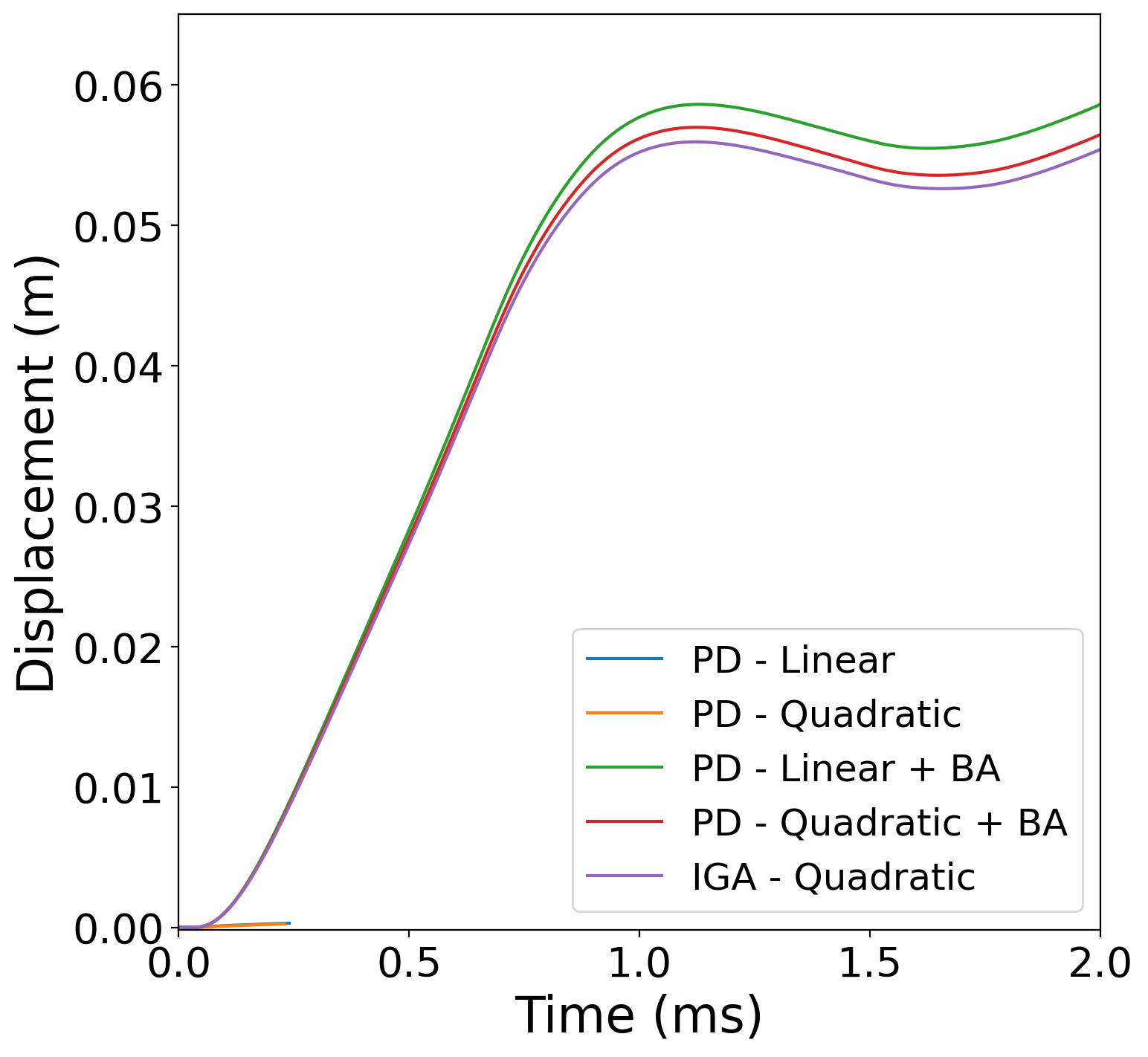}}
  \hspace{4pt}
  \subfloat[][]{\includegraphics[width=0.32\textwidth,height=0.31\textwidth,trim={0cm 0cm 0cm 0cm},clip]{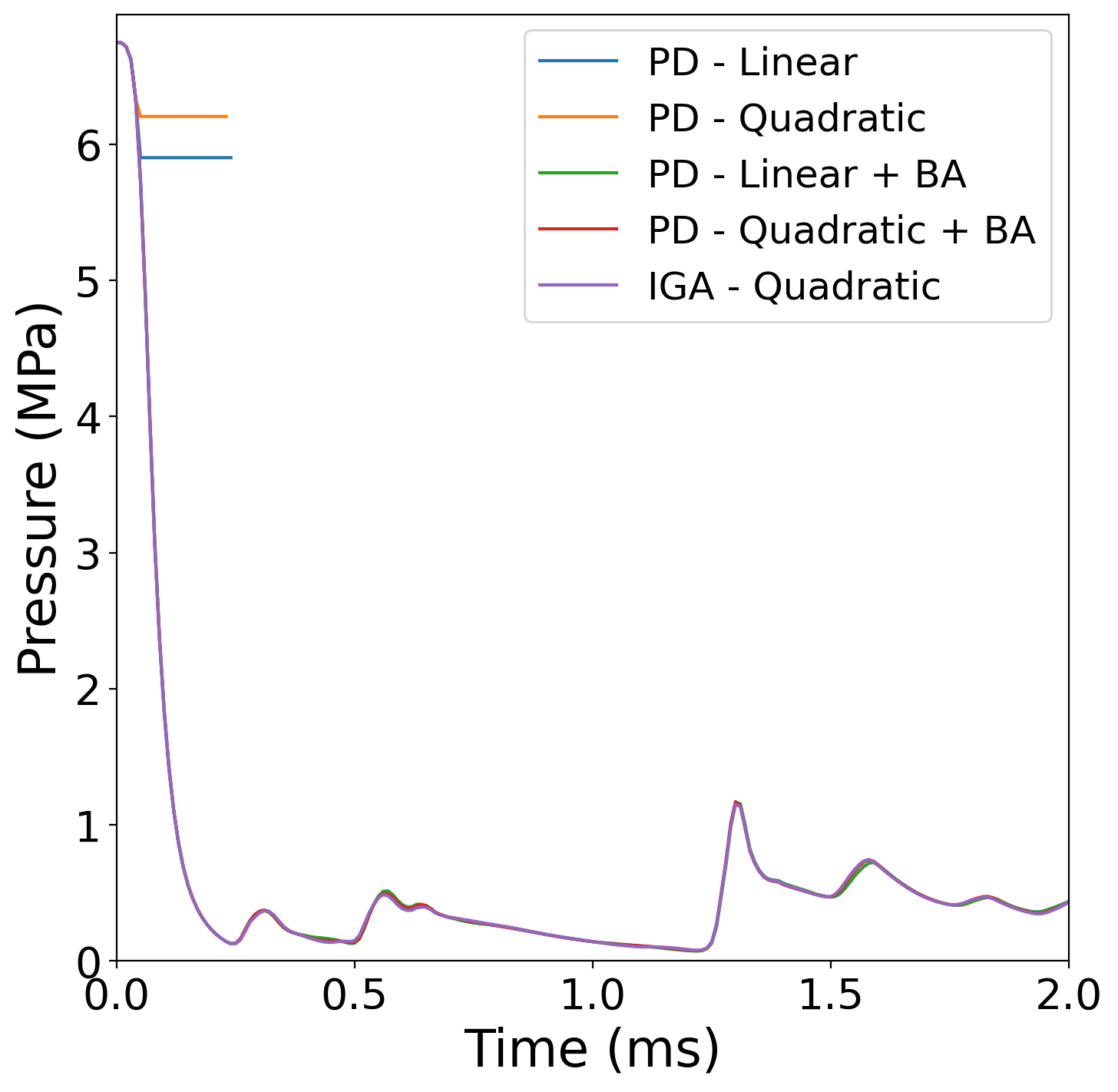}}
  \hspace{4pt}
  \subfloat[][]{\includegraphics[width=0.32\textwidth,height=0.31\textwidth,trim={0cm 0cm 0cm 0cm},clip]{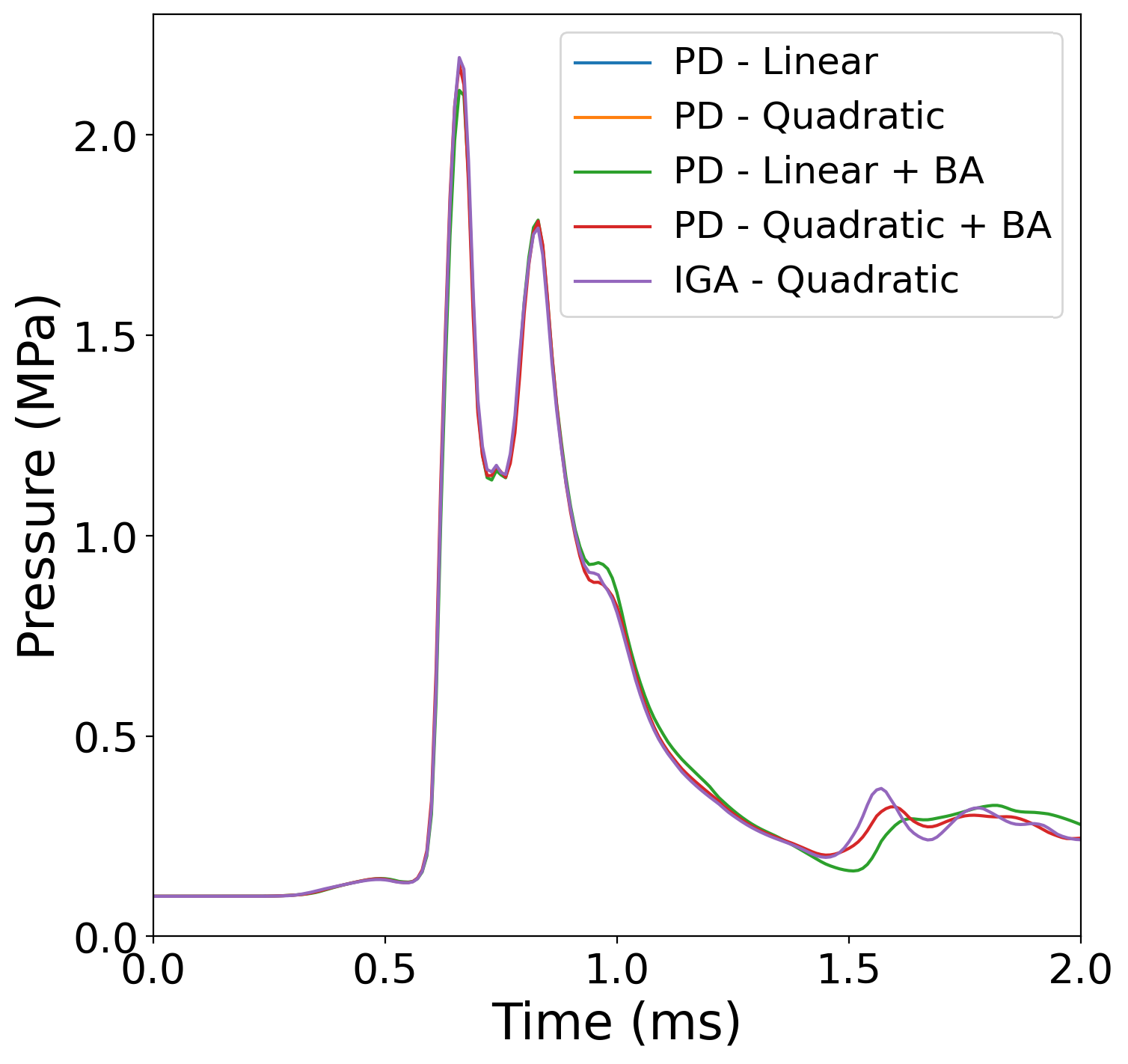}}
  \caption{Chamber detonation. Comparison of the simulation results involving IGA solid and different variants of PD solid (linear and quadratic gradient kernel functions, with and without bond-associative (BA) corrections). (a): Horizontal displacement of the bar center of mass. (b): Pressure at the center of detonation. (c): Pressure at the center of the right wall}
  \label{fig:detonation_result_variants}
\end{figure*}

We study the effects of the BA stabilization and order of accuracy on the PD solid behavior. The linear ($n=1$) and quadratic ($n=2$) PD variants with and without BA corrections are considered. The PD simulation results are compared with their quadratic IGA counterparts in \cref{fig:detonation_result_variants}. The non-BA PD results are of poor quality due to the presence of unstable modes that deteriorate the solution and eventually lead to divergence of the simulation (around $t=\SI{0.25}{ms}$ here). The results of quadratic PD model equipped with BA stabilization are in good agreement with the reference solution (better than the linear + BA version), which confirms the conclusions of~\cite{behzadinasab2021unifiedI} that the combination of BA stabilization and higher-order corrections provides the best quality results in correspondence-based PD. In the remainder of this paper, we utilize the quadratic + BA model for the PD calculations. 

\subsection{Brittle material subjected to internal explosion}
\label{sec:brittle}

\begin{figure*}[!hbpt]
  \centering
  \subfloat[][]{\includegraphics[width=0.7\textwidth,trim={0cm 0cm 0cm 0cm},clip]{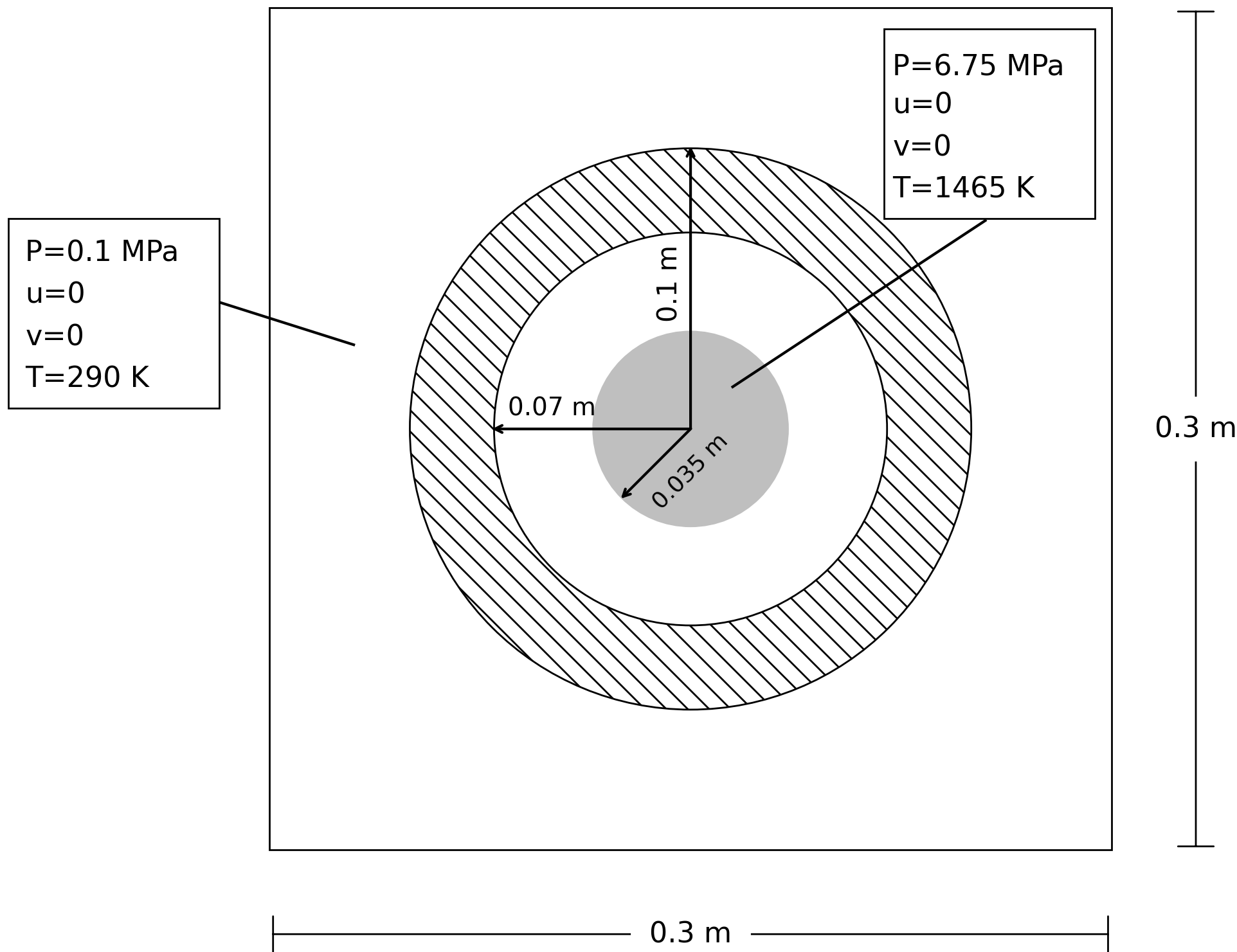}}
  \caption{Brittle fracture example. A hollow cylindrical block of elastic, brittle material subjected to internal detonation. Problem setup and geometry.}
  \label{fig:brittle_setup}
\end{figure*}

In this FSI example, we consider a hollow cylindrical block of elastic, brittle material subjected to internal detonation. As shown in \cref{fig:brittle_setup}, the detonation is initiated at the center of the hollow cylinder with the inner radius of $7 \, {\rm cm}$ and outer radius of $10 \, {\rm cm}$. The background domain is a $30 \, {\rm cm} \, \times \, 30 \, {\rm cm}$ square enclosing the cylinder domain. 
The solid elastic material has Young's modulus ${\rm E}=200 \, {\rm GPa}$, Poisson's ratio $\nu=0.29$, and initial density $\rho^s=7870 \, {\rm kg}/{\rm m}^3$. To simulate brittle fracture, the von Mises stress failure criterion is used with the critical stress $\sigma_{\rm cr} = 3 \, {\rm GPa}$. As noted in \cref{sec:failure}, a PD bond is broken once its associated von Mises stress exceeds $\sigma_{\rm cr}$. Initially, the air is at rest with $p = 0.1 \, {\rm MPa}$ and $T = 290 \, {\rm K}$. The detonation condition is enforced by setting the initial pressure $p = 6.75 \, {\rm MPa}$ and temperature $T = 1465 \, {\rm K}$ in a circular area with a radius of $3.5 \, {\rm cm}$, centered inside the hollow cylinder. 

In this problem, the background domain and foreground solid are discretized into a uniform mesh and semi-uniform nodal setting (uniform along the $\theta$-direction), respectively. We consider four discretization levels in this problem, D1--D4, with solid node spacings of $2.5 \, {\rm mm}$, $2 \, {\rm mm}$, $1.5 \, {\rm mm}$, and $1 \, {\rm mm}$, respectively. The fluid mesh size is set to three times of the solid node spacing in each case. The time step used for D1--D4 are $\SI{2.5}{\micro s}$, $\SI{2}{\micro s}$, $\SI{1.5}{\micro s}$, and $\SI{1}{\micro s}$, respectively. 

\begin{figure*}[!hbpt]
  \centering
  \subfloat{\includegraphics[width=0.6\textwidth,trim={0cm 1cm 0cm 3cm},clip]{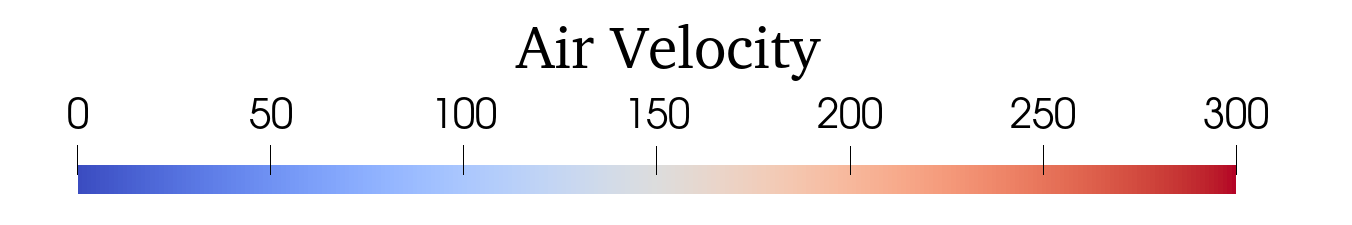}}

  \setcounter{subfigure}{0}
  \subfloat[][\SI{15}{\micro s}]{\includegraphics[width=0.32\textwidth,trim={0cm 0cm 0cm 0cm},clip]{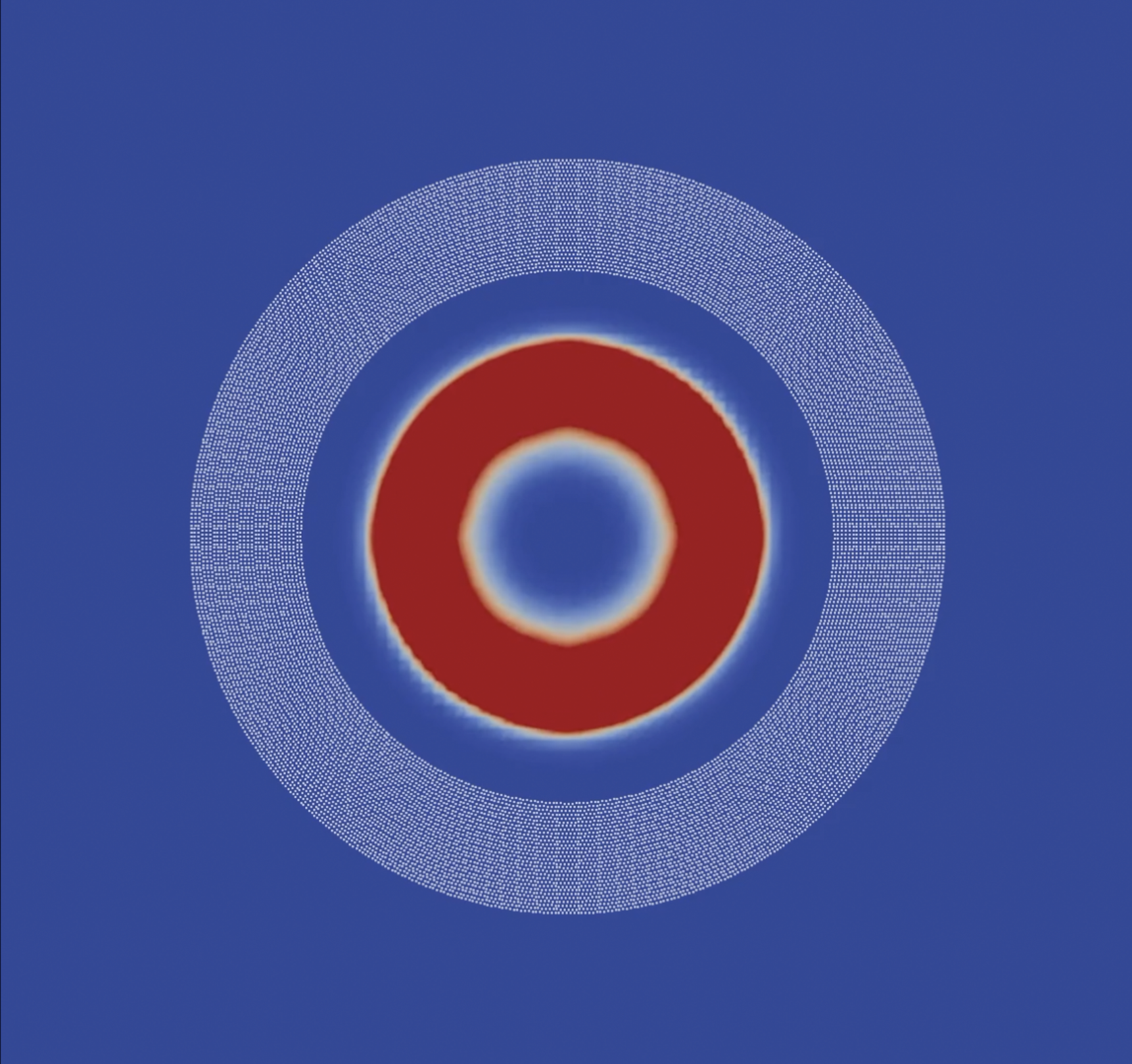}}
  \hspace{4pt}
  \subfloat[][\SI{30}{\micro s}]{\includegraphics[width=0.32\textwidth,trim={0cm 0cm 0cm 0cm},clip]{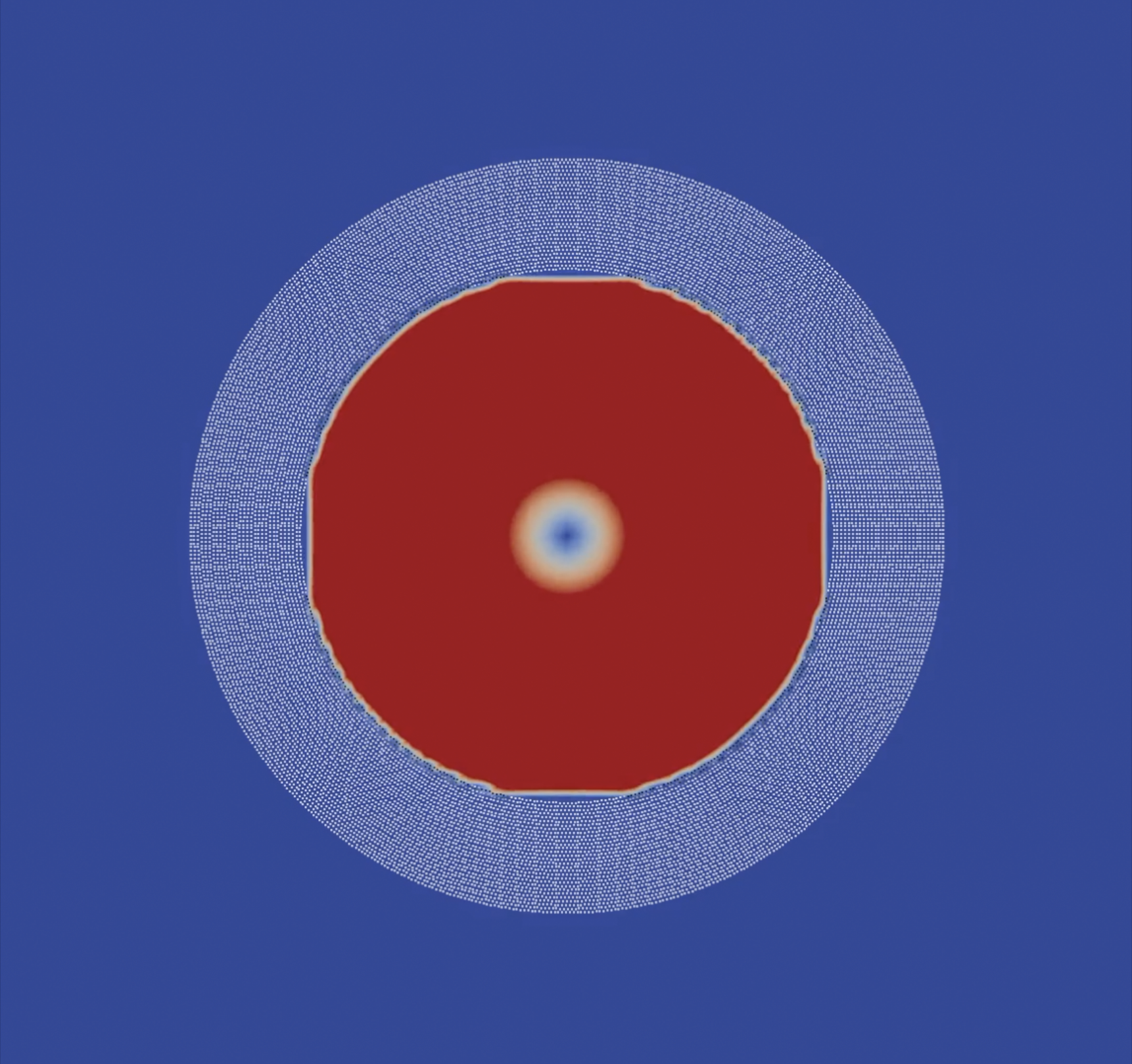}}
  \hspace{4pt}
  \subfloat[][\SI{70}{\micro s}]{\includegraphics[width=0.32\textwidth,trim={0cm 0cm 0cm 0cm},clip]{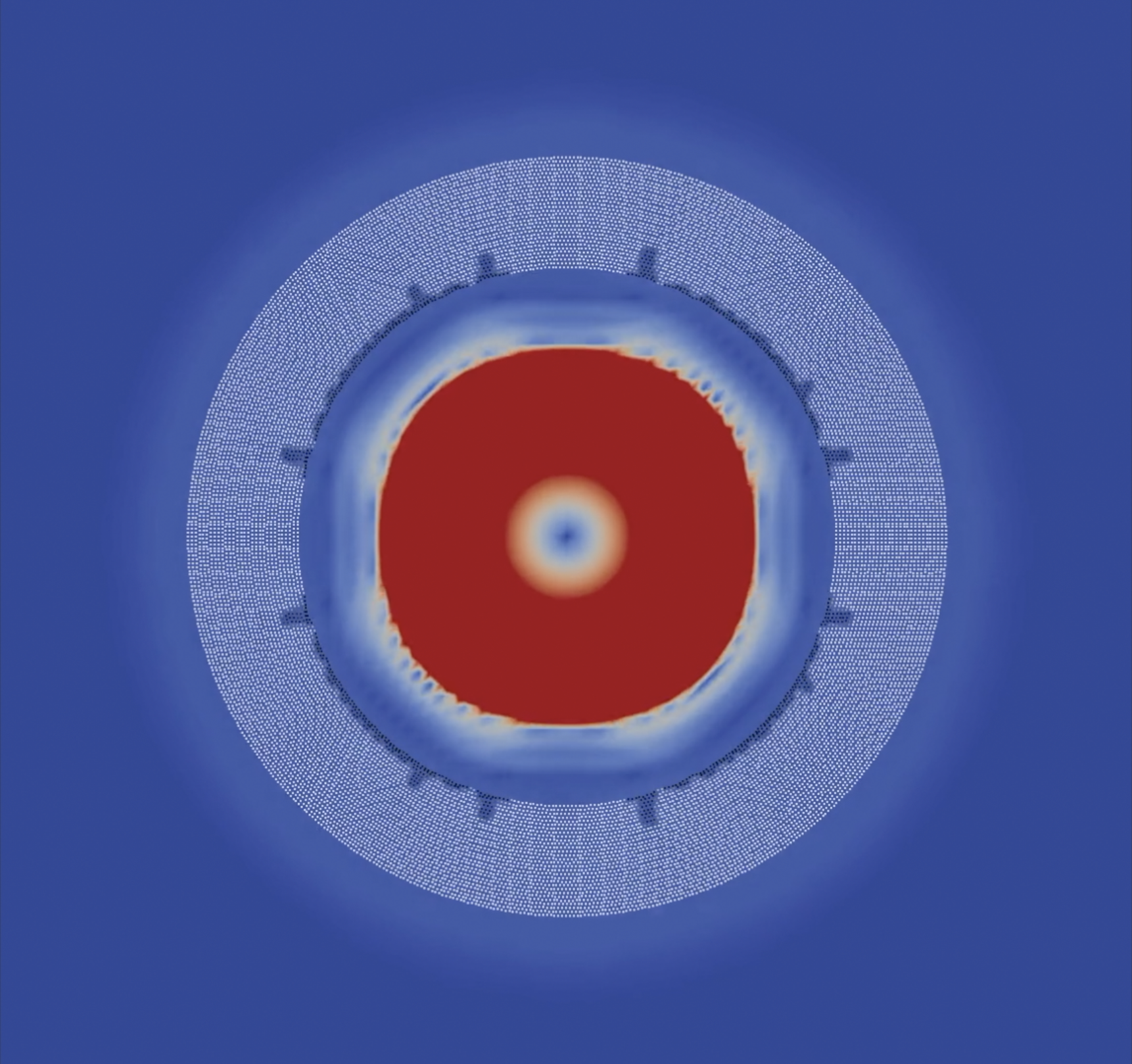}}

  \subfloat[][\SI{175}{\micro s}]{\includegraphics[width=0.32\textwidth,trim={0cm 0cm 0cm 0cm},clip]{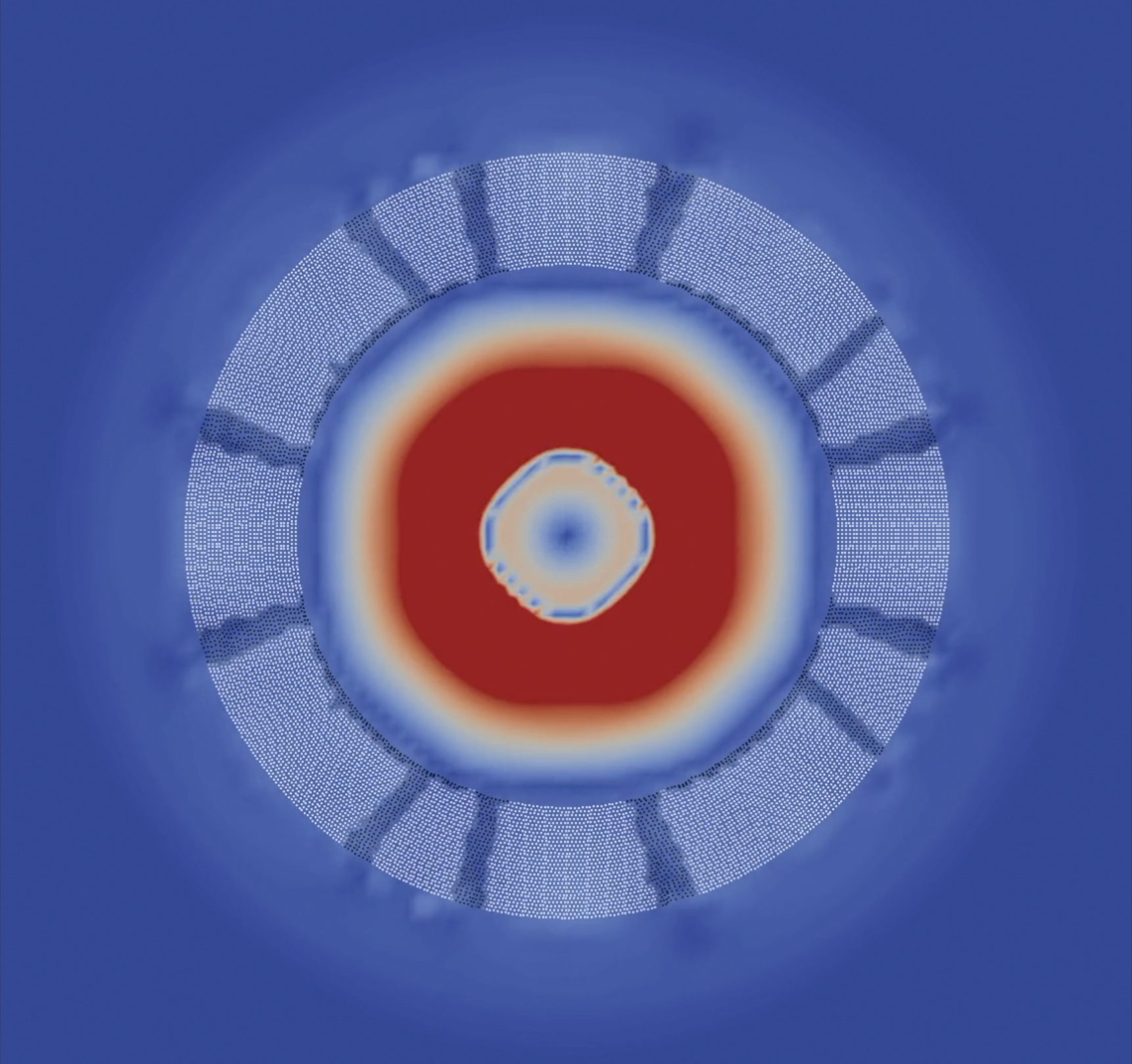}}
  \hspace{4pt}
  \subfloat[][\SI{150}{\micro s}]{\includegraphics[width=0.32\textwidth,trim={0cm 0cm 0cm 0cm},clip]{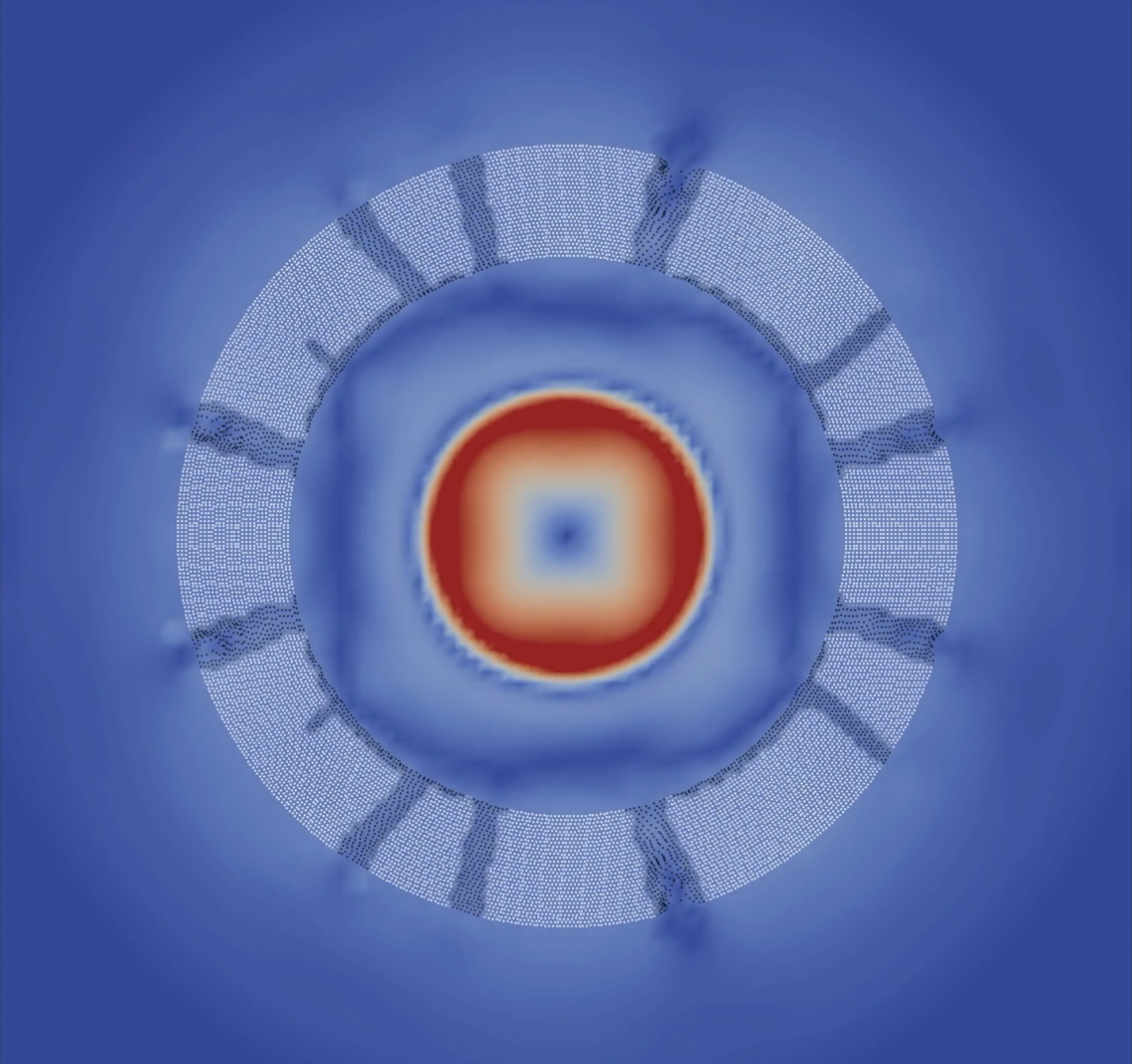}}
  \hspace{4pt}
  \subfloat[][\SI{300}{\micro s}]{\includegraphics[width=0.32\textwidth,trim={0cm 0cm 0cm 0cm},clip]{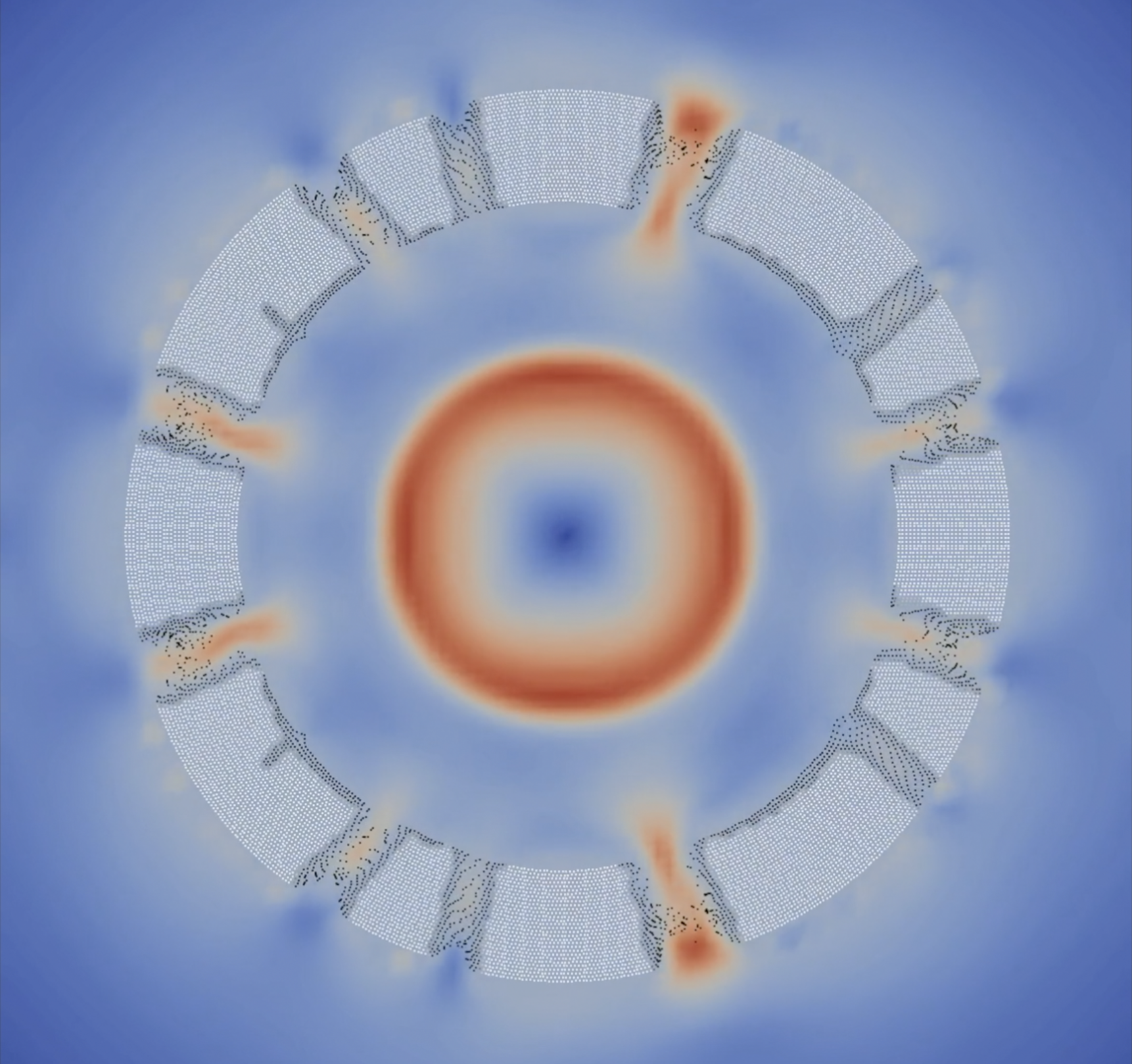}}
  \caption{Brittle fracture example. Snapshots of air speed (in m/s) and solid damage in the current configuration at different stages during the simulation.}
  \label{fig:brittle_velocity}
\end{figure*}

\begin{figure*}[!hbpt]
  \centering
  \subfloat[][D1]{\includegraphics[width=0.24\textwidth,trim={0cm 0cm 0cm 0cm},clip]{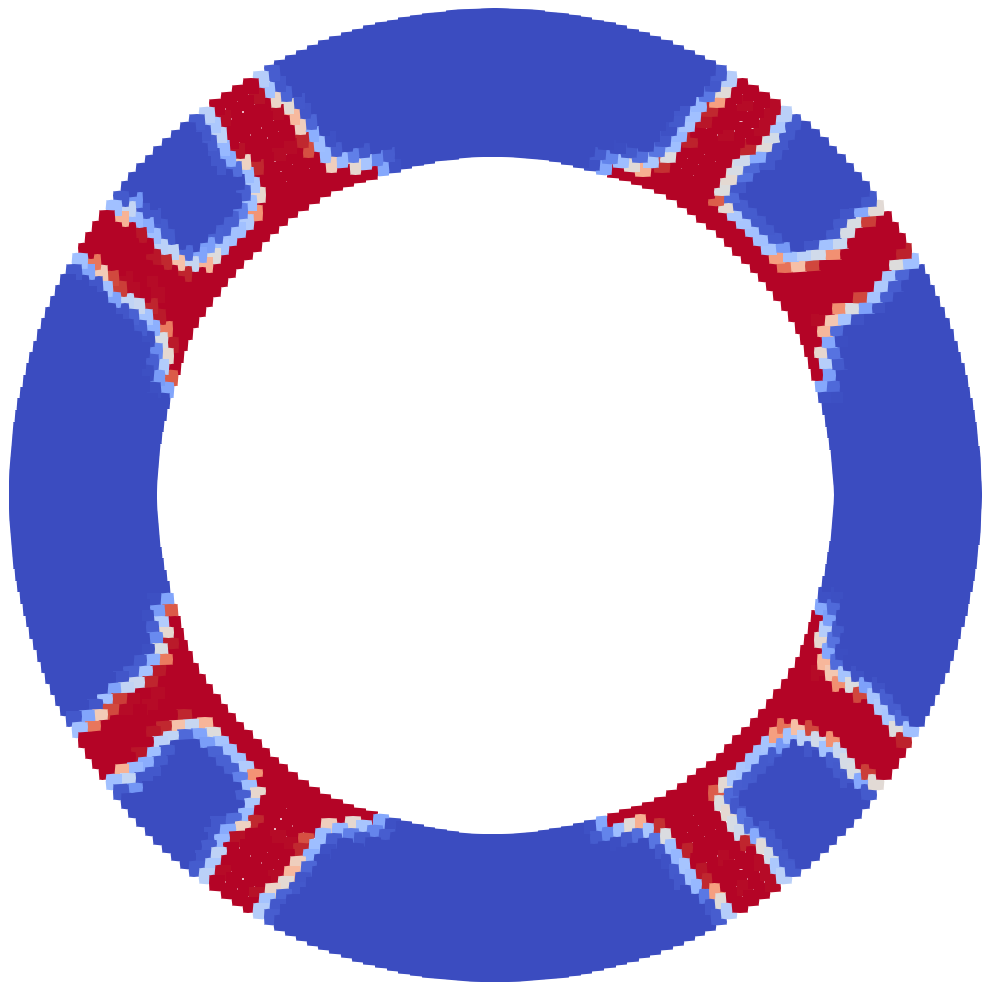}}
  \hspace{1pt}
  \subfloat[][D2]{\includegraphics[width=0.24\textwidth,trim={0cm 0cm 0cm 0cm},clip]{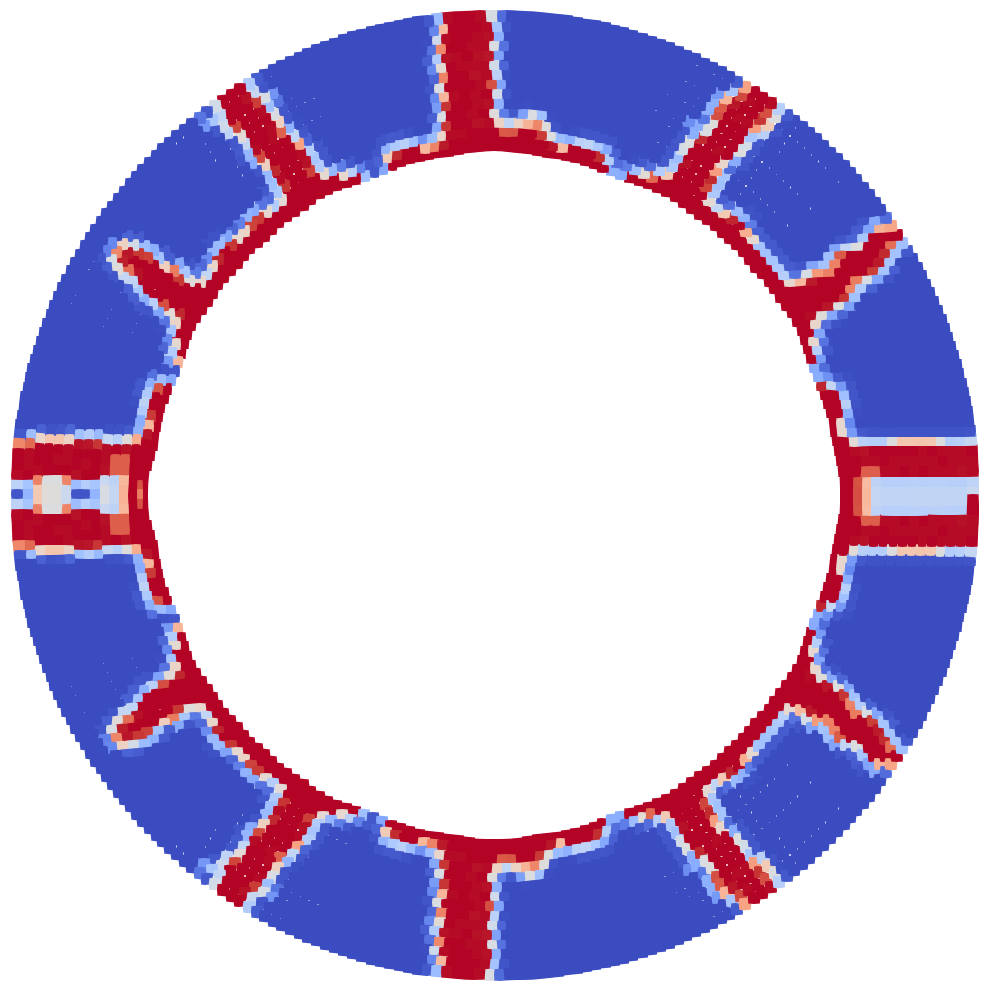}}
  \hspace{1pt}
  \subfloat[][D3]{\includegraphics[width=0.24\textwidth,trim={0cm 0cm 0cm 0cm},clip]{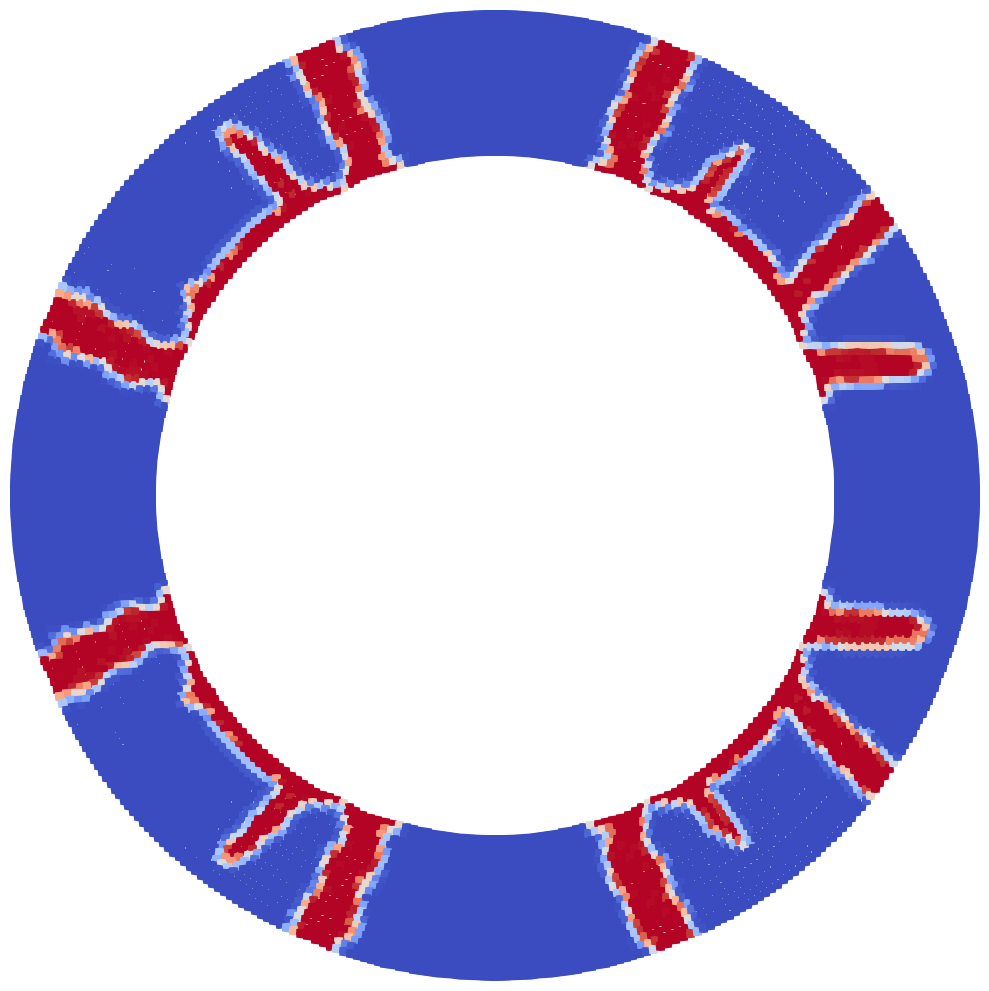}}
  \hspace{1pt}
  \subfloat[][D4]{\includegraphics[width=0.24\textwidth,trim={0cm 0cm 0cm 0cm},clip]{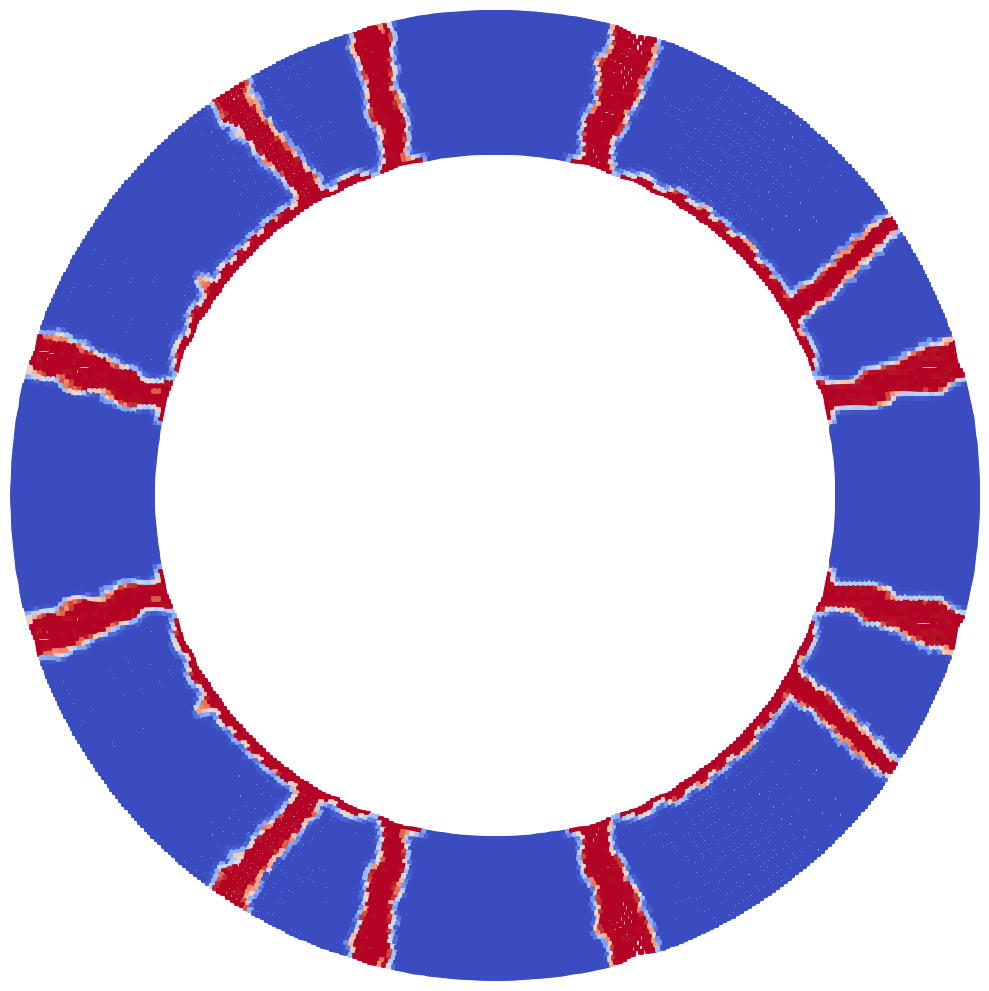}}
  \caption{Brittle fracture example. Solid damage contours in different discretization levels D1--D4 at $t = \SI{100}{\micro s}$. The fracture patterns are qualitatively similar for the different cases here.}
  \label{fig:brittle_fracture}
\end{figure*}

The results for D4 is shown in \cref{fig:brittle_velocity}, where the air speed is depicted on the background grid while the damage field is shown on the PD nodes. 
As expected, the impact of shock waves with the inner side of cylinder results in nucleation of several cracks at different locations along that side. Following the propagation of cracks through the solid structure, it is completely fractured and split into multiple fragments. The pressurized air passes through the empty space between fragments in the final snapshot in the figure. Because of the nature of the immersed methodology, fully-damaged solid nodes could still move within the computational domain without triggering mesh distortion issues that typically occur in a finite element foreground discretization. Although we do not present a comparison to experimental or other computational results, the predicted qualitative behavior appears to be physically reasonable. As shown in \cref{fig:brittle_fracture}, the fractures on the final configuration of the solid appear to converge toward a given state.

\subsection{Ductile material subjected to internal explosion}
\label{sec:ductile}

\begin{figure*}[!hbpt]
  \centering
  \subfloat[][]{\includegraphics[width=0.7\textwidth,trim={0cm 0cm 0cm 0cm},clip]{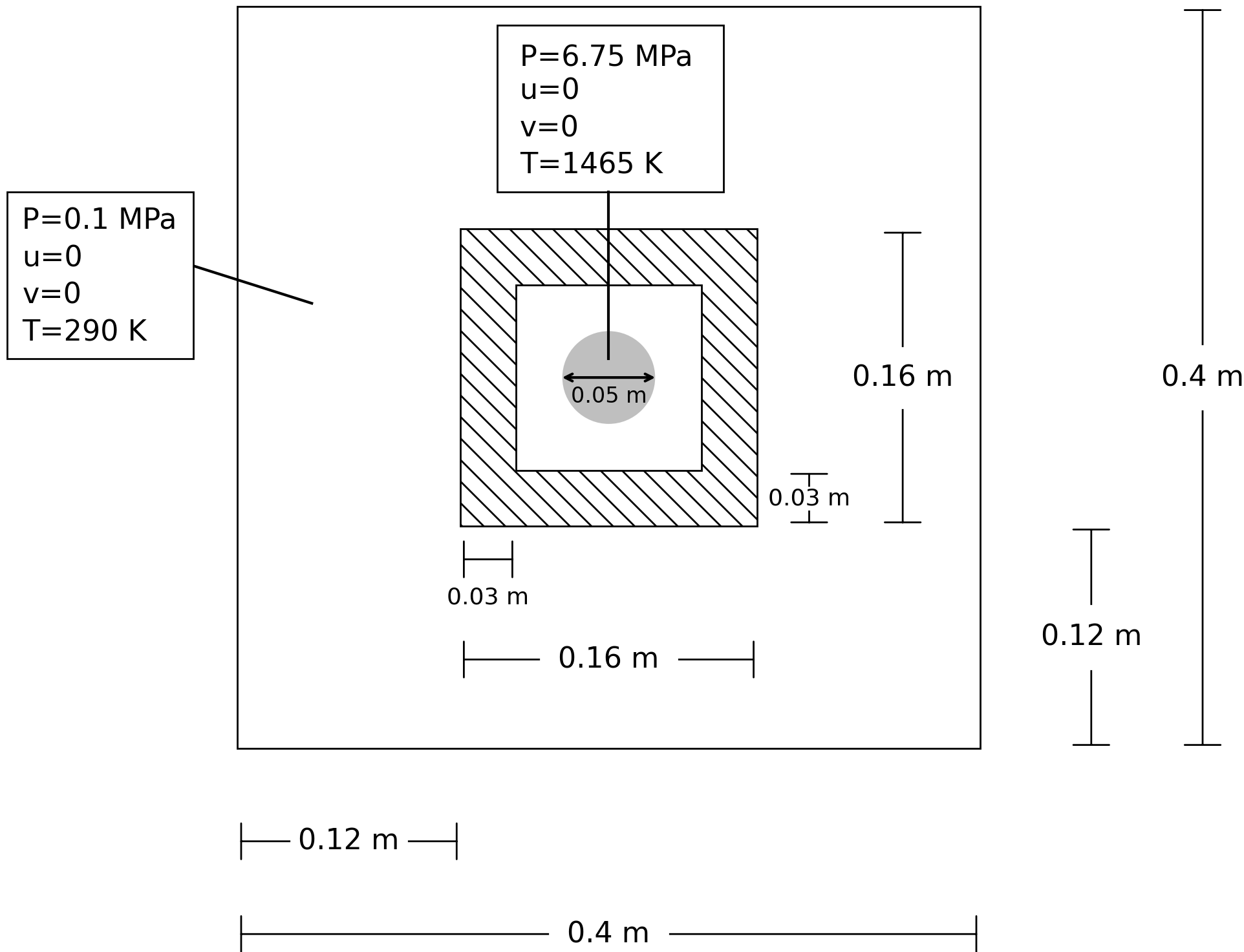}}
  \caption{Ductile fracture problem. A hollow square block of ductile material subjected to internal detonation. Problem setup and geometry.}
  \label{fig:ductile_setup}
\end{figure*}
The problem consists of a hollow square block of ductile material subjected to internal explosion. As shown in \cref{fig:ductile_setup}, the detonation is initiated at the center of the hollow square with inner dimension of $10 \, {\rm cm}$ and outer dimension of $16 \, {\rm cm}$. The background domain has the dimensions of $40 \, {\rm cm} \, \times \, 40 \, {\rm cm}$. Isotropic linear hardening rule is used for the solid ductile material, which has Young's modulus ${\rm E}=200 \, {\rm GPa}$, Poisson's ratio $\nu=0.29$, yield stress $\sigma_Y=0.4 \, {\rm GPa}$, hardening modulus $H = 0.1 \, {\rm GPa}$, and initial density $\rho=7870 \, {\rm kg}/{\rm m}^3$. To simulate ductile fracture, a plasticity-driven failure approach given in \cref{eqn:plastic_strain} is used with $\bar{\epsilon}^P_{\rm th} = 0.18$ and $\bar{\epsilon}^P_{\rm cr} = 0.2$. The air is at rest initially with $p = 0.1 \, {\rm MPa}$ and $T = 290 \, {\rm K}$. The detonation condition is initiated by setting the pressure $p = 6.75 \, {\rm MPa}$ and temperature $T = 1465 \, {\rm K}$ in a circular area with a radius of $5 \, {\rm cm}$, centered inside the hollow square. 

In this example, the background domain and foreground solid are discretized uniformly. Four discretization levels are considered here, D1--D4, with the solid node spacings of $2.5 \, {\rm mm}$, $2 \, {\rm mm}$, $1.5 \, {\rm mm}$, and $1 \, {\rm mm}$, respectively. In each case, the fluid mesh size is set to four times that of the solid node spacing. The time step used for the D1--D4 meshes is $\SI{1.8}{\micro s}$, $\SI{1.2}{\micro s}$, $\SI{0.9}{\micro s}$, and $\SI{0.6}{\micro s}$, respectively. 

\begin{figure*}[!hbpt]
  \centering
  \subfloat{\includegraphics[width=0.6\textwidth,trim={0cm 1cm 0cm 3cm},clip]{ductile_air_vel_scale.png}}

  \setcounter{subfigure}{0}
  \subfloat[][\SI{25}{\micro s}]{\includegraphics[width=0.32\textwidth,trim={0cm 0cm 0cm 0cm},clip]{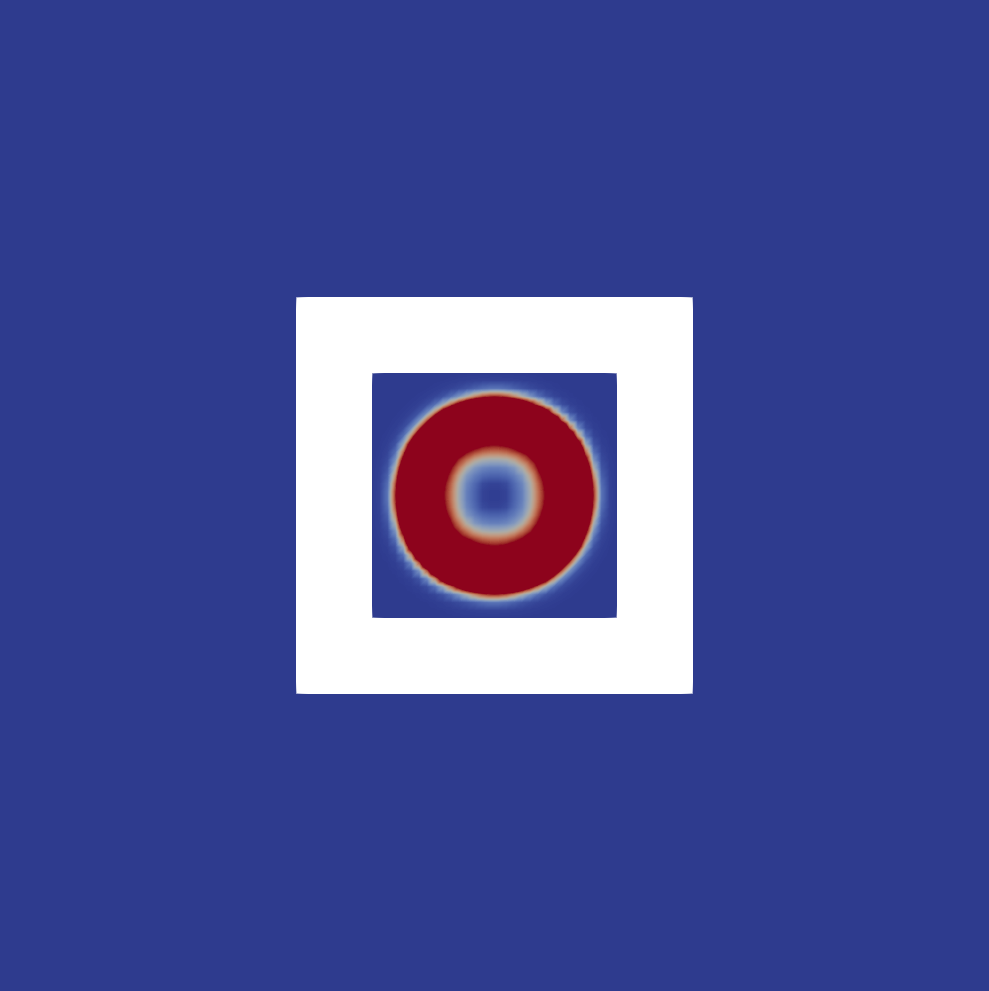}}
  \hspace{4pt}
  \subfloat[][\SI{75}{\micro s}]{\includegraphics[width=0.32\textwidth,trim={0cm 0cm 0cm 0cm},clip]{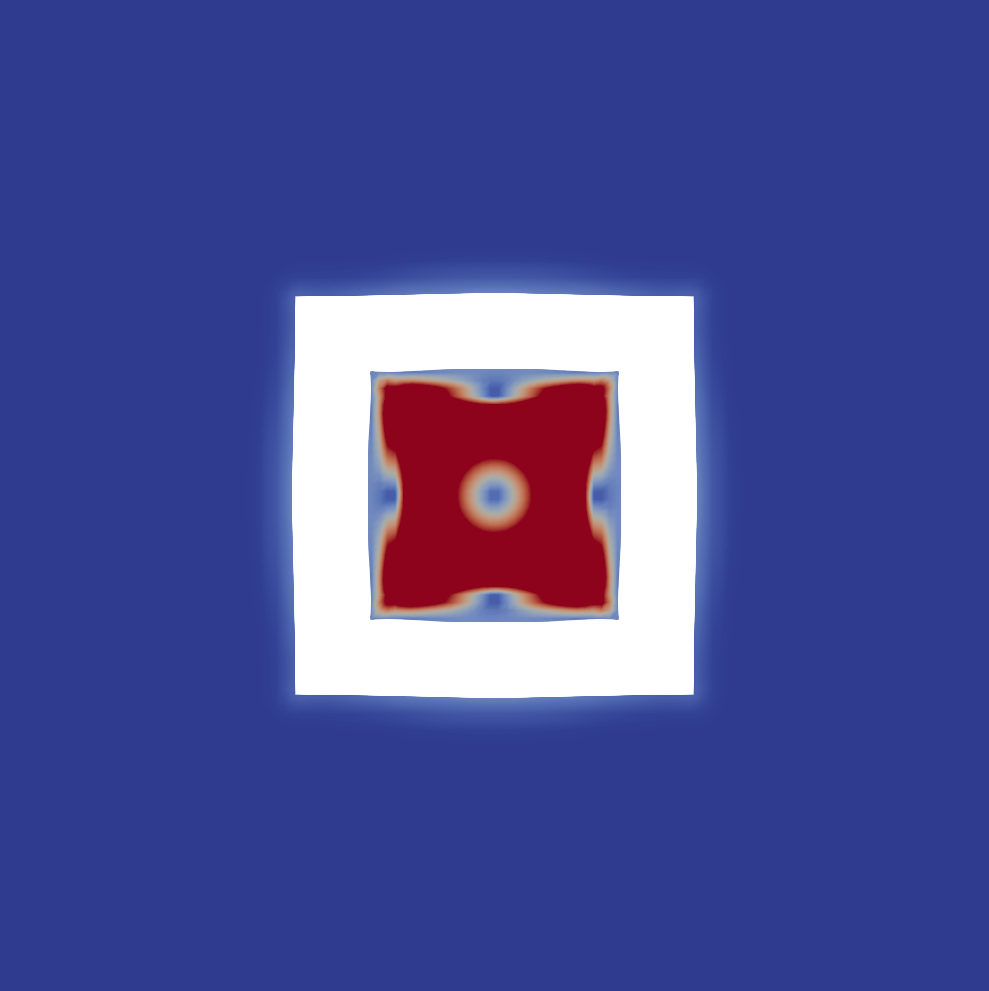}}
  \hspace{4pt}
  \subfloat[][\SI{125}{\micro s}]{\includegraphics[width=0.32\textwidth,trim={0cm 0cm 0cm 0cm},clip]{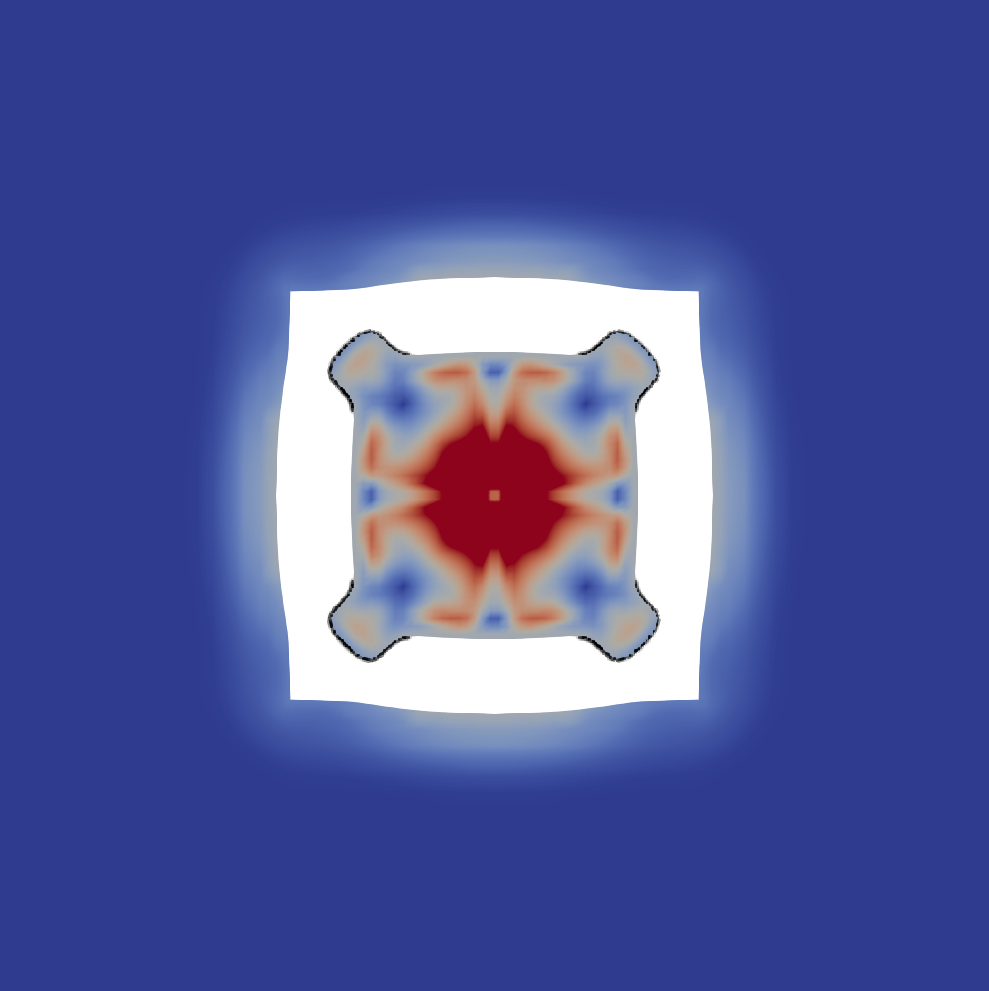}}

  \subfloat[][\SI{175}{\micro s}]{\includegraphics[width=0.32\textwidth,trim={0cm 0cm 0cm 0cm},clip]{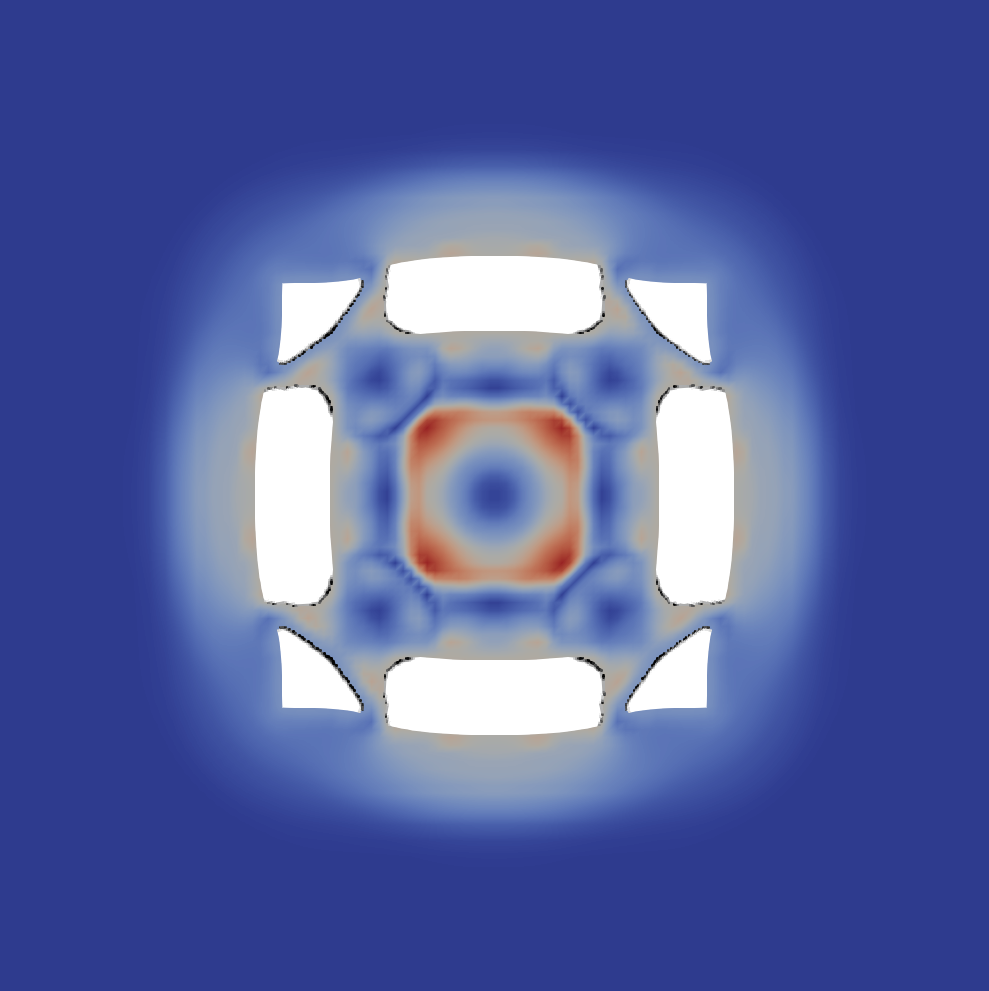}}
  \hspace{4pt}
  \subfloat[][\SI{275}{\micro s}]{\includegraphics[width=0.32\textwidth,trim={0cm 0cm 0cm 0cm},clip]{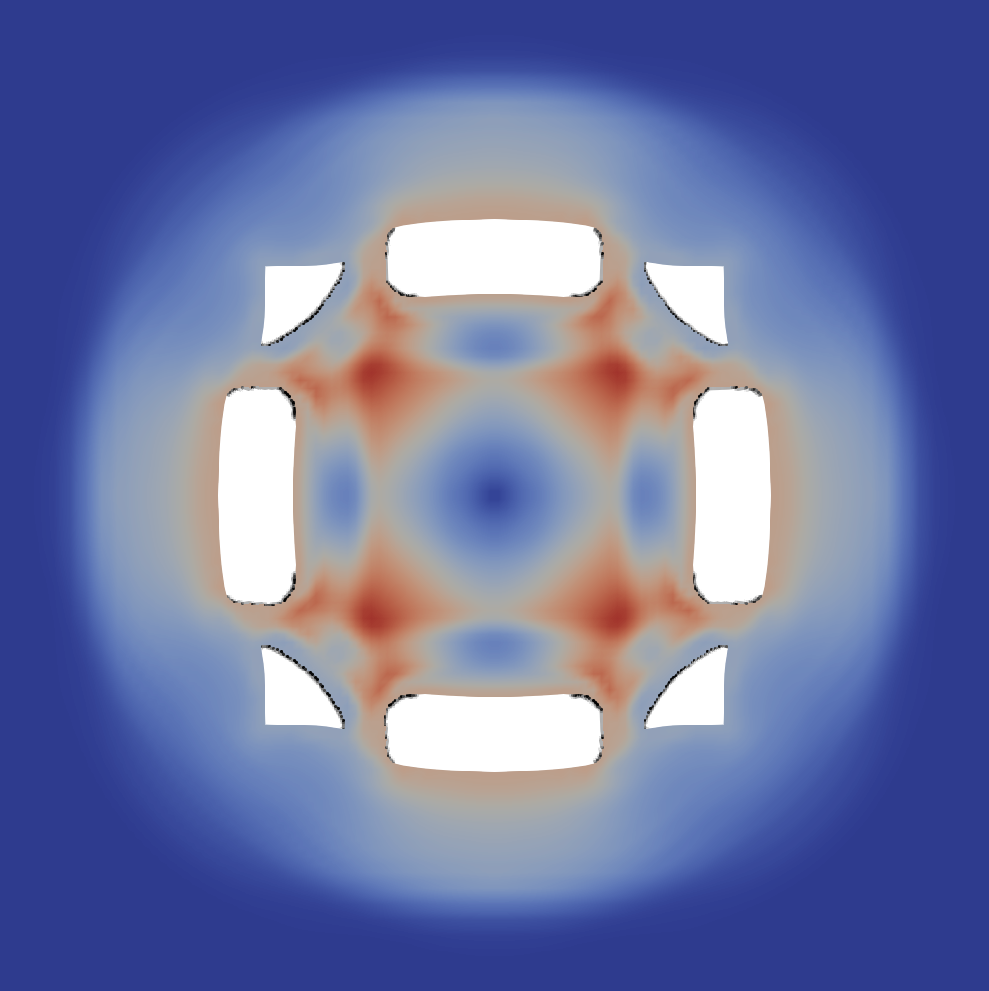}}
  \hspace{4pt}
  \subfloat[][\SI{375}{\micro s}]{\includegraphics[width=0.32\textwidth,trim={0cm 0cm 0cm 0cm},clip]{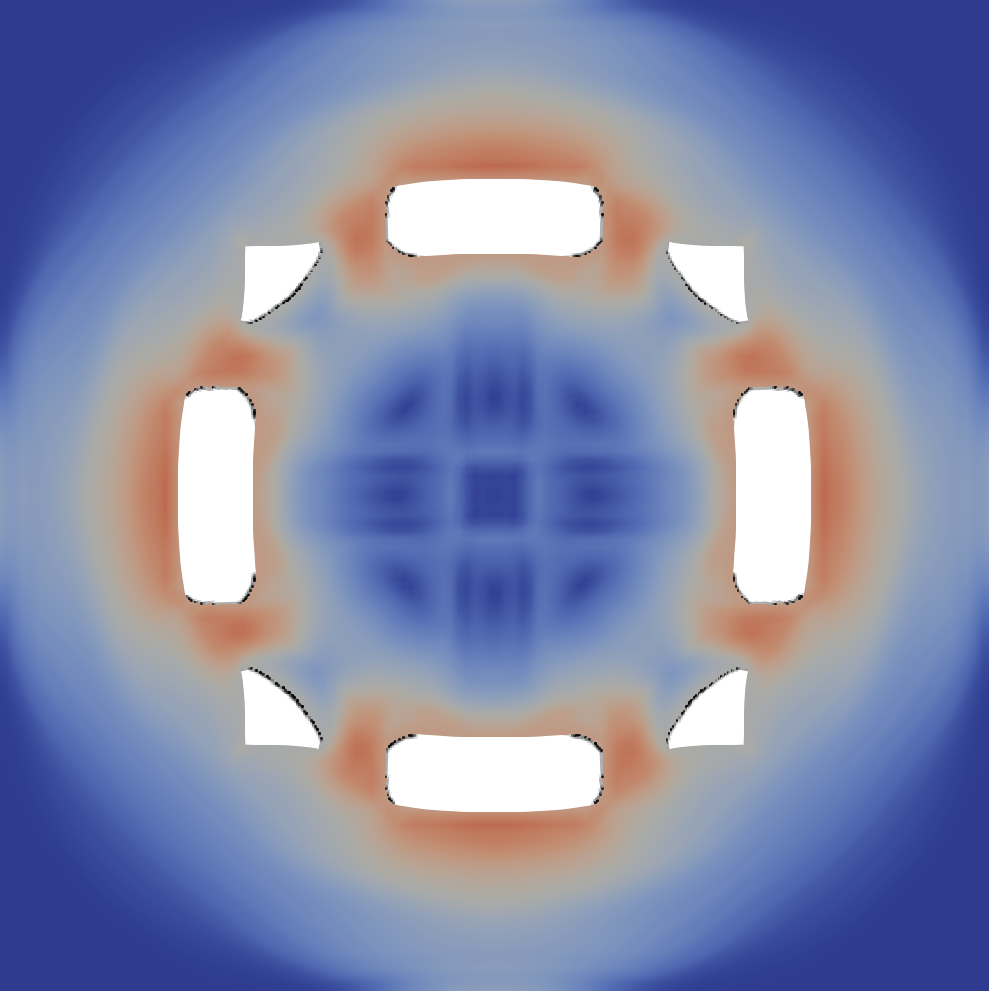}}
  \caption{Ductile fracture problem. Snapshots of air speed (in m/s) and solid damage in the current configuration at different stages during the simulation. For visualization purposes, the fully broken solid nodes are filtered out.}
  \label{fig:ductile_velocity}
\end{figure*}

\begin{figure*}[!hbpt]
  \centering
  \subfloat[][D1]{\includegraphics[width=0.24\textwidth,trim={0cm 0cm 0cm 0cm},clip]{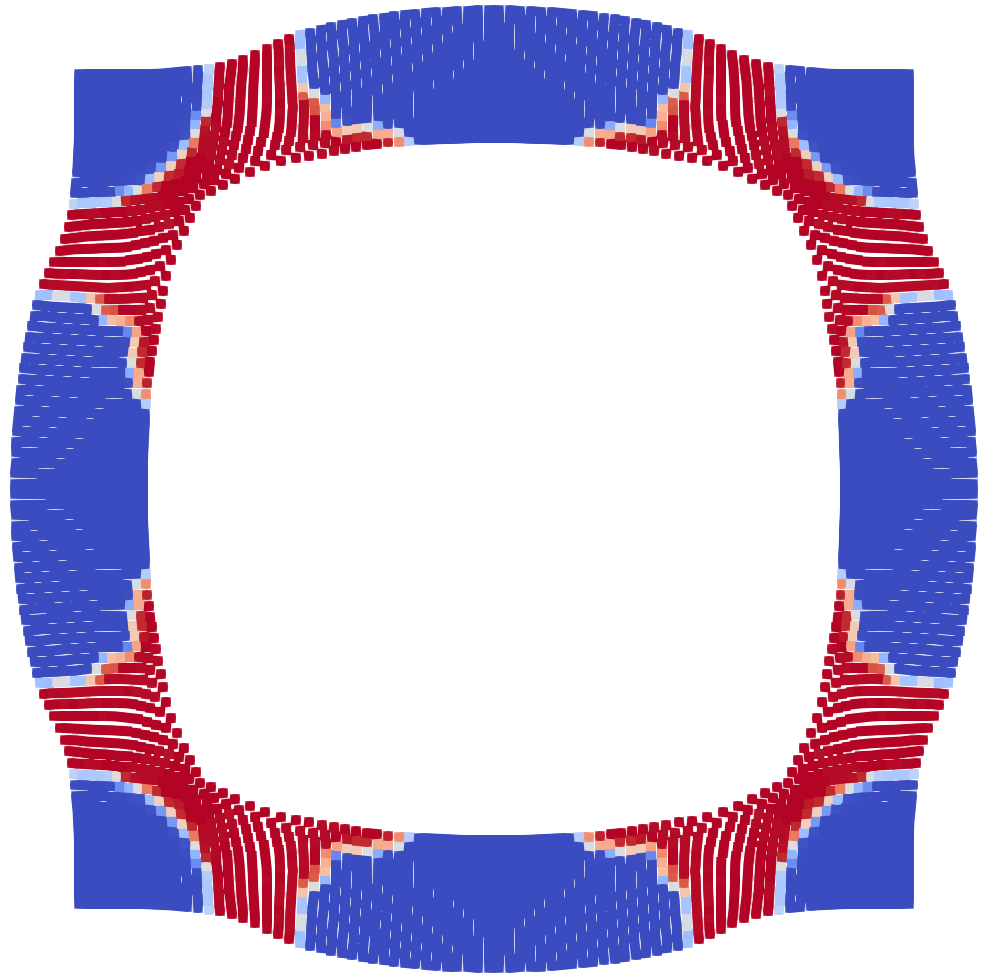}}
  \hspace{1pt}
  \subfloat[][D2]{\includegraphics[width=0.24\textwidth,trim={0cm 0cm 0cm 0cm},clip]{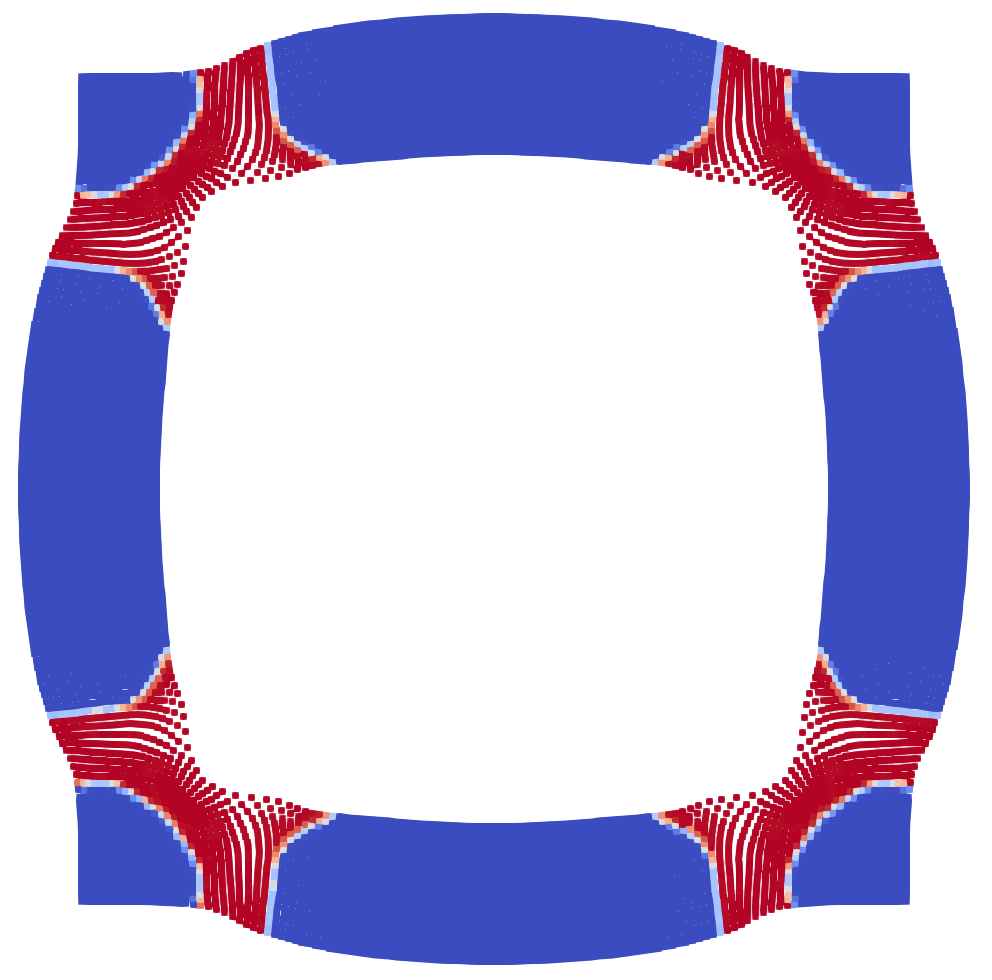}}
  \hspace{1pt}
  \subfloat[][D3]{\includegraphics[width=0.24\textwidth,trim={0cm 0cm 0cm 0cm},clip]{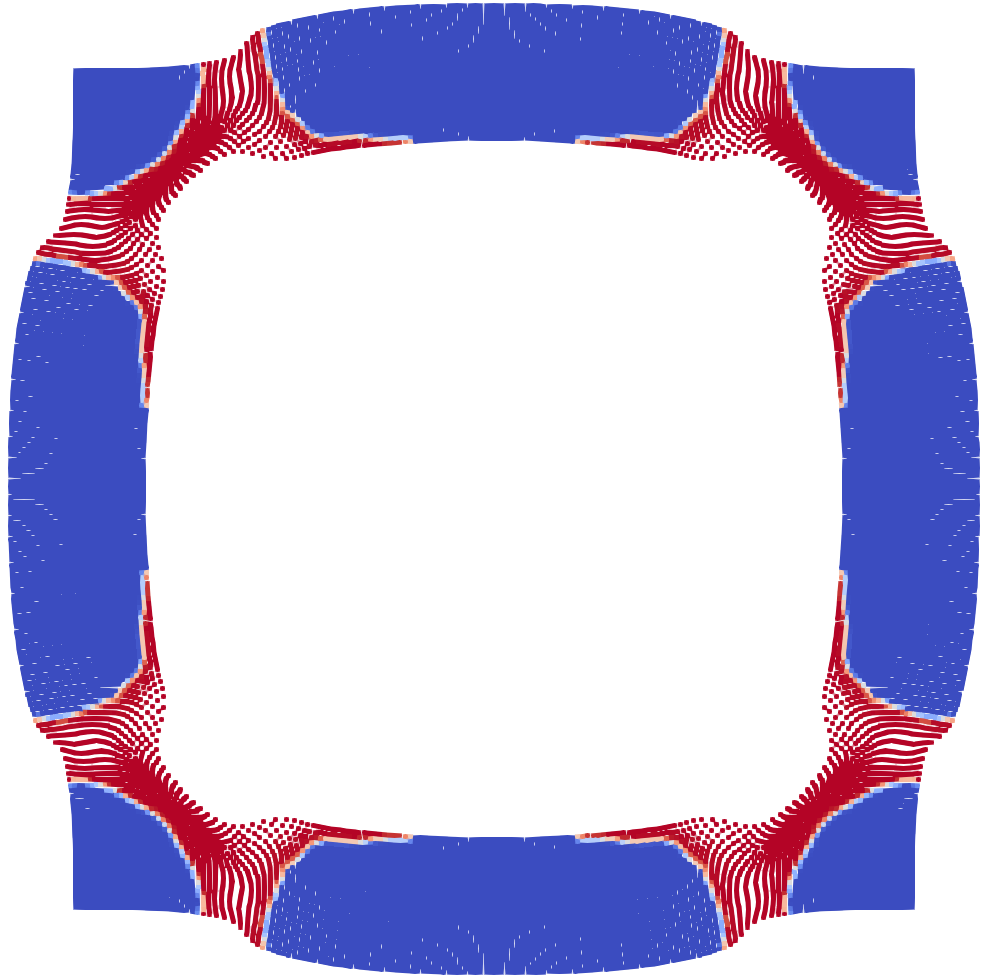}}
  \hspace{1pt}
  \subfloat[][D4]{\includegraphics[width=0.24\textwidth,trim={0cm 0cm 0cm 0cm},clip]{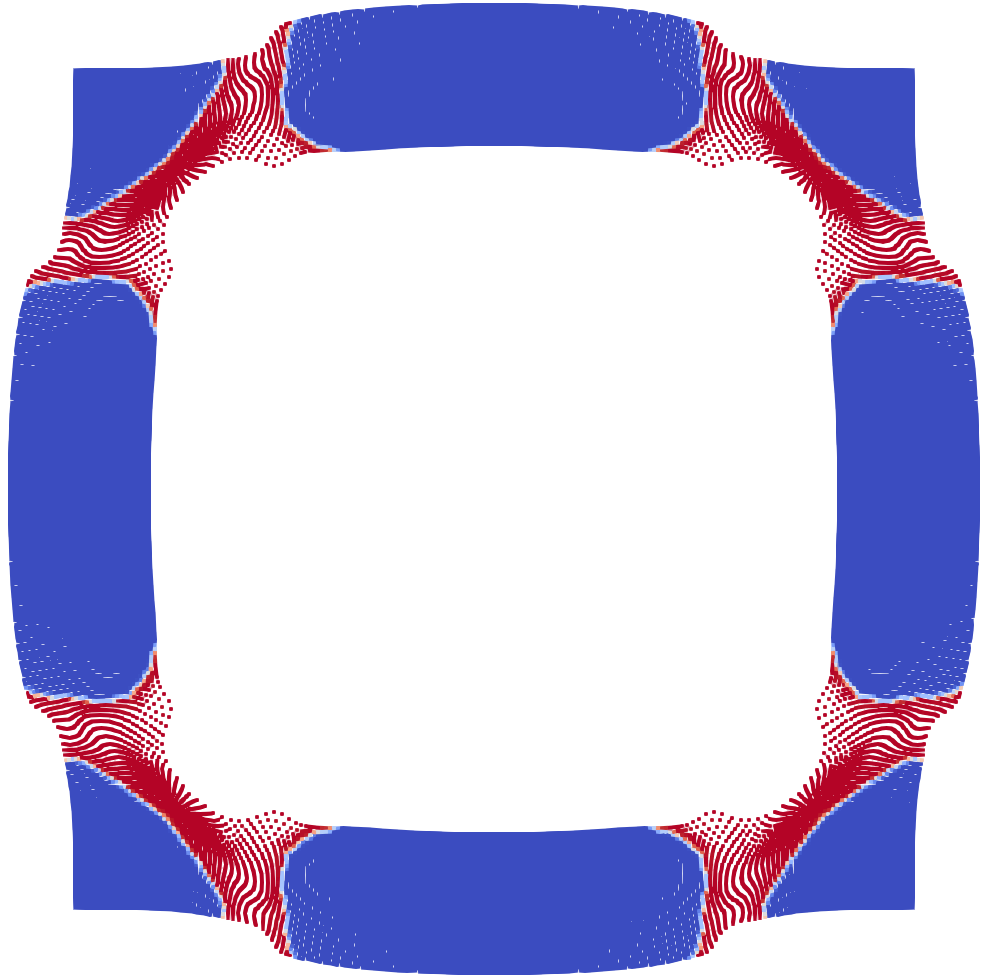}}
  \caption{Ductile fracture problem. Damage contours for different discretization levels D1--D4 at $t = \SI{225}{\micro s}$.}
  \label{fig:ductile_fracture}
\end{figure*}

\cref{fig:ductile_velocity} shows the simulation results of D4 discretization, where the air speed is shown on the background mesh while the PD damage field is plotted on the foreground PD nodes. 
The stress singularity created by the re-entrant corner of the hollow square leads to crack nucleation in that location, as expected. The cracks initially propagate along the diagonal, then broaden and branch. Eventually, the solid structure is completely fractured and split into multiple fragments with the air occupying the newly open empty space. Note the permanent deformation of the fragments occurring due to the plasticity effects. As shown in \cref{fig:ductile_fracture}, the fragmented solid structure deformed configuration converges toward a unique shape, which is hard to achieve for ill-posed local fracture models due to mesh-dependency issues.

\section{Discussion and conclusions}
\label{sec:conclusions}

In this paper, a computational FSI framework for simulating air blast events based on an immersed IGA-PD approach is presented. B-Splines are used to construct the trial and test functions in the background domain for the coupled FSI problem while the foreground solid is discretized using correspondence-based PD, which is a meshfree methodology. The PD nodes are constrained to follow the background-grid kinematics, while the solid internal forces are computed directly using the PD discretization and then transferred to the background grid to complete the coupled discrete FSI formulation. This novel approach to immersed FSI has several advantages over the existing methods. In particular, using PD eliminates the need to explicitly track the crack interfaces and makes modeling fragmentation a relatively easy task. Likewise, we do not require to track the fluid-structure interface due to the use of an immersed method. Also, the use of higher-order and smooth B-Spline functions in the background and RKPM functions for the PD solid in the foreground results is higher-order accurate solutions for both the fluid and solid fields.

The results of the proposed IGA-PD formulation match well with the solutions using other conforming and immersed methods, which provides good verification of the present methodology. Although no experimental comparison is presented for the air-blast–structure interaction problems involving dynamic fracture, the IGA–PD simulations provide physically reasonable results that also show convergence under mesh refinement.

We also note some challenges using the current approach. While we obtain the desired strong coupling by constraining the solid to the background kinematics, simulating fragmentation becomes complicated since the smooth background discretion is not designed to appropriately capture the material discontinuities in the foreground solid solution. In practice, the damage zones, whose size scales with the background element size (see \cref{fig:brittle_fracture,fig:ductile_fracture}), appear to be artificially thick leading to excessive structural damage. Refining the background mesh helps reduce the size of the damage zones, however, it considerably increases the computational costs. Relaxing the dependence of the foreground kinematics on the background solution, by leveraging such techniques as developed in~\cite{wang2021consistent}, may be a way forward to address this challenge. 

\begin{figure*}[!hbpt]
  \centering
  \subfloat{\includegraphics[width=0.24\textwidth,trim={0cm 0cm 0cm 0cm},clip]{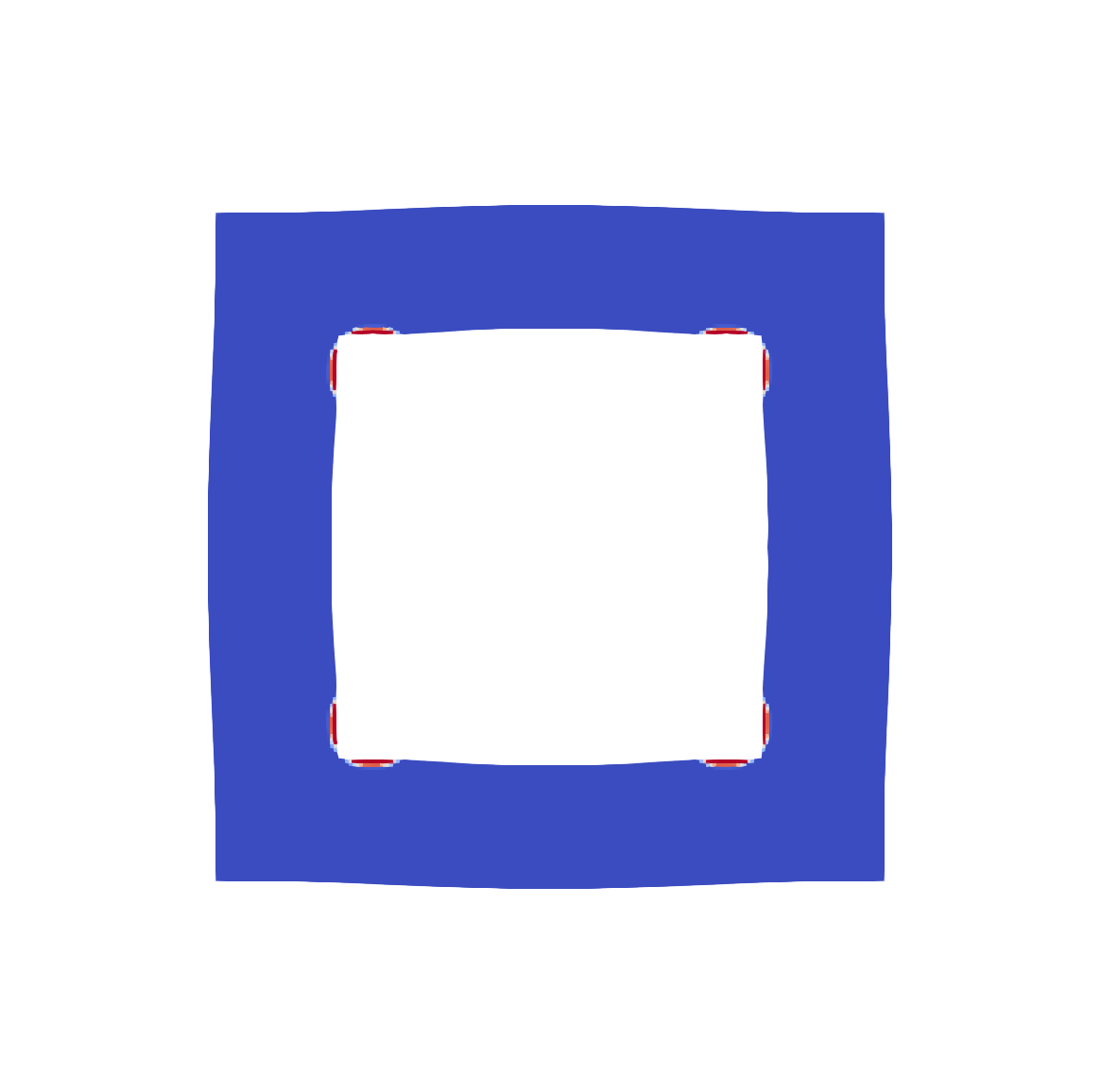}}
  \hspace{1pt}
  \subfloat{\includegraphics[width=0.24\textwidth,trim={0cm 0cm 0cm 0cm},clip]{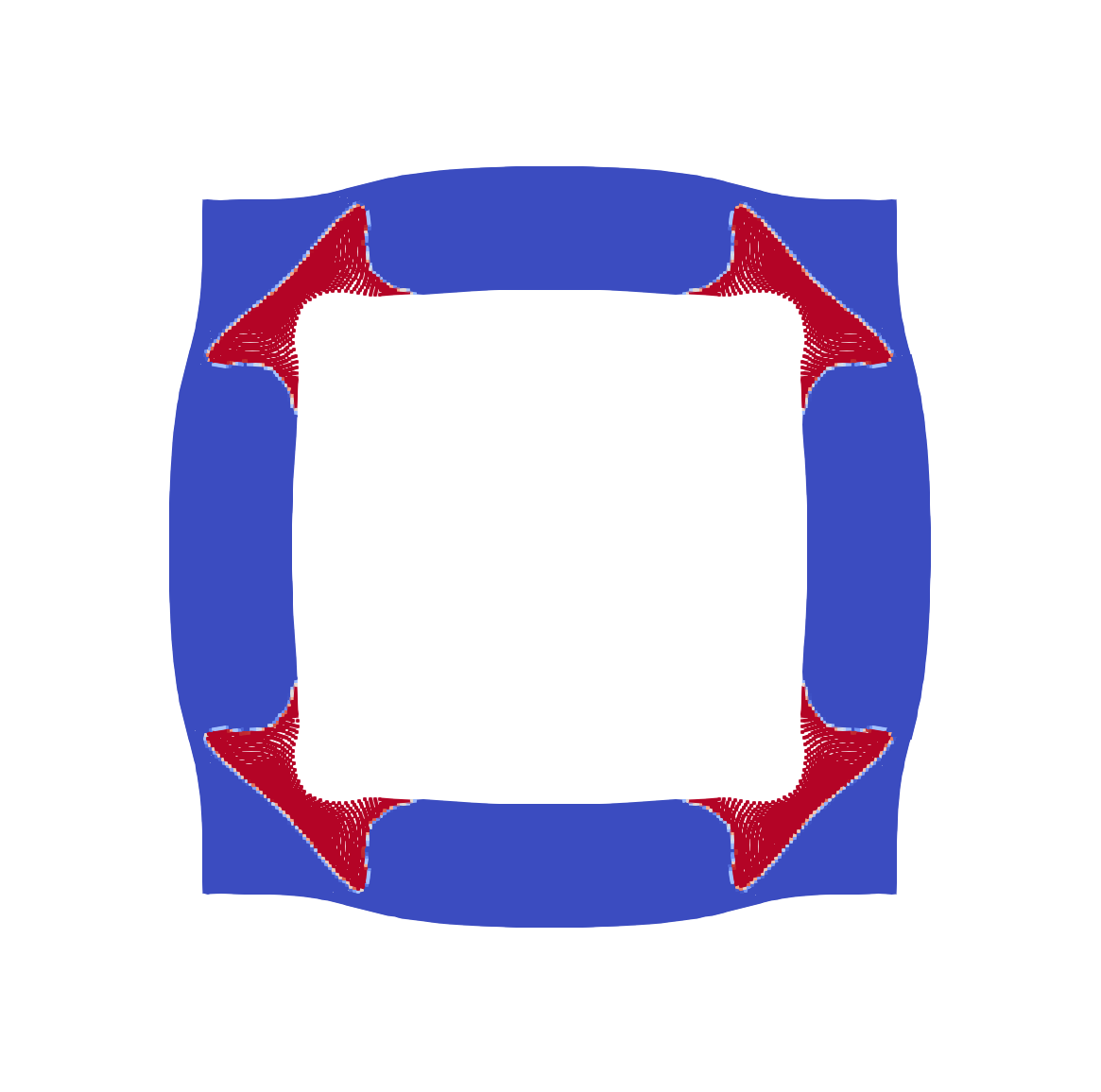}}
  \hspace{1pt}
  \subfloat{\includegraphics[width=0.24\textwidth,trim={0cm 0cm 0cm 0cm},clip]{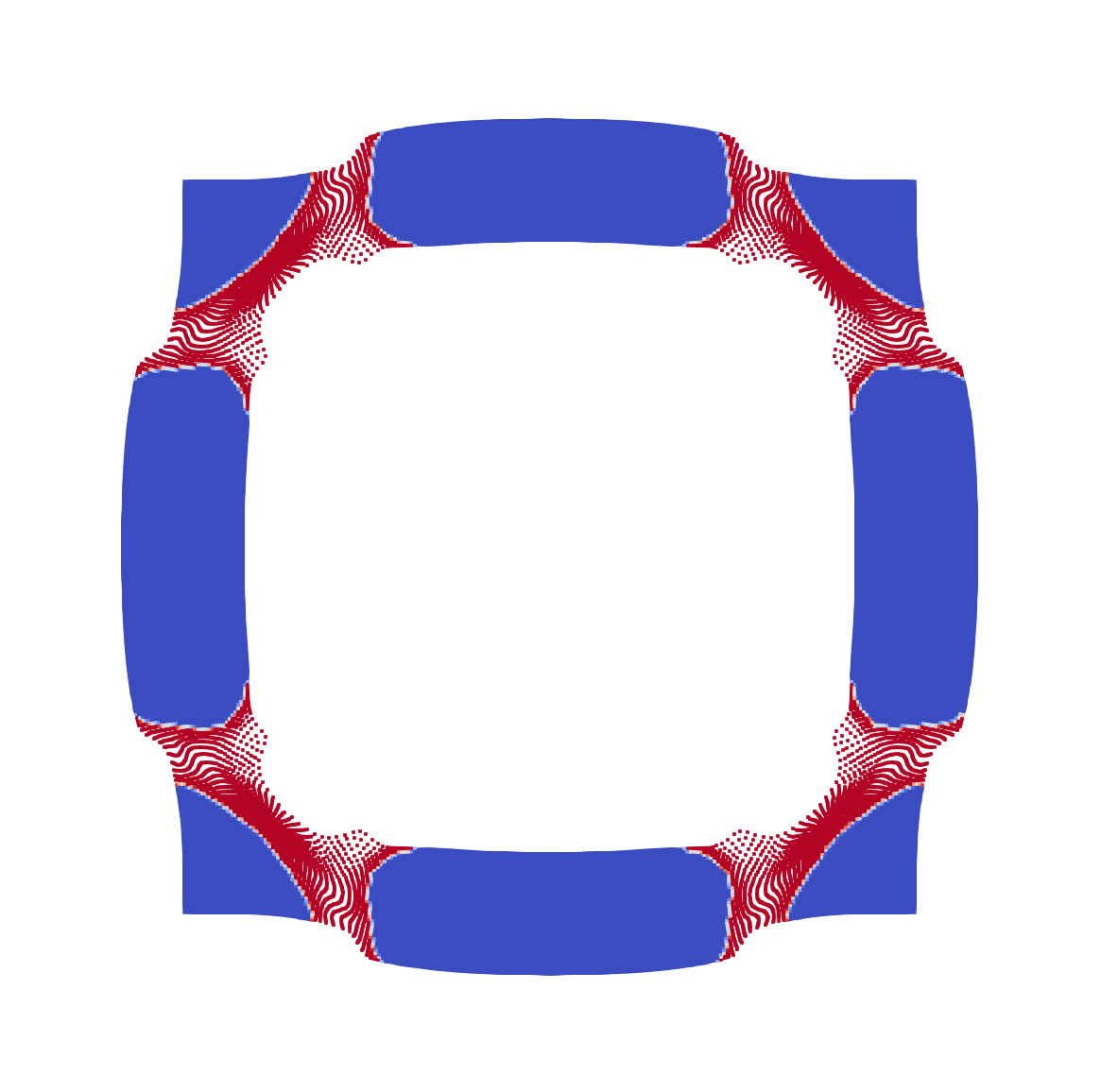}}
  \hspace{1pt}
  \subfloat{\includegraphics[width=0.24\textwidth,trim={0cm 0cm 0cm 0cm},clip]{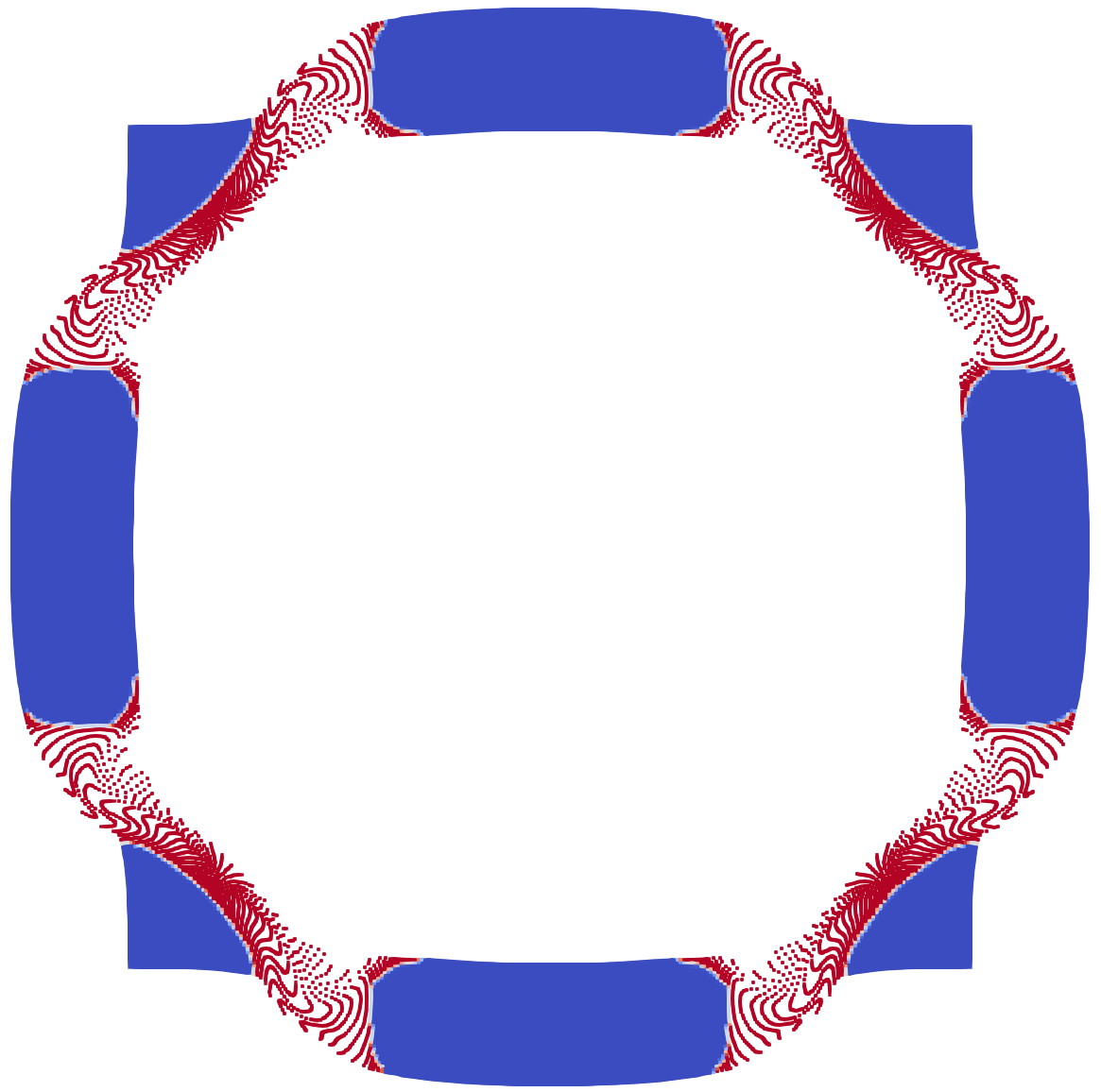}}

  \setcounter{subfigure}{0}
  \subfloat[][\SI{75}{\micro s}]{\includegraphics[width=0.24\textwidth,trim={0cm 0cm 0cm 0cm},clip]{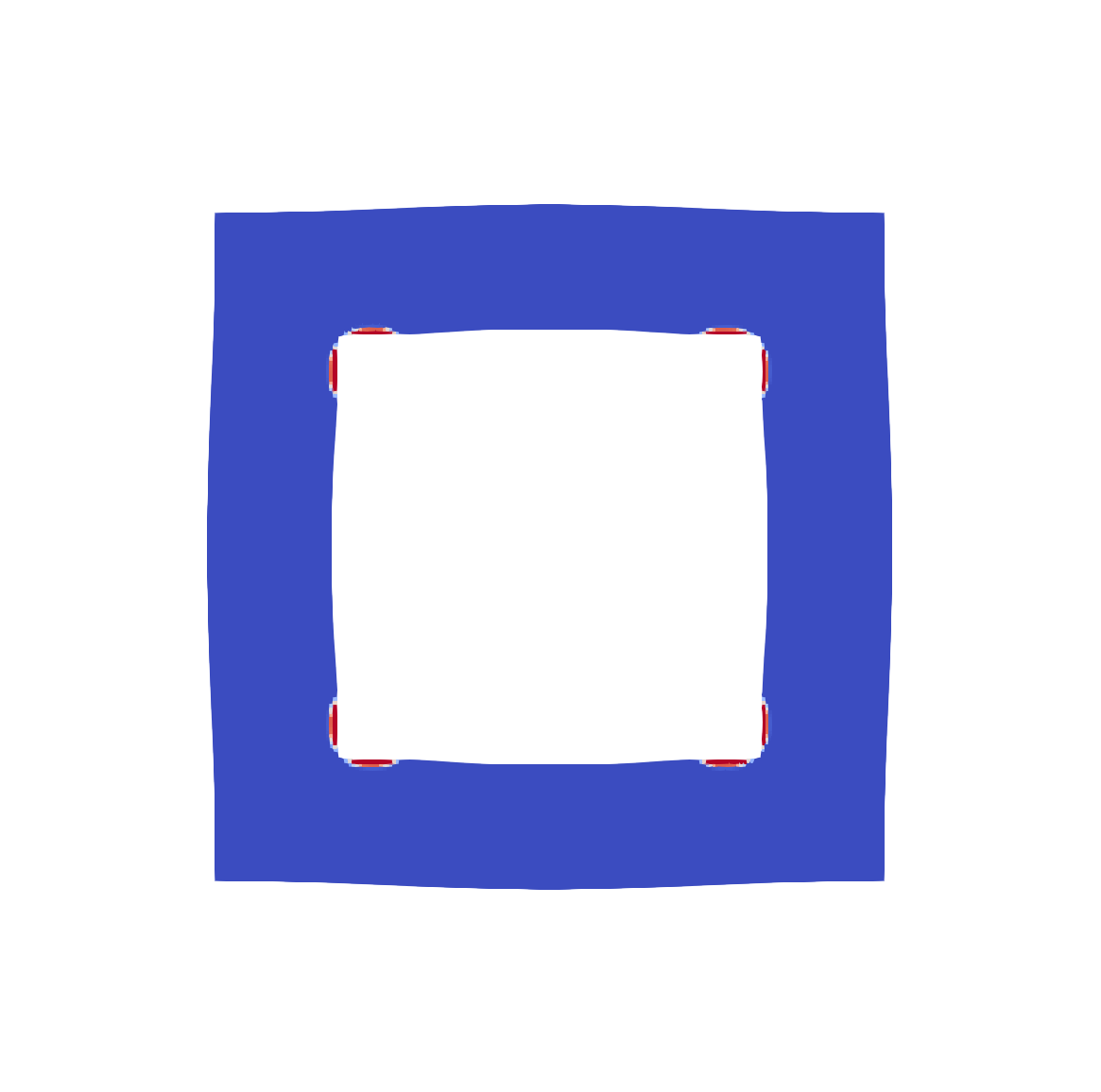}}
  \hspace{1pt}
  \subfloat[][\SI{150}{\micro s}]{\includegraphics[width=0.24\textwidth,trim={0cm 0cm 0cm 0cm},clip]{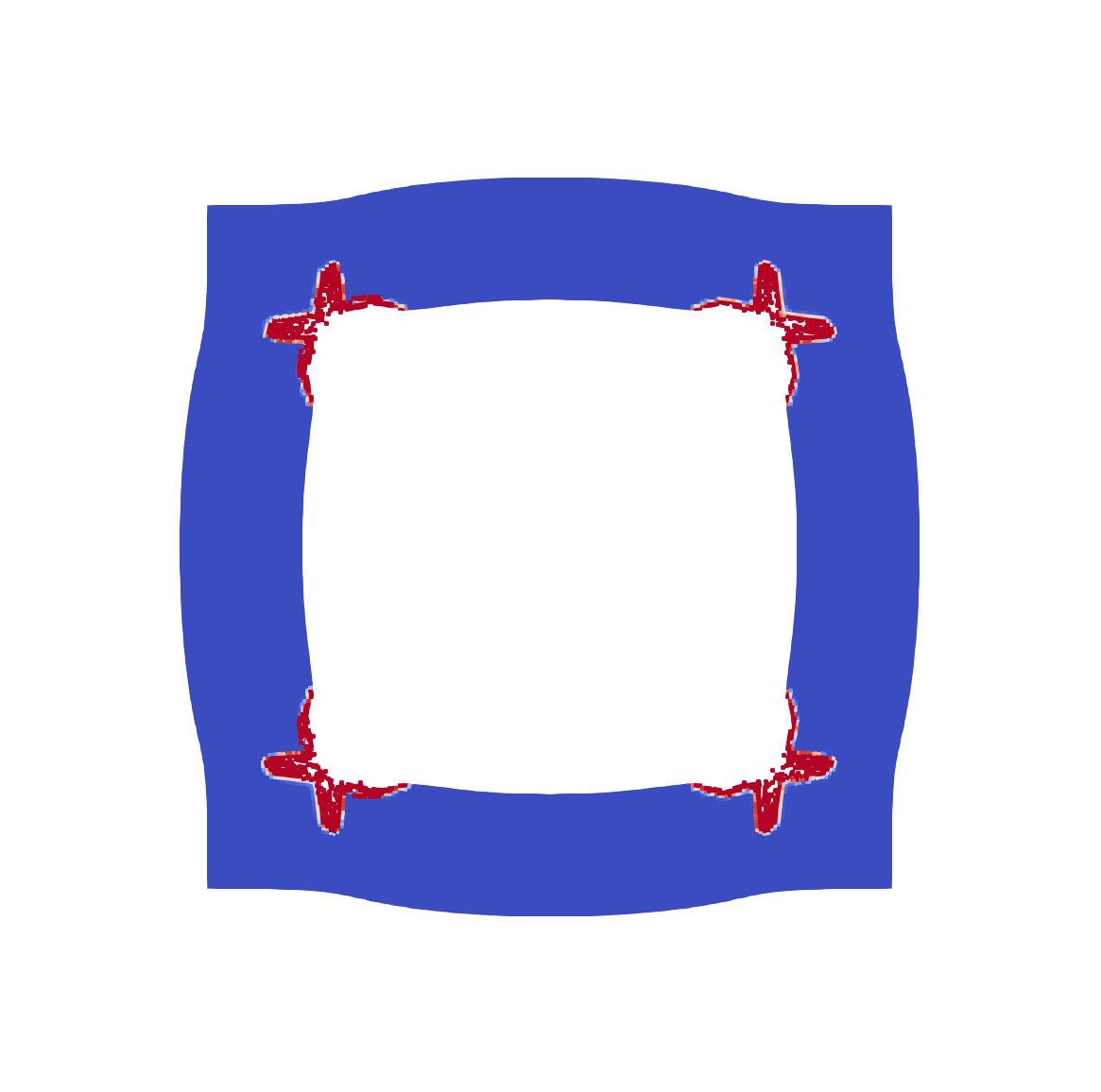}}
  \hspace{1pt}
  \subfloat[][\SI{225}{\micro s}]{\includegraphics[width=0.24\textwidth,trim={0cm 0cm 0cm 0cm},clip]{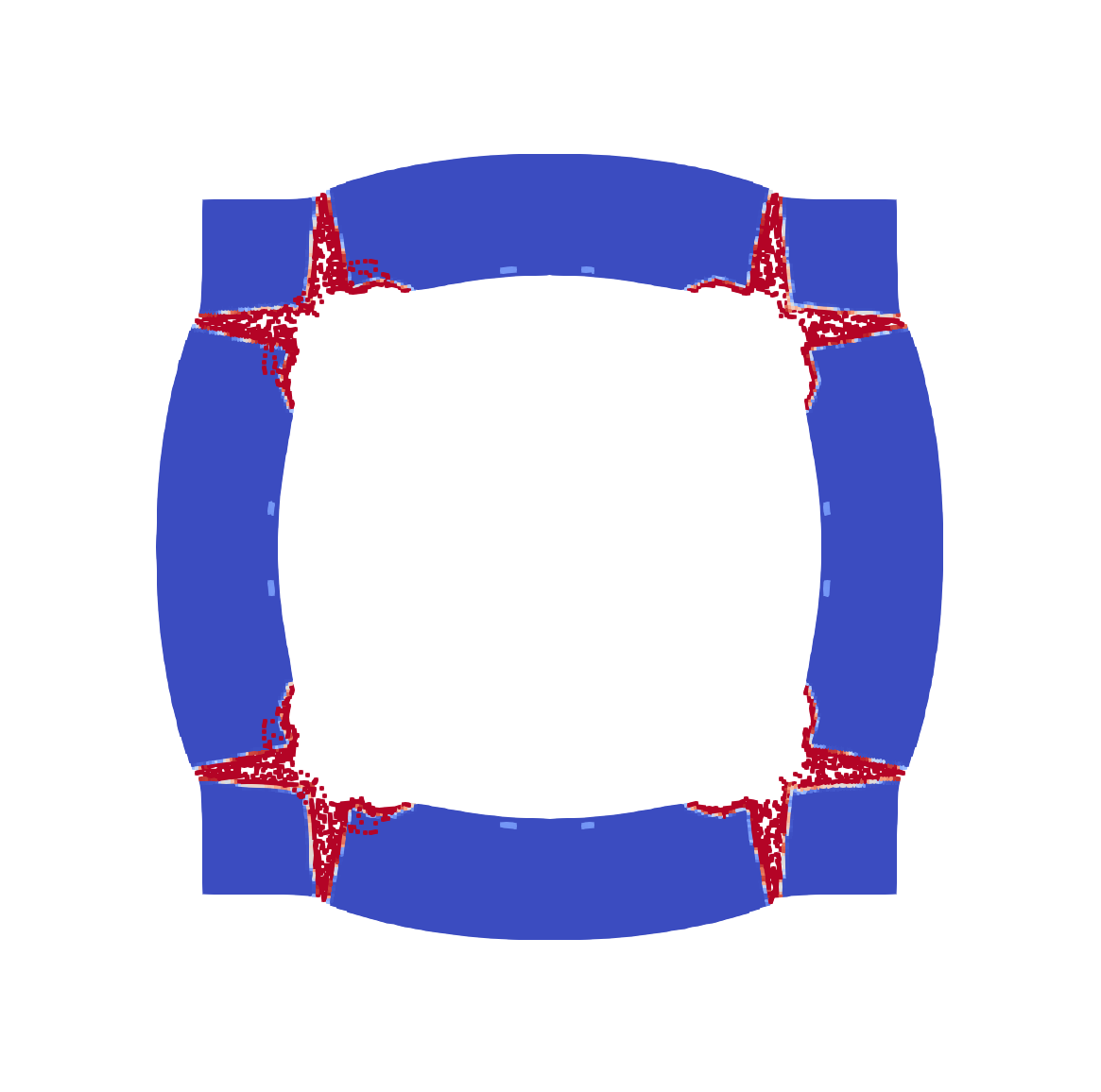}}
  \hspace{1pt}
  \subfloat[][\SI{375}{\micro s}]{\includegraphics[width=0.24\textwidth,trim={0cm 0cm 0cm 0cm},clip]{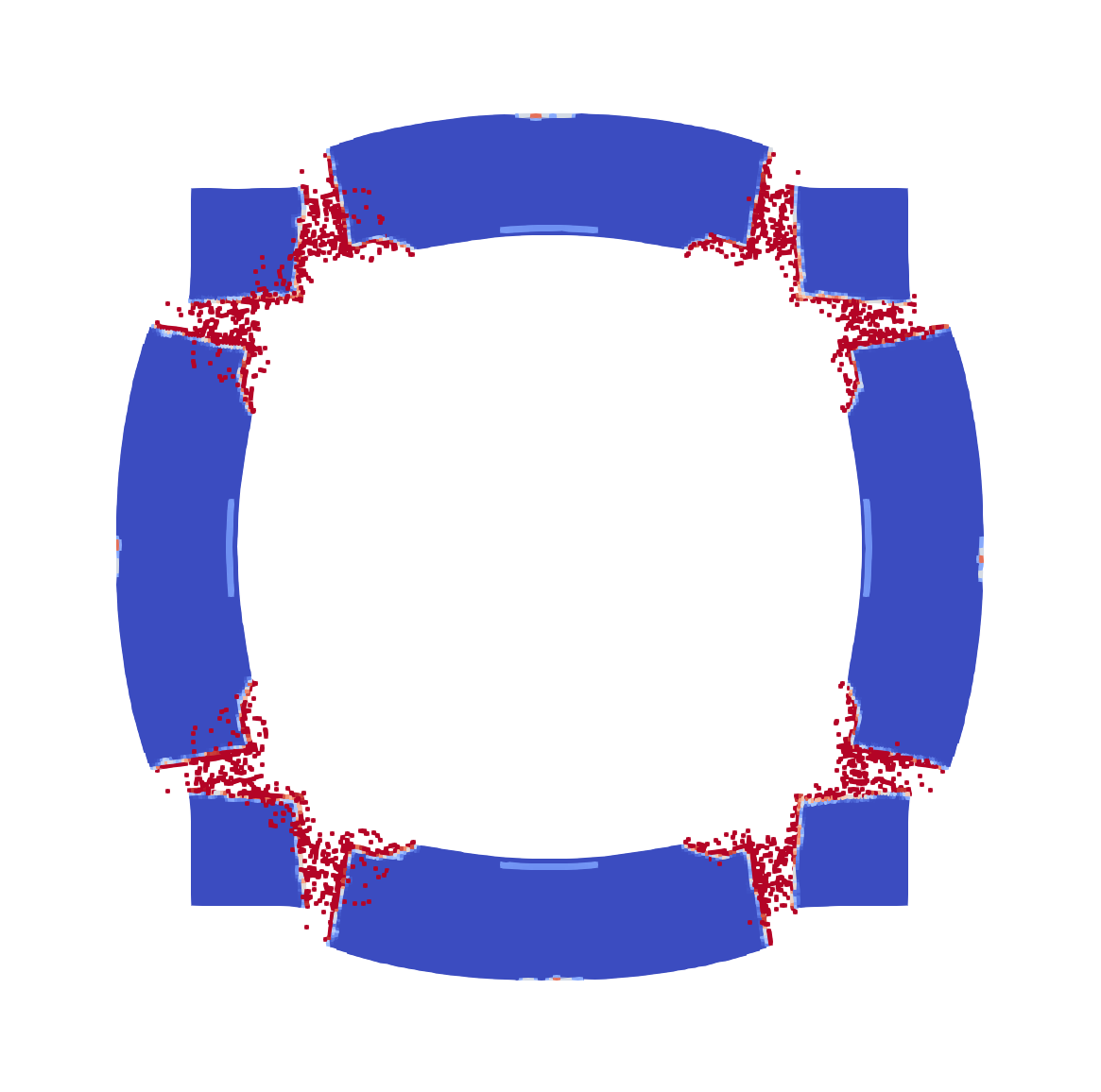}}

  \caption{Damage evolution in the ductile fracture problem. Top row corresponds to the fully coupled approach. Bottom row refers to the case where the dependence of the solid structure on the background discretization is removed at $t = \SI{75}{\micro s}$, thus letting the solid evolve on its own.}
  \label{fig:ductile_decoupling}
\end{figure*}

While we do not address this issue rigorously in the present work and leave it for the future developments, we explore a more practical way to simulate blast-induced fracture and fragmentation using the current formulation. In the blast events, the key factor that governs structural failure is the total impulse that is transmitted to the solid as a result of the explosion. Once enough loading is applied to the structure, it can nucleate cracks and potentially lead to complete disintegration. As a possible approach, we consider decoupling the foreground solution from the background discretization \textit{after} damage starts growing in the solid material and letting the PD solid evolve on its own afterwards. We explore this idea using the last numerical example (ductile structure under blast loading). As shown in \cref{fig:ductile_decoupling}, removing the constraint between the foreground and background fields once enough impulse is delivered to the structure and letting the standalone PD formulation govern the solid problem from that point on results in sharper cracks and thinner damage zones (bottom row in the figure) as opposed to the thicker damage regions in the fully coupled case (top row in the figure). While this approach involves some challenges in identifying when to decouple the solid, we observe that the unconstrained foreground solution is able to predict sharper fractures and better quality fragmentation. Future work will address this issue in a more rigorous way.

\section*{Acknowledgments}
\label{sec:acknowledge}
Y.~Bazilevs and M.~Behzadinasab were partially supported through the ONR Grant No. N00014-21-1-2670. Y.~Bazilevs was partially supported through the Sandia Contract No. 2111577. N.~Trask acknowledges funding under the DOE ASCR PhILMS center (Grant number DE-SC001924) and the Laboratory Directed Research and Development program at Sandia National Laboratories. J.T.~Foster acknowledges funding under Sandia National Laboratories contract No. 1885207. Sandia National Laboratories is a multi-mission laboratory managed and operated by National Technology and Engineering Solutions of Sandia, LLC., a wholly owned subsidiary of Honeywell International, Inc., for the U.S. Department of Energy’s National Nuclear Security Administration under contract DE-NA0003525. SAND Number: SAND2021-10191 O.

\appendix

\section{Semi-Lagrangian correspondence PD formulation of an inelastic solid}
\label{sec:semiPD}

 In the semi-Lagrangian PD model~\cite{behzadinasab2020semi}, the velocity gradient tensor is computed using only the current field variables and without a need for mapping to the reference configuration. As a prerequisite step for calculation of $\underline{\mathbf{L}}$, in this model, the point-level velocity gradient $\mathbf{L}$ is evaluated first using
\begin{equation} 
\mathbf{L} = \int_{\mathcal{H}} \underline{\mathbf{v}} \, \underline{\boldsymbol{\Phi}}^\intercal \, {\rm d}\mathcal{H} .
\label{eqn:L_node}
\end{equation}
Here, $\underline{\boldsymbol{\Phi}}$ is a set of gradient kernel functions determined using a higher-order meshfree method that is asymptotically compatible with the local gradient operator~\cite{behzadinasab2021unifiedI}. We employ the reproducing kernel (RK) shape functions~\cite{chi2013gradient,madenci2016peridynamic,hillman2020generalized}, as follows, to construct a meshfree differential operator equipped with an $n$-th order accuracy:
\begin{equation}
\underline{\boldsymbol{\Phi}}^{[n]} = \underline{\omega} \, \tilde{\boldsymbol{\Xi}}^{[n]} \, \mathbf{M}^{[n]^{-1}} \, \underline{\boldsymbol{\Xi}}^{[n]} ,
\label{eqn:Phi}
\end{equation}
where the so-called {\em influence state} $\underline{\omega}$ is a weighting function dependent on the relative distance between material points with respect to the support size. We use a cubic B-spline kernel to define the radial influence function, i.e.,
\begin{align}
\omega_r(\hat{\underline{\eta}}) = 
\begin{cases}
\dfrac{2}{3} - 4\,\hat{\underline{\eta}}^2 + 4\hat{\underline{\eta}}^3  ~~~ & \text{for} ~~~ 0 < \hat{\underline{\eta}} \leq \dfrac{1}{2} \\
\dfrac{4}{3} - 4\,\hat{\underline{\eta}} + 4\,\hat{\underline{\eta}}^2 - \dfrac{4}{3}\hat{\underline{\eta}}^3  ~~~ & \text{for} ~~~ \dfrac{1}{2} < \hat{\underline{\eta}} \leq 1 \\
0  ~~~ & \text{otherwise} 
\end{cases}
, 
\label{eqn:omega_r}
\end{align}
where $\hat{\underline{\eta}}$ is defined as the current length of bond normalized with the support size, i.e.,
\begin{equation}
\hat{\underline{\eta}} = \frac{\left \vert \underline{\mathbf{x}} \right \vert}{\delta} , 
\end{equation}
where $\underline{\mathbf{x}} = \mathbf{x}' - \mathbf{x}'$ is called the PD {\em position state}. To satisfy the identity condition in \cref{eqn:identity}, the normalized weighting function $\underline{\alpha}$ is defined as
\begin{equation}
\underline{\alpha} = \frac{ \ \underline{\omega} \ }{ \ \bar{\omega} \ } ,
\end{equation}
where
\begin{equation}
\bar{\omega} = \int_{\mathcal{H}} \underline{\omega} \, {\rm d}\mathcal{H} . 
\end{equation}

Back to \cref{eqn:Phi}, $\underline{\boldsymbol{\Xi}}^{[n]} = \boldsymbol{\Xi}^{[n]}\left(\underline{\mathbf{x}}\right)$ is a column vector of the complete set of monomials $\left\{\mathbf{x}^{\, \beta}\right\}_{|\beta|=1}^{n}$. The moment matrix $\mathbf{M}^{[n]}$ is defined as
\begin{equation}
\mathbf{M}^{[n]} = \int_{\mathcal{H}} \underline{\omega} \, \underline{\boldsymbol{\Xi}}^{[n]} \, \underline{\boldsymbol{\Xi}}^{[n]^{\intercal}} \, {\rm d}\mathcal{H} .
\end{equation}
For a $d$-dimensional Cartesian space with a $d_{\Xi}$-dimensional set $\boldsymbol{\Xi}^{[n]}$, the $d \times d_{\Xi}$ matrix $\tilde{\boldsymbol{\Xi}}^{[n]}$ is defined as
\begin{equation}
\tilde{\boldsymbol{\Xi}}^{[n]} = 
\begin{bmatrix} 
1, \ 0, \ 0, \ 0, \ 0 , \ 0, \ \dots \ , \ 0 \\
0, \ 1, \ 0, \ 0, \ 0 , \ 0, \ \dots \ , \ 0 
\end{bmatrix}.
\end{equation}

\begin{example}
\normalfont
For a 2-dimensional system and a cubic kernel function ($n=3)$, the involved terms are:
\begin{equation}
\boldsymbol{\Xi}^{[3]}(\boldsymbol\eta) = \left[ \eta_1, \, \eta_2, \, \eta_1^2, \, \eta_1 \, \eta_2, \, \eta_2^2, \, \eta_1^3, \, \eta_1^2 \, \eta_2, \, \eta_1 \, \eta_2^2, \, \eta_2^3 \right]^\intercal ,
\end{equation}
and
\begin{equation}
\tilde{\boldsymbol{\Xi}}^{[3]} = 
\begin{bmatrix} 
1, \ 0, \ 0, \ 0, \ 0 , \ 0, \ 0 , \ 0 , \ 0 \\
0, \ 1, \ 0, \ 0, \ 0 , \ 0, \ 0 , \ 0 , \ 0 
\end{bmatrix} .
\end{equation}
\end{example}

To obtain a naturally-stabilized PD correspondence model, the {\em bond-associative} (BA) stabilization method of~\cite{breitzman2018bond,behzadinasab2020semi} is utilized to compute the velocity gradient at the bond level:
\begin{equation} 
\underline{\mathbf{L}} = \frac{\mathbf{L} + \mathbf{L}'}{2} + \left( \underline{\mathbf{v}} - \frac{\mathbf{L} + \mathbf{L}'}{2}  \cdot \underline{\mathbf{x}} \right) \frac{\underline{\mathbf{x}}^\intercal}{\left|\underline{\mathbf{x}}\right|^2} ,
\label{eqn:bondL}
\end{equation}
where $\mathbf{L}' = \mathbf{L}(\mathbf{x})$. $\underline{\mathbf{L}}$ is shown to be a second-order operator if $\mathbf{L}$ is equipped with a second or higher order of accuracy~\cite{behzadinasab2021unifiedI}. At the {\bf bond level}, the velocity gradient is used to update the Cauchy stress using the classical theory, i.e.,
\begin{equation} 
\dot{\psi} \left( \mathbf{\nabla} \underline{\mathbf{L}} \right) = \underline{\boldsymbol\sigma} : \underline{\mathbf{L}}
\label{eqn:sigma-L}
\end{equation}
where $\underline{\boldsymbol\sigma}$ is the power conjugate of $\underline{\mathbf{L}}$, and a classical rate-based constitutive law governs the strain-power density function $\dot{\psi}$. The PD force state $\mathbf{T} \an{\mathbf{x}-\mathbf{x}'}$ in this model reads~\cite{behzadinasab2020semi}
\begin{equation} 
\underline{\mathbf{T}} = \frac{\underline{\omega} \, \underline{\boldsymbol\sigma} \, \underline{\mathbf{x}}}{\bar{\omega} \, \left|\underline{\mathbf{x}}\right|^2} + \mathbf{z} \, \underline{\boldsymbol{\Phi}} ,
\label{eqn:PD-T}
\end{equation}
where $\mathbf{z}(\mathbf{x})$ is given by
\begin{equation} 
\mathbf{z} = \left[ \int_{\mathcal{H}} \left(\frac{0.5}{\bar{\omega}} + \frac{0.5}{\bar{\omega}'}\right) \, \underline{\omega} \, \underline{\boldsymbol\sigma} \left( \mathbf{I} - \frac{\underline{\mathbf{x}} \, \underline{\mathbf{x}}^\intercal}{\left|\underline{\mathbf{x}}\right|^2} \right) \, {\rm d}\mathcal{H} \right]
\end{equation}
with the identity tensor $\mathbf{I}$.

\section{Solid constitutive equations}
\label{sec:solid_constitutive}
In this work, the rate-form constitutive model of~\cite{bazilevs2017new1} is adopted, which is based on the standard J2 flow theory with isotropic hardening~\cite{simo2006computational} and thus suitable for modeling metal plasticity. Note that the constitutive relation is applied at the {\bf bond level}. In this constitutive approach, the time history of the bond-associated velocity gradient derives the evolution of the bond-associated Cauchy stress. The rate-of-deformation tensor is additively decomposed into elastic and plastic parts, i.e.,
\begin{equation} 
\underline{\mathbf{D}} = \frac{1}{2} (\underline{\mathbf{L}} + \underline{\mathbf{L}}^\intercal) = \underline{\mathbf{D}}^e + \underline{\mathbf{D}}^p .
\end{equation}
To maintain objectivity, the Jaumann rate of the Cauchy stress~\cite{belytschko2013nonlinear} is used, i.e., 
\begin{equation} 
\overset{\nabla}{\underline{\boldsymbol\sigma}} = \mathbb{C} : \underline{\mathbf{D}}^e ,
\end{equation}
where $\mathbb{C}$ is the elasticity tensor. The material time derivative of the Cauchy stress $\dot{\underline{\boldsymbol\sigma}}$ is
\begin{equation} 
\dot{\underline{\boldsymbol\sigma}} = \overset{\nabla}{\underline{\boldsymbol\sigma}} + \underline{\boldsymbol\sigma} \, \underline{\boldsymbol\Omega}^\intercal + \underline{\boldsymbol\Omega} \, \underline{\boldsymbol\sigma},
\end{equation}
where the spin tensor
\begin{equation} 
\underline{\boldsymbol\Omega} = \frac{1}{2} (\underline{\mathbf{L}} - \underline{\mathbf{L}}^\intercal) .
\end{equation}
The yield function is given by
\begin{equation} 
f(\underline{\boldsymbol\sigma}, \bar{\underline{\epsilon}}^p) = \bar{\sigma} (\underline{\boldsymbol\sigma}) - \sigma_Y (\bar{\underline{\epsilon}}^p) \leq 0, 
\end{equation}
where the equivalent (von Mises) stress $\bar{\underline{\sigma}}$ is defined as
\begin{equation} 
\bar{\sigma} (\underline{\boldsymbol\sigma}) = \sqrt{\frac{3}{2} \underline{\boldsymbol\sigma}^\prime : \underline{\boldsymbol\sigma}^\prime},
\end{equation}
in which the deviatoric stress $\underline{\boldsymbol\sigma}^\prime$ is given by
\begin{equation} 
\underline{\boldsymbol\sigma}^\prime = \underline{\boldsymbol\sigma} - \frac{1}{3} (\text{tr}~\underline{\boldsymbol\sigma})~\mathbf{I} .
\end{equation}
$\sigma_Y$ is the yield stress that depends on the bond-associated equivalent plastic strain $\bar{\epsilon}^p$, which evolves according to the associative plasticity rule~\cite{simo2006computational}, i.e.,
\begin{equation} 
\dot{\bar{\underline{\epsilon}}}^p = \frac{ \sqrt{\dfrac{3}{2}} \dfrac{\underline{\boldsymbol\sigma}^\prime}{ \|\underline{\boldsymbol\sigma}^\prime \|} : \mathbb{C} : \underline{\mathbf{D}} }{ H + \dfrac{3}{2 } \dfrac{\underline{\boldsymbol\sigma}^\prime}{ \|\underline{\boldsymbol\sigma}^\prime \|} : \mathbb{C} : \dfrac{\underline{\boldsymbol\sigma}^\prime}{ \|\underline{\boldsymbol\sigma}^\prime \|} },
\end{equation}
where $H$ is the hardening modulus.

\section{Bond-associative damage correspondence modeling}
\label{sec:failure}
The bond-associative failure model~\cite{behzadinasab2020semi,behzadinasab2020revisiting} is used here. In this correspondence-based approach, a continuum damage theory is utilized to compute the state of material damage in the solid material. The influence state, using this technique, is modified as 
\begin{equation}
\underline{\omega} = \omega_r \big( |\underline{\mathbf{x}}| \big) \ \omega_D \big( \, \underline{D} \big) ,
\label{eqn:damagedOmega}
\end{equation}
where $\underline{\omega}_r$ is the conventional (radial) undamaged influence function as defined in \cref{eqn:omega_r}. $\underline{\omega}_D$ is the damage-dependent component of the influence state. A classical continuum damage model is incorporated to evolve the bond-associated damage parameter $\underline{D} \in [0,1]$ based on the bond-associated internal variables such as strain, stress, stretch, and temperature. It is required that $\underline{\omega}_D(0) = 1$ for an undamaged bond ($\underline{D}=0$) and $\underline{\omega}_D(1) = 0$ for a fully damaged bond ($\underline{D}=1$). Damage irreversibility is enforced by constraining $\underline{\omega}_D$ as a non-increasing function of $\underline{D}$. 

In this work, to model the brittle fracture phenomenon, we break a bond once its associated equivalent (von Mises) stress exceeds a critical value. When the failure criterion is met, the bond breakage is enforced sharply to respect the immediate crack propagation nature of brittle materials.

For simulating the ductile fracture process, the gradual aspect of material degradation in malleable materials is considered. A plasticity-driven failure approach is adopted to decay the load-carrying capacity of a bond in its transition from the undamaged to fully-damaged state, i.e.,
\begin{equation}
\begin{aligned}
& \omega_D \left( \underline{D} \right) = 1 - \underline{D} \ , \\
& \underline{D} = D \left( \bar{\underline{\epsilon}}^P \right) =
\begin{cases}
0 , \qquad \qquad \bar{\underline{\epsilon}}^P <= \bar{\epsilon}^P_{\rm th} , \\
\dfrac{\bar{\underline{\epsilon}}^P - \bar{\epsilon}^P_{\rm th}}{\bar{\epsilon}^P_{\rm cr} - \bar{\epsilon}^P_{\rm th}} , \ \quad \bar{\epsilon}^P_{\rm th} < \bar{\underline{\epsilon}}^P < \bar{\epsilon}^P_{\rm cr} , \\
1 , \qquad \qquad \bar{\underline{\epsilon}}^P >= \bar{\epsilon}^P_{\rm cr} ,
\end{cases}
\end{aligned}
\label{eqn:plastic_strain}
\end{equation}
where $\bar{\epsilon}^P_{\rm th}$ and $\bar{\epsilon}^P_{\rm cr}$ are the model parameters corresponding to the undamaged and completely damaged states, respectively.

\bibliographystyle{unsrt}
\bibliography{main}

\end{document}